%% file: RCG_FurtherResults.tex
\documentclass[]{article}

\usepackage{xr}
\input{RCG_preamble_biblio.tex}

\input{RCG_preamble_main.tex}

\numberwithin{equation}{section}

\title{\sffamily Further Results and Discussions on Random Cayley Graphs}

\author{\sffamily Jonathan Hermon\quad Sam Olesker-Taylor}%

\date{}

\usepackage{geometry}
\begin{document}

\maketitle

\acknofootnote

\vspace{-6ex}

\renewcommand{\abstractname}{\sffamily Abstract}
\begin{abstract}
Consider the random Cayley graph of a finite group $G$ with respect to $k$ generators chosen uniformly at random, with $1 \ll k \lesssim \log |G|$. The results of this article supplement those in the three main papers on random Cayley graphs, namely \cite{HOt:rcg:abe:cutoff,HOt:rcg:matrix,HOt:rcg:abe:geom}. The majority of the results are inspired by a `universality' conjecture of Aldous and Diaconis~\cite{AD:conjecture}.

To start, we study the limit profile of cutoff for the simple random walk on this random graph, as well as a detailed investigation into mixing properties when $G = \mathbb Z_p^d$ with $p$ prime. We then exposit a proof of Diaconis and Saloff-Coste~\cite{DSc:growth-rw} establishing lack of cutoff when $k \asymp 1$. We move onto discussing material from \cite{HOt:rcg:matrix} on matrix groups. We then study distance of a typical element of $G$ from the identity in an $L_q$-type graph distance in the Abelian set-up. Finally, we give necessary and sufficient conditions for $k$ independent uniform elements of $G$ to generate $G$, ie for the random Cayley graph to be connected, based on work of Pomerance~\cite{P:gen-abelian}.

The aforementioned results all hold with high probability over the random Cayley graph.
\end{abstract}

%
%
%

\small
\begin{quote}
\begin{description}
	\item [Keywords:]
	cutoff, mixing times, random walk, random Cayley graphs, universality, entropy
	
	\item [MSC 2020 subject classifications:]
	05C12, 05C80, 05C81; 20D15; 60J27, 60K37
\end{description}
\end{quote}
\normalsize


\vspace{2ex}
\numberingroman

\vfill
\printtoc{\value{tocdepth}}
\vspace*{2ex}

\newpage
\section{Introduction and Statement of Results}
\label{sec-p5:intro}

\subsection{Motivation, Brief Overview of Results and Notation}

\subsubsection{Motivating Conjecture of \citeauthor{AD:conjecture}}

We analyse properties of the random walk (abbreviated \textit{RW}) on a \textit{Cayley graph} of a finite group.
The generators of this graph are chosen independently and uniformly at random.
Precise definitions are given in \S\ref{sec-p2:intro:cayley-def}; for now, let $G$ be a finite group, let $k$ be an integer (allowed to depend on $G$) and denote by $G_k$ the Cayley graph of $G$ with respect to $k$ independently and uniformly random generators.
We consider values of $k$ with $1 \ll \log k \ll \log \abs G$ for which $G_k$ is connected with high probability (abbreviated \textit{whp}), ie with probability tending to 1 as $\abs G$ grows.


\medskip

Since pioneering work of Erd\H{o}s, it has been  understood that the typical behaviour of \emph{random} objects in some class can shed valuable light on the class as a whole.
Thus, when considering some class of combinatorial objects, it is natural to ask questions such as the following.
\begin{itemize}[noitemsep, topsep = \smallskipamount, label = \bcdot] \slshape
	\item What does a typical object in this class `look like'?
	\item If an object is chosen uniformly at random, which properties hold with high probability?
\end{itemize}
\textcite{AD:conjecture} applied this philosophy to the study of random walks on groups.

\smallskip

\textcite{AD:conjecture,AD:shuff-stop} coin the phrase \textit{cutoff phenomenon}:
	this occurs when the total variation distance (TV) between the law of the RW and its invariant distribution drops abruptly from close to $1$ to close to $0$ in a time-interval of smaller order than the mixing time.
The material in this article is motivated by a conjecture of theirs regarding `universality of cutoff' for the RW on the random Cayley graph $G_k$; it is given in \cite[Page~40]{AD:conjecture}, which is an extended version of \cite{AD:shuff-stop}.

\begin{introconj*}[\citeauthor{AD:conjecture}, \citeyear{AD:conjecture}]
	For any group $G$,
	if $k \gg \log \abs G$ and $\log k \ll \log \abs G$, then
	the random walk on $G_k$ exhibits cutoff whp.
	Further, the cutoff time, to leading order, is independent of the algebraic structure of the group: it can be written as a function only of~$k$~and~$\abs G$.
\end{introconj*}

This conjecture spawned a large body of work, including \cite{D:phd,DH:enumeration-rws,H:cutoff-cayley->,H:cutoff-cayley-survey,H:cutoff-cayley-<,R:random-random-walks,W:rws-hypercube}; see \S\ref{sec-p2:intro:previous-work}.
It has been established in the Abelian set-up
by \citeauthor{DH:enumeration-rws} \cite{DH:enumeration-rws,H:cutoff-cayley->}; see \S\ref{sec-p2:intro:previous-work:ad-conj} and \S\ref{sec-p2:conc-rmks:roichman} where we give a short proof.
Save \cite{H:cutoff-cayley-<} which considers the cyclic group $\mbz_p$ for prime $p$ and \cite{W:rws-hypercube} which considers $\mbz_2^d$ (which enforces $k \gtrsim \log \abs G$),
focus has been on $k \gg \log \abs G$.

We extend consideration to $1 \ll k \lesssim \log \abs G$, establishing cutoff and the limit profile for a large class of Abelian groups. We study other statistics, of a geometric flavour, in this regime too.

\subsubsection{Brief Overview of Results}

Our focus is on mixing (for the RW) and geometric properties of the random Cayley graph.
We consider the limit as $n \cq \abs G \to \infty$ under the assumption that $1 \ll \log k \ll \log \abs G$.
The condition $1 \ll \log k \ll \log \abs G$ is necessary for cutoff on $G_k^\pm$ for all nilpotent $G$; see \cref{rmk-p5:intro:cutoff:nec-con}.

\begin{itemize}[itemsep = 0pt, topsep = \smallskipamount, label = \bcdot]
	\item 
	\textit{Limit Profile.}
	For $k \gtrsim \log \abs G$,
	under some conditions on the group,
	we find the limit profile of the convergence to equilibrium of the RW.
	(We cover the case $1 \ll k \ll \log \abs G$ in \cite[\S\ref{sec-p2:cutoff1}]{HOt:rcg:abe:cutoff}.)
	
	\item 
	\textit{Detailed Investigation of $\mbz_\pp^d$.}
	For $1 \ll k \lesssim \log \abs G$ and $G \cq \mbz_\pp^d$ with $\pp$ prime,
	we analyse some cases not covered by the previous cutoff analysis of \cite{HOt:rcg:abe:cutoff}. Eg, we allow $k$ to be very close to $d$.
	
	\item 
	\textit{$L_\qq$-Type Typical Distances.}
	For $1 \ll k \lesssim \log \abs G$, we study an $L_\qq$-type graph distance between a typical point $U \sim \Unif(G)$ and the identity in the random Cayley graph.
	
	\item 
	\textit{Additional Discussions.}
	We explain briefly
		how to adapt our cutoff arguments for the upper triangular group studied in \cite[\S\ref{sec-p1:cutoff}]{HOt:rcg:matrix} can be extended to other, more general, nilpotent groups
	and
		how there is no cutoff when $k \asymp 1$ (which is a result due to \textcite{DSc:growth-rw}).
	For Abelian groups, we give (simple) necessary and sufficient conditions on $k$ in terms of $G$ for the group to be generated by independent, uniform generators (based on \textcite{P:gen-abelian}).
\end{itemize}

Introduced by \textcite{AD:conjecture}, there has been a great deal of research into these random Cayley graphs.
Motivation for this model and an overview of historical work is referenced~in~\S\ref{sec-p5:intro:previous-work}.

\subsubsection{Notation and Terminology}

Cayley graphs are either directed or undirected; we emphasise this by writing $G_k^+$ and $G_k^-$, respectively.
When we write $G_k$ or $G^\pm_k$, this means ``either $G^-_k$ or $G^+_k$'', corresponding to the undirected, respectively directed, graphs with generators chosen independently and uniformly at~random.

Conditional on being simple, $G^+_k$ is uniformly distributed over the set of all simple degree-$k$ Cayley graphs. Up to a slightly adjusted definition of \textit{simple} for undirected Cayley graphs, our results hold with $G_k$ replaced by a uniformly chosen simple Cayley graph of degree $k$; see \S\ref{sec-p5:intro:rmks:typ-simp}.

Our results are for sequences $(G_N)_\Ninn$ of finite groups with \toinf{\abs{G_N}} \asinf N.
For ease of presentation, we write statements like ``let $G$ be a group'' instead of ``let $(G_N)_\Ninn$ be a sequence of groups''.
Likewise, the quantities $d(G)$, $p$ and, of course, $k$ appearing in the statements all correspond to sequences, which need not be fixed (or bounded) unless we explicitly say otherwise.
In the same vein, an event holds \textit{with high probability} (abbreviated \textit{\whp}) if its probability tends~to~1.

We use standard asymptotic notation:
	``$\ll$'' or ``$\oh{\cdot}$'' means ``of smaller order'';
	``$\lesssim$'' or $\Oh{\cdot}$'' means ``of order at most'';
	``$\asymp$'' means ``of the same order'';
	``$\eqsim$'' means ``asymptotically equivalent''.

\subsection{Statements of Main Results}
\label{sec-p5:intro:res}

We analyse mixing in the \textit{total variation} (abbreviated \textit{TV}) distance.
The uniform distribution on $G$, denoted $\pi_G$, is invariant for the RW.
Let $S = (S(t))_{t\ge0}$ denote the RW on $G_k$; its law is denoted $\pr[G_k]{S(t) \in \cdot}$.
For $t \ge 0$,
denote the TV distance between the law of $S(t)$ and $\pi_G$ by
\[
	d_{G_k}(t)
\cq
	\tvb{ \pr[G_k]{ S(t) \in \cdot } - \pi_G }
=
	\MAX{A \subseteq G}
	\absb{ \pr[G_k]{ S(t) \in A } - \abs A / \abs G }.
\]
Throughout, unless explicitly specified otherwise, we use continuous time: $t \ge 0$ means $t \in [0,\infty)$.

\subsubsection{Cutoff: Limit Profile for Random Walks on Abelian Groups}
\label{sec-p5:intro:res:gen}

\nextresult

We use standard notation and definitions for \textit{mixing} and \textit{cutoff}; see, eg, \cite[\S 4 and \S18]{LPW:markov-mixing}.

\begin{introdefn*}
	A sequence $(X^N)_\Ninn$ of Markov chains is said to exhibit \textit{cutoff} when, in a short time-interval, known as the \textit{cutoff window}, the TV distance of the distribution of the chain from equilibrium drops from close to $1$ to close to $0$, or more precisely if there exists $(t_N)_\Ninn$ with
	\[
		\LIMINF{\Ninf} d_N\rbb{ t_N (1 - \eps) } = 1
	\Qand
		\LIMSUP{\Ninf} d_N\rbb{ t_N (1 + \eps) } = 0
	\Qforall
		\eps \in (0,1),
	\]
	where $d_N(\cdot)$ is the TV distance of $X^N(\cdot)$ from its equilibrium distribution for each $\Ninn$.
	
	We say that a RW on a sequence of random graphs $(H_N)_\Ninn$ \textit{exhibits cutoff around time $(t_N)_\Ninn$ whp} if,
	for all fixed $\eps$, in the limit $\Ninf$,
	the TV distance at time $(1 + \eps) t_N$ converges in distribution to $0$ and at time $(1 - \eps) t_N$ to $1$,
	where the randomness is over $H_N$.
\end{introdefn*}

We use an \emph{entropic method}, which involves defining \emph{entropic times}; see \cite[\S\ref{sec-p2:intro:previous-work:generic-ent}]{HOt:rcg:abe:cutoff} for a high-level description of the method and \S\ref{sec-p5:profile:ent} for a more specific description.
The main idea is to use an auxiliary process $W$ to generate the walk $S$; one then studies the entropy of the process $W$.
	Here, $W_i(t)$ is, for each $i$, the number of times generator $Z_i$ has been applied minus the number of times $Z_i^{-1}$ has been applied;
	$W$ is a rate-1 RW on $\mbz^k$.
	Then $S(t) = W(t) \bcdot Z$ when the group is Abelian.

For undirected graphs, $W$ is the usual simple RW (abbreviated \textit{SRW}):
	a coordinate is selected uniformly at random and incremented/decremented by 1 each with probability $\tfrac12$.
For directed graphs, inverses are never applied, so a step of $W$ is as follows:
	a coordinate is selected uniformly at random and incremented by 1;
	we term this the \textit{directed} RW (abbreviated~\textit{DRW}).

For $t \ge 0$, write $\mu_t$ for the law of $W(t)$.
Define $Q(t) \cq - \log \mu_t\rbr{ W(t) }$.

\begin{introdefn}
\label{def-p5:intro:ent-time}
	For all $k,n \in \mbn$ and all $\alpha \in \mbr$,
	define
	\(
		t_\alpha \cq t_\alpha(k,n)
	\)
	so that
	\[
		\ex{ Q\rbr{t_\alpha} } = \rbb{ \log n + \alpha \sqrt{v k} }
	\Qwhere
		v \cq \Varb{Q\rbr{t_0}}/k.
	\]
	We call $t_0$ the \textit{entropic time} and the $\bra{t_\alpha}_{\alpha\in\mbr}$ \textit{cutoff times}.
	
	An asymptotic evaluation of these times is given in \cref{res-p5:profile:t0a}.
\end{introdefn}

For an Abelian group $G$, write $d(G)$ for the minimal size of a generating subset of $G$ and
\[
	m_*(G)
\cq
	\max\brb{ \mint{j \in [d]} m_j \midb \oplus_1^d \: \mbz_{m_j} \text{ is a decomposition of } G }.
\]
Also write $\Psi : \mbr \to (0,1)$ for the standard Gaussian tail:
\(
	\Psi(\alpha) \cq (2\pi)^{-1/2} \intt[\infty]{\alpha} e^{-x^2/2} dx
\)
for $\alpha \in \mbr$.

Recall that $1 \ll \log k \ll \log \abs G$ is necessary for cutoff for nilpotent $G$, eg Abelian $G$;
see \cref{rmk-p5:intro:cutoff:nec-con}.
A more refined statement is given in \cref{res-p5:profile:res}. 

\begin{introthm}
\label{res-p5:intro:cutoff}
	Let $G$ be an Abelian group.
	Suppose that
		$k \gtrsim \log \abs G$,
		$\log k \ll \log \abs G$,
		$d(G) \le \tfrac15 \log \abs G / \log \log k$
	and
		$m_*(G) \ge \log \abs G / \log \log \abs G$.
	Abbreviate $t_\alpha \cq t_\alpha(k, \abs G)$~for~all~\text{$\alpha \in \mbr$}.
	
	Then, \whp, the RW on $G_k$ exhibits cutoff at $t_0$ with Gaussian profile given by $\bra{t_\alpha}_{\alpha \in \mbr}$:
	\begin{gather*}
		d_{G_k}(t_\alpha) \to^\mbp \Psi(\alpha)
	\Quad{(in probability)}
		\text{for all}
	\quad
		\alpha \in \mbr;
	\\
		t_0 \asymp k
	\Qand
		\abs{ t_\alpha - t_0 } \asymp \alpha t_0 / \sqrt k \ll t_0
	\Quad{for all}
		\alpha \in \mbr.
	\end{gather*}
\end{introthm}

\begin{subtheorem-num}{intrormkT}
	\label{rmk-p5:intro:cutoff}

\begin{intrormkt}
	We can write the cutoff statement in terms of the mixing time, rather than the TV distance:
	writing $\tmix(\eps)$ for the $\eps$-mixing time,
	for all $\eps \in (0,1)$,
	we have
	\[
		\rbb{ \tmix(\eps) - t_0 } / w
	\to^\mbp
		\Psi^{-1}(\eps),
	\]
	where $t_0$ is the mixing time and $w$ is the cutoff window defined via $\bra{t_\alpha - t_0}_{\alpha \in \mbr}$.
	For a more explicit formula, using asymptotic evaluation of the cutoff times, see \cref{rmk-p5:profile:mixing-profile}.
\end{intrormkt}

\begin{intrormkt}
\label{rmk-p5:intro:cutoff:nec-con}
This article establishes cutoff in a variety of set-ups, but always in the regime $1 \ll \log k \ll \log \abs G$.
This leaves the regimes $k \asymp 1$ and $\log k \asymp \log \abs G$, for which there is no cutoff for any choice of generators:
	when $k \asymp 1$, this holds whenever the group is nilpotent;
	when $\log k \asymp \log \abs G$, this holds for all groups.
The former result is due to \textcite{DSc:growth-rw}; we give a short exposition of this in \S\ref{sec-p5:const-k} below.
The latter result is proved in \cite[\S\ref{sec-p2:conc-rmks:roichman}]{HOt:rcg:abe:cutoff}; the mixing time~is~order~1.
\textcite[Theorems~3.3.1 and~3.4.7]{D:phd} establishes a more general result for $\log k \asymp \log \abs G$.
\end{intrormkt}

\end{subtheorem-num}

\subsubsection{Cutoff: A Detailed Investigation of $\mathbb Z_p^d$}
\label{sec-p5:intro:res:p}

\nextresult

For our next theorem, we specialise to the case $G \cq \mbz_\pp^d$ with $\pp$ prime.
This specialisation allows us to derive some very refined results.
In particular, before we could not allow $d$ to be very close to $k$; here we consider any $k \ge d$.
Now every element of $G$ has order $\pp$; as such, need only consider the auxiliary $W$ mod $\pp$.
We redefine the entropic times to take this into account.

\begin{introdefn}
	Define $t_0 \cq t_0(k,\pp,d)$ to be the time at which the entropy of the RW on $\mbz_\pp^k$ is~$\log(\pp^d)$.
	
	An asymptotic evaluation of this time is given in Proposition~\ref{res-p5:p:t0a}.
\end{introdefn}

A more refined statement is given in \cref{res-p5:p:res}.
In particular, by defining $t_\alpha$ appropriately, we can also also consider the cutoff window; see \cref{def-p5:p:entropic-time,res-p5:p:t0a}.

\begin{introthm}
\label{res-p5:intro:p}
	Let $G \cq \mbz_\pp^d$ with $\pp$ prime.
	Assume that $1 \ll k \lesssim \log \abs G = d \log \pp$.
	Then, whp, the RW on $G_k$ exhibits cutoff at time $t_0$ conditional on generating the group, ie $\langle Z_1, ..., Z_k \rangle = G$.
	Further, if $(k - d) \pp \gg 1$ then the group is generated whp.
\end{introthm}

%
%

\begin{intrormkt}
	The RW on the usual degree-$d$ Cayley graph of $G \cq \mbz_\pp^d$ exhibits cutoff at the same time $\tfrac12 d \log d / (1 - \cos(2 \pi / \pp))$ as for the random degree-$k$ Cayley graph above with $k - d \asymp 1$.
	
	This is not overly surprising:
		one can check that, \emph{for $p$ prime}, if $\langle Z_1, ..., Z_k \rangle = G$ then there exists an $I \subseteq [k]$ with $\abs I = d$ such that $\langle Z_i \rangle_{i \in I} = G$ and the Cayley graph $G([Z_i]_{i \in I})$ is isomorphic to the usual degree-$d$ Cayley graph.
	Thus, conditional on generating the group, the random Cayley graph is simply the usual one with an additional $k - d \asymp 1$ uniformly chosen generators.
\end{intrormkt}

\subsubsection{Cutoff: No Cutoff when $k$ Is Constant}
\label{sec-p5:intro:res:const-k}

\nextresult

It is natural to ask what happens when $k$ is a fixed constant.
This regime has already been analysed by \textcite{DSc:growth-rw}.
We give an exposition of their results, using the language which we have developed.
We emphasise that this is a result of \citeauthor{DSc:growth-rw}.

A more refined statement is given in \cref{res-p5:const-k:main}.

\begin{introthm}[{cf \cite[Corollary~5.3]{DSc:growth-rw}}]
	Let $G$ be a finite, nilpotent group of bounded step.
	Suppose that $k \asymp 1$.
	Then the RW on $G^-(Z)$ does not exhibit cutoff for any choice of $Z$ with $\abs Z = k$.
\end{introthm}

%
%
%

\subsubsection{Typical Distance: Generalised Graph Distance}
\label{sec-p5:intro:res:typdist}

\nextresult

Our final result regards graph distances between the identity and a uniformly chosen element of $G$.
For a Cayley graph $\mcg$, for $R \ge 0$ and $\beta \in [0,1]$, define the \textit{$\beta$-typical distance} $\mcd_\mcg(\beta)$ via
\[
	\mcb_\mcg(R)
\cq
	\brb{ x \in \mcg \midb \dist(\id, x) \le R }
\Qand
	\mcd_\mcg(\beta)
\cq
	\min\brb{ R \ge 0 \midb \abs{\mcb_\mcg(R)} \ge \beta \abs \mcg }.
\]

Locally, when $\log k \ll \log \abs G$, typical degree-$k$ Cayley graphs of an Abelian group look like $\mbz^k$.
In a lattice, graph distance corresponds to $L_1$ distance; this can be extended to an $L_\qq$ distance, for $\qq \in [1,\infty]$.
Analogously, we can extend the usual $L_1$ graph distances to an $L_\qq$-type, for $\qq \in [1,\infty]$.

Consider a collection $z = [z_1, ..., z_k]$ of generators and distances in the Cayley graph $G(z)$.
For a path $\rho$ in $G(z)$,
for each $i \in [k]$,
write
	$\rho_{i,+}$ for the number of times $z_i$ is used,
	$\rho_{i,-}$ for the number of times $z_i^{-1}$ is used (if in the undirected case, otherwise $\rho_{i,-} \cq 0$)
and
	$\rho_i \cq \rho_{i,+} - \rho_{i,-}$.
The path connects the identity with $\rho \bcdot z$.
Then the ($L_1$) graph distance
of $\rho$ is $\norm{\rho}_1 \cq \sumt[k]{1} \rbr{ \rho_{i,+} + \rho_{i,-} }$.

For any $\qq \in [1,\infty)$,
define the \textit{$L_\qq$ graph distance} of $\rho$ by
\(
	\norm{\rho}_\qq^\qq
\cq
	\sumt[k]{1} \rbr{ \rho_{i,+}^\qq + \rho_{i,-}^\qq }.
\)
For the \textit{$L_\infty$ graph distance},
define
\(
	\norm{\rho}_\infty
\cq
	\maxt{i} \bra{ \rho_{i,+} + \rho_{i,-}}.
\)
(The usual graph distance is given by $\qq = 1$.)

For Abelian groups, clearly for any $\qq \in [1,\infty)$ an \textit{$L_\qq$ geodesic}, ie a path of minimal $L_\qq$ weight, will only use either $z_i$ or $z_i^{-1}$, not both (since the terms in the product can be reordered), ie $\rho_{i,+} \rho_{i,-} = 0$ for all $i$.
Thus $\norm{\rho}_\qq^\qq = \sumt[k]{1} \abs{\rho_i}^\qq$.
Similarly, any $L_\infty$ geodesic $\rho$ can be adjusted into a new path $\rho'$ with $\rho \bcdot z = \rho' \bcdot z$ and $\norm{\rho}_\infty = \norm{\rho'}_\infty$ satisfying $\rho'_{i,+} \rho'_{i,-} = 0$ for all $i$.

\smallskip

We define the \textit{$L_\qq$ typical distance} $\mcd_{G(z),\qq}(\cdot)$ analogously to $\mcd_{G(z)}(\cdot)$, ie the $\qq = 1$ case.

\smallskip

For an Abelian group $G$, we write $d(G)$ for the minimal size of a generating subset of $G$ and
\[
	m_*(G)
\cq
	\max\brb{ \mint{j \in [d]} m_j \midb \oplus_1^d \: \mbz_{m_j} \text{ is a decomposition of } G }.
\]
Finally we set up a little more notation.
Make the following definitions for $q \in [1,\infty]$:
\[
	C^-_\qq \cq 2 \, \Gamma(1/\qq+1) (\qq e)^{1/\qq},
\quad
	C_\qq^+ \cq \tfrac12 C^-_\qq
\Qand
	\mfd^\pm_\qq(k,n) \cq k^{1/\qq} n^{1/k} / C^\pm_\qq,
\]
where the case $\qq = \infty$ is to be interpreted as the limit $\qq \to \infty$;
eg,
\(
	C^-_\infty
=
	2
\)
and
\(
	\mfd^+_\infty(k,n)
=
	n^{1/k}.
\)

A more refined statement is given in \cref{res-p5:typdist:res}.
Use the convention $k^{1/\infty} \cq 1$.

\begin{introthm}
	Let $G$ be an Abelian group and $q \in [1,\infty]$.
	Suppose that $1 \ll k \ll \log \abs G$.
	Suppose that $m_*(G) \gg k^{1/q} \abs G^{1/k}$ and if $q \in (1,\infty)$ then additionally require $k \le \log \abs G / \log \log \abs G$.
	Suppose that $\limsup d(G)/k < 1$ for undirected graphs and $\limsup d(G)/k < \tfrac12$ for directed graphs.
	
	Then, \whp, the $L_\qq$ typical distance on $G_k$ concentrates at $\mfd_\qq^\pm$, namely
	\[
		\mcd_{G_k,\qq}^\pm(\beta) / \mfd_\qq^\pm \to^\mbp 1
	\Quad{(in probability)}
		\text{for all $\beta \in (0,1)$}.
	\]
	(The randomness is over the uniform choice of generators $Z = [Z_1, ..., Z_k]$.)
\end{introthm}

\subsection{Historic Overview}
\label{sec-p5:intro:previous-work}

The study of these random Cayley graphs has a rich history.
Rather than repeat it here, we refer the reader to the appropriate points in the main articles.

\begin{itemize}[itemsep = 0pt, topsep = \smallskipamount, label = \bcdot]
	\item 
	For motivation, discussions and history on universality of cutoff for this model, and a conjecture of \textcite{AD:conjecture,AD:shuff-stop} which inspired this work, see
		\cite[\S\ref{sec-p2:intro:previous-work:ad-conj}]{HOt:rcg:abe:cutoff}
	and
		\cite[\S\ref{sec-p1:intro:previous-work:ad-conj}]{HOt:rcg:matrix}.
	\item 
	For discussions on using an entropic method to establish cutoff for `generic' instances of Markov chains, see \cite[\S\ref{sec-p2:intro:previous-work:generic-ent}]{HOt:rcg:abe:cutoff}.
	
	\item 
	For history of diameter-flavour results, see \cite[\S\ref{sec-p3:intro:previous-work:typdist}]{HOt:rcg:abe:geom} or \cite[\S\ref{sec-p1:intro:previous-work:typdist}]{HOt:rcg:matrix}.
\end{itemize}

\subsection{Additional Remarks}
\label{sec-p5:intro:rmks}

\subsubsection{Precise Definition of Cayley Graphs}
\label{sec-p5:intro:cayley-def}

Consider a finite group $G$.
Let $Z$ be a multisubset of $G$.
We consider geometric properties,
namely through distance metrics and the spectral gap,
of the \textit{Cayley graph} of $(G,Z)$;
we call $Z$ the \textit{generators}.
The
	\textit{undirected}, respectively \textit{directed},
\textit{Cayley graph of $G$ generated by $Z$}, denoted
	$G^-(Z)$, respectively $G^+(Z)$,
is the multigraph whose vertex set is $G$ and whose edge multiset is
\[
	\sbb{ \bra{ g,g \cdot z } \mid g \in G, \, z \in Z },
\Quad{respectively}
	\sbb{ \rbr{ g,g \cdot z } \mid g \in G, \, z \in Z }.
\]
(We do not assume that the Cayley graph is connected; that is, $Z$ may not generate $G$.)
If the walk is at $g \in G$, then a step in $G^+(Z)$, respectively $G^-(Z)$, involves choosing a generator $z \in Z$ uniformly at random and moving to $g z$, respectively one of $g z$ or $g z^{-1}$ each with probability~$\tfrac12$.

We focus attention on the \emph{random} Cayley graph defined by choosing $Z_1, ..., Z_k \sim^\iid \Unif(G)$; when this is the case, denote $G^+_k \cq G^+(Z)$ and $G^-_k \cq G^-(Z)$.
While we do not assume that the Cayley graph is connected (ie, $Z$ may not generate $G$), in the Abelian set-up the random Cayley graph $G_k$ is connected whp whenever $k - d(G) \gg 1$; see \cref{res-p5:gen:dichotomy} below.
In the nilpotent set-up, this is the case whenever $k - d(G/[G,G]) \gg 1$; see \cite[\cref{rmk-p1:intro:typdist:nil-gen}]{HOt:rcg:matrix}.


The graph depends on the choice of $Z$.
Sometimes it is convenient to emphasise this;
we use a subscript, writing $\pr[G(z)]{\cdot}$ if the graph is generated by the group $G$ and multiset~$z$.
Analogously, $\pr[G_k]{\cdot}$ stands for the \emph{random} law $\pr[G(Z)]{\cdot}$ where $Z = [Z_1, ..., Z_k]$ with $Z_1, ..., Z_k \sim^\iid \Unif(G)$.

\subsubsection{Typical and Simple Cayley Graphs}
\label{sec-p5:intro:rmks:typ-simp}

The directed Cayley graph $G^+(z)$ is simple if and only if no generator is picked twice, ie $z_i \ne z_j$ for all $i \ne j$.
The undirected Cayley graph $G^-(z)$ is simple if in addition no generator is the inverse of any other, ie $z_i \ne z_j^{-1}$ for all $i,j \in [k]$.
In particular, this means that no generator is of order 2, as any $s \in G$ of order 2 satisfies $s = s^{-1}$---this gives a multiedge between $g$ and $g s$ for each $g \in G$.

The RW on $G^-(z)$ is equivalent to an adjusted RW on $G^+(z)$ where,
when a generator $s \in z$ is chosen,
instead of applying a generator $s$,
either $s$ or $s^{-1}$ is applied, each with probability $\tfrac12$.
Abusing terminology, we relax the definition of simple Cayley graphs to allow order 2 generators, ie remove the condition $z_i \ne z_i^{-1}$ for all $i$.

Given a group $G$ and an integer $k$,
we are drawing the generators $Z_1, ..., Z_k$ independently and uniformly at random.
It is not difficult to see that the probability of drawing a given multiset depends only on the number of repetitions in that multiset.
Thus, conditional on being simple, $G_k$ is uniformly distributed on all simple degree-$k$ Cayley graphs.
Since $k \ll \sqrt{\abs G}$, the probability of simplicity tends to 1 \asinf{\abs G}.
So when we say that our results hold ``\whp (over $Z$)'', we could equivalently say that the result holds ``for almost all degree-$k$ simple Cayley graphs of $G$''.

Our asymptotic evaluation does not depend on the particular choice of $Z$, so the statistics in question depend very weakly on the particular choice of generators for almost all choices.
In many cases, the statistics depend only on $G$ via $\abs G$ and $d(G)$.
This is a strong sense of `universality'.

\subsubsection{Overview of Random Cayley Graphs Project}
\label{sec-po:intro:rmks:advert}

This paper is one part of an extensive project on random Cayley graphs.
There are
	three main articles \cite{HOt:rcg:abe:cutoff,HOt:rcg:matrix,HOt:rcg:abe:geom},
	this technical report \cite{HOt:rcg:abe:extra}
and
	a supplementary document \cite{HOt:rcg:supp}.
\textit{Each main article is readable independently.}

The main objective of the project is to establish cutoff for the random walk and determining whether this can be written in a way that, up to subleading order terms, depends only on $k$ and $\abs G$; we also study universal mixing bounds, valid for all, or large classes of, groups.
Separately, we study the distance of a uniformly chosen element from the identity, ie typical distance, and the diameter; the main objective is to show that these distances concentrate and to determine whether the value at which these distances concentrate depends only on $k$ and $\abs G$.

\begin{itemize}[noitemsep, topsep = \smallskipamount, label = \bcdot]
	\item [\cite{HOt:rcg:abe:cutoff}]
	Cutoff phenomenon (and Aldous--Diaconis conjecture) for general Abelian groups; also, for nilpotent groups, expander graphs and comparison of mixing times with Abelian groups.
	
	\item [\cite{HOt:rcg:abe:geom}]
	Typical distance, diameter and spectral gap for general Abelian groups.
	
	\item [\cite{HOt:rcg:matrix}]
	Cutoff phenomenon and typical distance for upper triangular matrix groups.
	
	\item [\cite{HOt:rcg:abe:extra}]
	Additional results on cutoff and typical distance for general Abelian groups.
	
	\item [\cite{HOt:rcg:supp}]
	Deferred technical results mainly regarding random walk on $\mbz$ and the volume of lattice balls.
\end{itemize}

\subsubsection{Acknowledgements}
\label{sec-p5:intro:ackno}

This whole random random Cayley graphs project has benefited greatly from advice, discussions and suggestions of many of our peers and colleagues.
We thank a few of them specifically here.

\begin{itemize}[itemsep = 0pt, topsep = \smallskipamount, label = $\bcdot$]
	\item 
	Allan Sly for suggesting the underlying entropy idea for cutoff (used in \S\ref{sec-p5:profile} and \S\ref{sec-p5:p}).
	
	\item 
	Justin Salez for reading this paper in detail and giving many helpful and insightful comments as well as stimulating discussions ranging across the entire random Cayley graphs project.
	
	\item 
	Itai Benjamini for discussions on typical distance.
	
	\item 
	P\'eter Varj\'u for multiple insightful discussions on mixing for the upper triangular group and for more general nilpotent groups, particularly other step-2 nilpotent groups.
	
	\item 
	Evita Nestoridi and Persi Diaconis for general discussions, consultation and advice.
\end{itemize}

\section{Cutoff: Limit Profile for Random Walks on Abelian Groups}
\label{sec-p5:profile}

In \cite[\cref{res-p2:intro:tv}]{HOt:rcg:abe:cutoff} we established cutoff for arbitrary Abelian groups when $1 \ll \log k \ll \log \abs G$ and $k - d(G) \gg 1$.
(The subregime where $k \gg \log \abs G$ was already established.)
Further, in \cite[\cref{res-p2:cutoff1:res}]{HOt:rcg:abe:cutoff} we determined \textit{cutoff profile} when some conditions on $k$ were imposed:
\begin{alignat*}{2}
	&\text{to consider any} \ k - d(G) \gg 1,&
\quad
	&\text{we needed} \ k \ll \sqrt{\log \abs G / \log\log\log \abs G};
\\
	&\text{to consider any} \ k - d(G) \asymp k,&
\quad
	&\text{we needed} \ k \ll \log \abs G / \log\log\log \abs G;
\\
	&\text{to consider any} \ k \ll \log \abs G,&
\quad
	&\text{we needed} \ d \ll \log \abs G / \log\log\log \abs G;
\end{alignat*}
see \cite[\cref{hyp-p2:cutoff1,rmk-p2:cutoff1:hyp}]{HOt:rcg:abe:cutoff}.
In particular, we could never consider general $k \gtrsim \log \abs G$.
Recall that cutoff had already been established for arbitrary Abelian groups when $k \gg \log \abs G$, but the window, never mind the profile, was not known.
In this section, we outline how to alleviate the conditions on $k$, at the cost of some conditions on the group.

\subsection{Entropic Times: Methodology, Definition and Concentration}
\label{sec-p5:profile:ent}

We use an `entropic method', as mentioned in \S\ref{sec-p5:intro:res:gen}; cf \cite{BLPS:giant-mixing,BCS:cutoff-entropic,BL:cutoff-entropic-covered,Ck:cutoff-lifts}.
The method is fairly general; we now explain the specific application in a little more depth.

\smallskip

We define an auxiliary random process $(W(t))_{t\ge0}$, recording how many times each generator has been used:
	for $t \ge 0$, for each generator $i = 1,...,k$, write $W_i(t)$ for the number of times that it has been picked by time $t$.
By independence, $W(\cdot)$ forms a rate-1 DRW on $\mbz_+^k$.
For the undirected case, recall that we either apply a generator or its inverse; when we apply the inverse of generator $i$, increment $W_i \to W_i - 1$ (rather than $W_i \to W_i + 1$).
In this case, $W(\cdot)$ is a SRW on $\mbz^k$.

Since the underlying group is Abelian, the order in which the generators are applied is irrelevant and generator-inverse pairs cancel; hence we can write
\[
	S(t)
=
	\sumt[k]{i=1} W_i(t) Z_i
=
	W(t) \bcdot Z.
\]

Recall that the invariant distribution is uniform on $G$, giving mass $1/n$ to each vertex.
The proposed mixing time is then the time at which the auxiliary process $W$ obtains entropy $\log n$.
This time can be calculated fairly precisely in many situations; see \cref{res-p5:profile:t0a}.

\medskip

We now define precisely the notion of \textit{entropic times}.
Write $\mu_t$, respectively $\nu_s$, for the law of $W(t)$, respectively $W_1(sk)$;
so $\mu_t = \nu_{t/k}^{\otimes k}$.
Define
\[
	Q_i(t) \cq - \log \nu_{t/k}\rbb{W_i(t)},
\Quad{and set}
	Q(t) \cq - \log \mu_t\rbb{W(t)} = \sumt[k]{1} Q_i(t).
\]
So $\ex{Q(t)}$ and $\ex{Q_1(t)}$ are the entropies of $W(t)$ and $W_1(t)$, respectively. Observe that $t \mapsto \ex{Q(t)} : [0,\infty) \to [0,\infty)$ is a smooth, increasing bijection.

\begin{defn}[Entropic and Times]
\label{def-p5:profile:t0a}
	For all $k,n \in \mbn$ and all $\alpha \in \mbr$,
	define
	\(
		t_\alpha \cq t_\alpha(k,n)
	\)
	so that
	\[
		\ex{ Q_1\rbr{t_\alpha} } = \rbb{ \log n + \alpha \sqrt{v k} }/k
	\Qand
		s_\alpha \cq t_\alpha/k,
	\Qwhere
		v \cq \Varb{Q_1\rbr{t_0}},
	\]
	assuming that $\log n + \alpha \sqrt{v k} \ge 0$.
	We call $t_0$ the \textit{entropic time} and the $\bra{t_\alpha}_{\alpha\in\mbr}$ \textit{cutoff times}.
\end{defn}

Direct calculation with the Poisson distribution and SRW on $\mbz$ gives the following relations.
A sketch is given below; the rigorous details are given in \cite[\S\ref{sec-p0:se}]{HOt:rcg:supp}.

\begin{prop}[Entropic and Cutoff Times, {\cite[\cref{res-p0:se:t0a}]{HOt:rcg:supp}}]
\label{res-p5:profile:t0a}
Assume that $1 \ll \log k \ll \log n$.
Write $\kappa \cq k/\log n$.
For all $\alpha \in \mbr$,
we have
\(
	t_\alpha \eqsim t_0
\)
and furthermore,
for some functions $f$ and $g$ and all $\lambda > 0$,
the following relations hold:
\begin{subequations}
	\label{eq-p5:profile:t0a}
\begin{alignat*}{3}
	\text{if}
\quad
	k &\eqmathsbox{gen:t0a}{\ll} \log n,&
\Quad{then}
	t_\alpha &\eqsim k \cdot n^{2/k} / (2 \pi e)
&\Qand
	\rbr{ t_\alpha - t_0 }/t_0 &\eqsim \alpha \sqrt 2 / \sqrt k;
\label{eq-p5:profile:t0a:k<logn}
\nt
\\
	\text{if}
\quad
	k &\eqmathsbox{gen:t0a}{\eqsim} \lambda \log n,&
\Quad{then}
	t_\alpha &\eqsim k \cdot f(\lambda)
&\Qand
	\rbr{ t_\alpha - t_0 }/t_0 &\eqsim \alpha g(\lambda) / \sqrt k;
\label{eq-p5:profile:t0a:k=logn}
\nt
\\
	\text{if}
\quad
	k &\eqmathsbox{gen:t0a}{\gg} \log n,&
\Quad{then}
	t_\alpha &\eqsim k \cdot 1/(\kappa \log \kappa)
&\Qand
	\rbr{ t_\alpha - t_0 }/t_0 &\eqsim \alpha \sqrt{\kappa \log \kappa} / \sqrt k.
\label{eq-p5:profile:t0a:k>logn}
\nt
\end{alignat*}
\end{subequations}
Moreover, $f,g : (0,\infty)$ are continuous bijections, whose value differs between SRW and DRW.
\end{prop}

\begin{Proof}[Sketch of Proof]
In \cite[\S\ref{sec-p2:cutoff1:ent:eval}]{HOt:rcg:abe:cutoff}, we sketched the argument for $k \ll \log n$.
For $k \asymp \log n$, the target entropy is order 1, and so all the random variables are bounded in probability, away from both 0 and $\infty$.
For $k \gg \log n$, we have $t_0 \ll k$, so approximate the RW by a Bernoulli distribution (with a uniformly chosen sign for the SRW); in \cite[\S\ref{sec-p2:cutoff1:ent:eval}]{HOt:rcg:abe:cutoff}, for $k \ll \log n$, we had $t_0 \gg k$ and so approximated by a normal distribution.
With this adaptation, the sketch from \cite[\S\ref{sec-p2:cutoff1:ent:eval}]{HOt:rcg:abe:cutoff} passes~over.
\end{Proof}

Since the $W_i$, and hence the $Q_i$, are iid, $Q$ is a sum of $k$ iid random variables.
Also, it turns out that
\(
	\Var{ Q(t) } \eqsim \Var{ Q(t_0) } \gg 1
\Qwhere
	t \eqsim t_0;
\)
see \cite[\cref{res-p0:se:var:t0}]{HOt:rcg:supp}.
It then stands to reason that a CLT holds for $Q = \sumt[k]{1} Q_i$; this is indeed the case.
The following CLT, which will be of great importance, is proved in \cite[\S\ref{sec-p0:se}]{HOt:rcg:supp}.

\begin{prop}[{\cite[\cref{res-p0:se:CLT}]{HOt:rcg:supp}}]
\label{res-p5:profile:CLT}
	Assume that $1 \ll \log k \ll \log n$.
	For all $\alpha \in \mbr$,
	we have
	\[
		\pr{ Q(t_\alpha) \le \log n \pm \omega } \to \Psi(\alpha)
	\Qfor
		\omega \cq \Varb{ Q(t_0) }^{1/4} = (vk)^{1/4}.
	\]
	(There is no specific reason for choosing this $\omega$. We just need some $\omega$ with $1 \ll \omega \ll (vk)^{1/2}$.)
\end{prop}

\subsection{Precise Statement and Remarks}

In this section we give the more refined version of \cref{res-p5:intro:cutoff}.
Recall that,
for an Abelian group $G$,
we write $d(G)$ for the minimal size of a generating subset of $G$ and
\[
	m_*(G)
\cq
	\max\brb{ \mint{j \in [d]} m_j \midb \oplus_1^d \: \mbz_{m_j} \text{ is a decomposition of } G }.
\]

\begin{hyp}
\label{hyp-p5:profile}
	The sequence $(k_N, G_N)_\Ninn$ satisfies \textit{\cref{hyp-p5:profile}} if
	$\abs{G_N} \to \infty$ as $\Ninf$ and
	there exists a constant $c > 0$ so that the following inequalities hold for all $\Ninn$:
	\begin{gather*}
		k_N \ge c \log \abs{G_N},
	\quad
		m_*(G_N) \ge \log k_N / \log \log k_N
	\Qand
		d(G_N) \le \tfrac15 \log \abs{G_N} / \log \log k_N.
	\end{gather*}
\end{hyp}

Recall that we write $d^\pm_{G_k, N}(t)$ for the TV distance from uniform at time $t$ for the walk on $G_k^\pm$ and $\Psi$ for the standard Gaussian tail.
Throughout the proofs, we drop the subscript-$N$ from the notation, considering sequences implicitly.
We now state the main theorem of this section.

\begin{thm}
\label{res-p5:profile:res}
	Let $(k_N)_\Ninn$ be a sequence of positive integers and $(G_N)_\Ninn$ a sequence of finite, Abelian groups;
	for each $\Ninn$, define $Z_{(N)} \cq [Z_1, ..., Z_{k_N}]$ by drawing $Z_1, ..., Z_{k_N} \sim^\iid \Unif(G_N)$.
	
	Suppose that the sequence $(k_N, G_N)_\Ninn$ satisfies \cref{hyp-p5:profile}.
	For all $\alpha \in \mbr$ and $N \in \mbn$,
	write $t_{\alpha,N} \cq t_\alpha(k_N, \abs{G_N})$.
	Then,
	for all $\alpha \in \mbr$,
	we have
	\[
		t_{\alpha,N} / t_{0,N}
	\to
		1
	\Qand
		d^\pm_{G_k, N}\rbr{ t_{\alpha,N} }
	\to^\mbp
		\Psi(\alpha)
	\Quad{(in probability)}
		\asinf N.
	\]
	That is, whp there is cutoff at time $t_0$ with profile given by $\bra{t_\alpha}_{\alpha \in \mathbb R}$: for all $\eps \in (0,1)$, the difference in the mixing times $\tmix(\eps) - \tmix(\tfrac12)$ is given, up to subleading order terms, by $t_{\Psi^{-1}(\eps)} - t_0$.
	Moreover, the implicit lower bound on the TV distance holds deterministically, ie for all choices of generators.
\end{thm}

\begin{rmkt}
\label{rmk-p5:profile:mixing-profile}
We can write the cutoff statement, emphasising the $N$-dependence, in the form
\[
	\rbb{ \tmix^{Z,N}(\eps) - t_{0,N} } / w_N
\to^\mbp
	\Psi^{-1}(\eps)
\Qforall
	\eps \in (0,1),
\]
where $(t_{0,N})_\Ninn$ is the mixing time and $(w_N)_\Ninn$ is the window, defined by \cref{res-p5:profile:t0a}.
Namely,
for all $\lambda \in (0,\infty)$ and all $\eps \in (0,1)$,
combining \cref{res-p5:profile:t0a,res-p5:profile:res},
we have
\[
	\frac{ \tmix^Z(\eps) - f(\lambda) k }{ g(\lambda) \sqrt k }
\to
	\Psi^{-1}(\eps)
\Qforall
	\eps \in (0,1)
\Qwhere
	k \eqsim \lambda \log \abs G.
\qedhere
\]
\end{rmkt}

\begin{rmkt*}
The CLT, \cref{res-p5:profile:CLT}, gives the dominating term in the TV distance \cref{res-p5:profile:res}:
\begin{itemize}[noitemsep, topsep = 0pt, label = \bcdot]
	\item 
	on the event $\bra{Q(t_\alpha) \le \log n - \omega}$, we lower bound the TV distance by $1 - \oh1$;

	\item 
	on the event $\bra{Q(t_\alpha) \ge \log n + \omega}$, we upper bound the expected TV distance by $\oh1$.
\end{itemize}
Combined with the CLT, we deduce that the $d_{G_k}(t_\alpha) \to \Psi(\alpha)$ in probability.
\end{rmkt*}

\begin{rmkt*}
\label{rmk-p5:profile:k=logn-interest}
The regime $k \asymp \log n$ is of particular interest.
It can be thought of as a `critical regime':
	if $k \ll \log n$, then $\tmix \gg k$;
	if $k \asymp \log n$, then $\tmix \asymp k$;
	if $k \gg \log n$, then $\tmix \ll k$.

From \cite[\cref{res-p2:intro:tv}]{HOt:rcg:abe:cutoff}, the following statements hold.
In order to have the cutoff time independent of the algebraic structure of the group, for $k$ with $1 \ll \log k \ll \log n$, it is sufficient to have either $k \gg \log n$ or $d(G) \ll \log n$. Further, there exist counter-examples when $k \asymp \log \abs G \asymp d(G)$:
	eg, the mixing times for $\mbz_2^{2r}$ and $\mbz_4^r$ are different if $k \eqsim \log n$, but not if $k \gg 2 \log_2 n$;
	note that $d(\mbz_2^{2r}) = 2r = \log_2 \abs{\mbz_2^{2r}}$ and $d(\mbz_4^r) = r = \tfrac12 \log_2 \abs{\mbz_4^r}$.
In this article, we place additional conditions on the group which are sufficient for us to not only establish cutoff but further obtain the limit profile, which is again independent of the algebraic structure of the group.
\end{rmkt*}

\subsection{Outline of Proof}
\renewcommand{\mm}{\ensuremath{\gamma}}

The outline here is very similar to that from the main article.
For a detailed outline, see \cite[\S\ref{sec-p2:cutoff1:outline}]{HOt:rcg:abe:cutoff} there; here we outline the difference.
Note that the lower bound in \cite[\S\ref{sec-p2:cutoff1:lower}]{HOt:rcg:abe:cutoff} was valid for all groups; we repeat it here for convenience.

For the upper bound, we were trying to bound the expectation of a $d$-th power of a gcd. Issues arose when $k$ became too large while $k - d$ is fairly small; see the proof of \cite[\cref{res-p2:cutoff1:gcd-ex}]{HOt:rcg:abe:cutoff}.
This arose from the fact that we used the following estimate from \cite[\cref{res-p2:cutoff1:divis}]{HOt:rcg:abe:cutoff}:
\[
	\pr{ V_1 \in \mm \mbz }
\le
	\pr{ V_1 \in \mm \mbz \mid V_1 \ne 0 } + \pr{ V_1 = 0 }
\le
	1/\mm + 2/n^{1/k}.
\]
The second term becomes an issue when $k$ gets close to $\log n$---eg if $k \gtrsim \log n$ then $n^{1/k}$ is a constant and the upper bound is simply ``at most some constant'', rather than decaying like $1/\mm$.
We alleviate this by defining
\[
	\mci
\cq
	\brb{ i \in [k] \mid V_i \ne 0 }
\]
and studying $\pr{ V_i \in \mm \mbz \mid i \in \mci }$; the problematic term $2/n^{1/k}$ then does not exist as we consider only non-zero coordinates of $V$.
If $G = \oplus_1^d \: \mbz_{m_j}$, then we are actually interested in $V_i \mod m_j$ for each $j$.
Recall that $m_* = \min_j m_j$.
`Typically', one has $\abs{V_i} \le m_*$. We suppose initially that $m_*$ is large enough so that $\max_i \abs{V_i} < m_*$ whp.
Thus looking at $V_i = 0$ or $V_i \equiv 0 \mod m_j$ is no different.

For large $\abs I$, the gcd analysis goes through similarly to before.
When $\abs I$ is small, eg smaller than $d$, it is more difficult to control; in this case, we use a fairly naive bound on the gcd, but control carefully the probability of realising such an $\mci$.
The case $\mci = \emptyset$ corresponding to $V = 0$, is handled using the concentration around the entropic time in exactly the same way as before.

\subsection{Lower Bound on Mixing}
\label{sec-p5:profile:lower}

In this subsection,
we prove the lower bound on mixing, which holds for every choice of $Z$.

In \cite[\S\ref{sec-p2:cutoff1}]{HOt:rcg:abe:cutoff}, we only considered $1 \ll k \ll \log n$.
As such, we only stated the entropic results for this regime.
Above, in \cref{res-p5:profile:t0a,res-p5:profile:CLT}, we stated analogous results for the full regime $1 \ll \log k \ll \log n$.
In the lower bound given in \cite[\S\ref{sec-p2:cutoff1:lower}]{HOt:rcg:abe:cutoff}, valid for arbitrary groups, there were no conditions on $k$ beyond those required for the entropic concentration, namely \cite[\cref{res-p2:cutoff1:CLT}]{HOt:rcg:abe:cutoff}.
As such, the identical proof passes over to the full regime $1 \ll \log k \ll \log n$ unchanged.
We include it below for completeness.

\begin{Proof}[Proof of Lower Bound in \cref{res-p5:profile:res}]
For this proof, assume that $Z$ is given, and suppress it.

We convert the CLT, \cref{res-p5:profile:CLT}, from a concentration statement about $Q$ into one about $W$:
for all $\alpha \in \mbr$,
by the CLT,
we have
\[
	\pr{\mce_\alpha}
\eqsim
	\Psi(\alpha)
\Qwhere
	\mce_\alpha
\cq
	\brb{ \mu\rbb{W(t_\alpha)} \ge n^{-1} e^\omega }
=
	\brb{ Q(t_\alpha) \le \log n - \omega };
\]
recall that $\omega \gg 1$.
Fix $\alpha \in \mbr$.
Consider the set
\[
	A_\alpha
\cq
	\brb{ x \in G \midb \exists \, w \in \mbz^d \st \mu_{t_\alpha}(w) \ge n^{-1} e^\omega \text{ and } x = w \bcdot Z }.
\]
Since we use $W$ to generate $S$, we have $\pr{ S(t_\alpha) \in A_\alpha \mid \mce_\alpha } = 1$.
Every element $x \in A_\alpha$ can be realised as $x = w_x \bcdot Z$ for some $w_x \in \mbz^k$ with $\mu_{t_\alpha}(w_x) \ge n^{-1} e^\omega$.
Hence, for all $x \in A_\alpha$, we have
\[
	\pr{ S(t_\alpha) = x }
\ge
	\pr{ W(t_\alpha) = w_x }
=
	\mu_{t_\alpha}(w_x)
\ge
	n^{-1} e^\omega.
\]
From this we deduce that
\[
	1
\ge
	\sumt{x \in A_\alpha}
	\pr{ S(t_\alpha) = x }
\ge
	\abs{A_\alpha} \cdot n^{-1} e^\omega,
\Quad{and hence}
	\abs{A_\alpha}/n \le e^{-\omega} = \oh1.
\]
Finally we deduce the lower bound from the definition of TV distance:
\[
	\tvb{ \pr{ S(t_\alpha) \in \cdot \mid Z } - \pi_G }
\ge
	\pr{ S(t_\alpha) \in A_\alpha } - \pi_G(A_\alpha)
\ge
	\prt{ \mce_\alpha } - \tfrac1n \abs{A_\alpha}
\ge
	\Psi(\alpha) - \oh1.
\qedhere
\]
\end{Proof}

\subsection{Upper Bound on Mixing}
\label{sec-p5:profile:upper}

Throughout this section we implicitly assume the conditions of \cref{hyp-p5:profile}, without repeating this in the statements below.
Similarly, we implicitly assume that $n = \abs G$ is sufficiently large.

We use a modified $L_2$ calculation, as in \cref{res-p5:p:mod-l2,def-p5:p:typ} in \S\ref{sec-p5:p:upper} above.
There we only bounded the order of the cutoff window; now we desire the profile.
We use definitions analogous to \cite[\cref{res-p2:cutoff1:mod-l2,def-p2:cutoff1:typ} in \S\ref{sec-p2:cutoff1:upper}]{HOt:rcg:abe:cutoff}, where the profile is studied.
Herein, we often suppress the time and $\alpha$-subscripts, eg writing $W$ for $W(t_\alpha)$ or $W(t)$, depending on context.


Let $W'$ be an independent copy of $W$; then $S' \cq W' \bcdot Z$ is an independent copy of $S$.
We recall the modified $L_2$ calculation; the following lemma is the same as \cref{res-p5:p:mod-l2}.

\begin{lem}
\label{res-p5:profile:mod-l2}
	For all $t \ge 0$ and all $\mcw \subseteq \mbz^k$,
	the following inequalities hold:
	\begin{subequations}
		\label{eq-p5:profile:mod-l2}
	\begin{gather*}
		\tvb{ \pr[Z]{ S(t) \in \cdot } - \pi_G }
	\le
		\tvb{ \pr[Z]{ S(t) \in \cdot \mid W(t) \in \mcw } - \pi_G }
	+	\pr{ W(t) \notin \mcw };
	\label{eq-p5:profile:mod-l2:triangle}
	\nt
	\\
		4 \, \ex{ \tvb{ \pr[Z]{ S(t) \in \cdot \mid W(t) \in \mcw } - \pi_G }^2 }
	\le
		\pr{ S(t) = S'(t) \mid W(t), W'(t) \in \mcw } - 1.
	\label{eq-p5:profile:mod-l2:tv-l2}
	\nt
	\end{gather*}
	\end{subequations}
\end{lem}



We now make the specific choice of the `typical' set $\mcw$; we make a different choice for each $\alpha \in \mbr$.
Write $\Psi$ for the standard Gaussian tail.
The collection $\bra{\mcw_\alpha}_{\alpha \in \mbr}$ of sets will satisfy
\[
	\pr{ W(t_\alpha) \notin \mcw_\alpha } \eqsim \Psi(\alpha),
\]
using the CLT (\cref{res-p5:profile:CLT}).
We show that the expression in \cref{eq-p5:profile:mod-l2:tv-l2} is $\oh1$.
Then applying \cref{eq-p5:profile:mod-l2:triangle} gives
\(
	d_{G_k}(t_\alpha)
\le
	\Psi(\alpha) + \oh1
\)
\whp.
This matches the lower bound in \S\ref{sec-p5:profile:lower}.

By considering all $\alpha \in \mbr$, we are able to find the shape of the cutoff. If we only desire the order of the window, then we need only consider the limit $\alpha \to \infty$; in this case, $\pr{ W(t_\alpha) \notin \mcw_\alpha } \approx \Psi(\alpha) \approx 0$, which explains the use of the word `typically' in describing $\mcw_\alpha$.

In order to define precisely the set $\mcw_\alpha$ here,
we first define two parameters, $r_\alpha$ and $p_\alpha$.

\begin{subtheorem}{thm}

\begin{defn}
\label{def-p5:profile:rp:def}
	For all $\alpha \in \mbr$, define $r_\alpha(k,n)$ and $p_\alpha(k,n)$ as follows:
	\[
		r_\alpha(k,n)
	&\cq
		\min
		\brb{ r \in \mbz_+ \midb \pr{ \absb{ W_1(t_\alpha) - \ex{W_1(t_\alpha)} } > r } \le 1/k^{3/2} };
	\\
		p_\alpha(k,n)
	&\cq
		\min
		\brb{ \pr{ W_1(t_\alpha) - \ex{W_1(t_\alpha)} = j } \midb \abs{j} \le r_\alpha(k,n) }.
	\]
	Also define
	\(
		r_*(k,n) \cq \tfrac12 \log k
	\Quad{and}
		p_*(k,n) \cq k^{-2}.
	\)
	(We consider $k \gtrsim \log n \gg 1$ in any application.)
\end{defn}

The typicality conditions will be a combination of `local' (coordinate-wise) and `global' ones.

\begin{defn}
\label{def-p5:profile:typ}
	For all $\alpha \in \mbr$, define the \textit{local} and \textit{global} typicality conditions, respectively:
	\[
		\mcw_{\alpha,\ell}
	&\cq
		\brb{ w \in \mbz^k \midb \absb{ w_i - \ex{W_1(t_\alpha)} } \le r_\alpha \: \forall \, i = 1,...,k };
	\\
		\mcw_{\alpha,g}
	&\cq
		\brb{ w \in \mbz^k \midb \pr{ W(t_\alpha) = w } \le n^{-1} e^{-\omega} }.
	\]
	Define $\mcw_\alpha \cq \mcw_{\alpha,\ell} \cap \mcw_{\alpha,g}$, and say that $w \in \mbz^k$ is ($\alpha$-)\textit{typical} if $w \in \mcw_\alpha$.
\end{defn}

\end{subtheorem}

The following proposition says that $r_\alpha \ge r_*$ and $p_\alpha \ge p_*$.

\begin{subtheorem}{thm}

\begin{prop}[{\cite[\cref{res-p0:rp:res,rmk-p0:rp:k-large}]{HOt:rcg:supp}}]
\label{res-p5:profile:rp:res}
	For all $\alpha \in \mbr$,
	we have
	\[
		r_\alpha(k,n) \ge r_*(k,n)
	\Qand
		p_\alpha(k,n) \ge p_*(k,n).
	\label{eq-p5:profile:rp:res}
	\nt
	\]
	(Recall that here we are considering $k \gtrsim \log n \gg 1$, by \cref{hyp-p5:profile}.)
\end{prop}

\begin{Proof}
	This follows from standard large deviation theory.
	Its proof can be found in \cite[\S\ref{sec-p0:rp}]{HOt:rcg:supp}.
\end{Proof}

The following proposition determines the probability that $W(t_\alpha)$ lies in $\mcw_\alpha$, ie of typicality.

\begin{prop}
\label{res-p5:profile:typ}
	For each $\alpha \in \mbr$,
	we have
	\[
		\pr{ W(t_\alpha) \notin \mcw_\alpha } \to \Psi(\alpha).
	\]
\end{prop}

\begin{Proof}
The probability that the global conditions hold converges to $1 - \Psi(\alpha)$ by our CLT, \cref{res-p5:profile:CLT}.
The probability that a single coordinate fails the local condition is at most $k^{-2}$ by \cref{def-p5:profile:rp:def,res-p5:profile:rp:res}.
Thus the probability that local typicality fails to hold is then at most $k^{-1} = \oh1$ by the union bound.
The claim follows.
\end{Proof}

\end{subtheorem}

Herein, we fix $\alpha \in \mbr$ and frequently suppress the $t_\alpha$ from the notation, eg writing $W_\cdot$ for $W_\cdot(t_\alpha)$ or $\mcw$ for $\mcw_\alpha$.
Let $V \cq W - W'$,
so
\(
	\bra{ W \bcdot Z = W' \bcdot Z }
=
	\bra{ V \bcdot Z = 0 }.
\)
Write
\[
	D
\cq
	D_\alpha
\cq
	n \, \pr{ V(t_\alpha) \bcdot Z = 0 \mid \typ_\alpha } - 1
\Qwhere
	\typ 
\cq
	\typ_\alpha
\cq
	\brb{ W(t_\alpha), W'(t_\alpha) \in \mcw_\alpha) }.
\]
It remains to show that $D_\alpha = \oh1$ for all $\alpha \in \mbr$.
Recall the conditions of \cref{hyp-p5:profile}:
\[
	\mint{j} m_j \gg \log k,
\quad
	k \gtrsim \log \abs G
\Qand
	d \le \tfrac1{30} \log \abs G / \log k.
\]

\begin{prop}
\label{res-p5:profile:l2}
	Suppose that $(k, G)$ jointly satisfy \cref{hyp-p5:profile}.
	(Recall that, implicitly, $(k, G)$ is a sequence of Abelian groups and integers.)
	Then, for all $\alpha \in \mbr$, we have
	\(
		D_\alpha = \oh1.
	\)
\end{prop}

Given this proposition, we can prove the upper bound in the main theorem, \cref{res-p5:profile:res}.

\begin{Proof}[Proof of Upper Bound in \cref{res-p5:profile:res} Given \cref{res-p5:profile:l2}]
\cref{hyp-p5:profile} implies that conditions required for \cref{res-p5:profile:l2}.
Apply the modified $L_2$ calculation, \cref{res-p5:profile:mod-l2,def-p5:profile:typ}, and use \cref{res-p5:profile:typ,res-p5:profile:l2} to control the two resulting terms.
Combined, these says that $d_{G_k}(t_\alpha) \le \Psi(\alpha) + \oh1$ \whp.
\end{Proof}

It remains to prove \cref{res-p5:profile:l2}, ie to bound the modified $L_2$ distance.
The remainder of the section is dedicated to this goal.

\medskip

Clearly we need to control the law of $V \bcdot Z$.
First we analyse the case $V = 0$, which immediately implies that $V \bcdot Z = 0$.
The global typicality condition is designed precisely to handle this.
The proof of the following lemma is deferred until the end of the subsection.

\begin{lem}
\label{res-p5:profile:V=0}
	We have
	\[
		n \, \pr{ V = 0 \mid \typ }
	\le
		e^{-\omega} / \prt{\typ}
	\]
\end{lem}

We move onto $V \ne 0$.
Linear combinations of independent uniform random variables in an Abelian group are themselves uniform on their support.
Hence the distribution of $v \bcdot Z$ is uniform on $\gcd(v_1, ..., v_k, n) G$.
Here we are using the convention
\(
	\gcd(r_1, ..., r_\ell,0) \cq \gcd(\abs{r_1}, ..., \abs{r_\ell})
\)
for $r_1, ..., r_\ell \in \mbz \setminus \bra{0}$.
The following lemma is proved rigorously in \cite[\cref{res-p0:deferred:unif-gcd}]{HOt:rcg:supp}.

\begin{lem}
\label{res-p5:profile:unif-gcd}
	For all $v \in \mbz^k$,
	we have
	\[
		v \bcdot Z \sim \Unif\rbb{ \mm G }
	\Qwhere
		\mm \cq \gcd(v_1, ..., v_k, n).
	\]
\end{lem}

We thus need to control $\abs{\mm G}$,
since \cref{res-p5:profile:unif-gcd} implies that
\[
	\pr{ V \bcdot Z = 0 \mid \typ }
=
	\sumt{\mm \in \mbn}
	\prt{ \mfgcd = \mm \mid \typ } / \abs{\mm G}
\Qwhere
	\mfgcd \cq \gcd\rbb{ V_1, ..., V_k, n }.
\]

\begin{lem}
\label{res-p5:profile:G/gammaG}
	For all Abelian groups $G$ and all $\mm \in \mbn$,
	we have
	\[
		\abs G / \abs{\mm G}
	\le
		\mm^{d(G)}.
	\]
\end{lem}

\begin{Proof}
	Decompose $G$ as $\oplus_1^d \: \mbz_{m_j}$ with $d = d(G)$ and some $m_1, ..., m_d \in \mbn$.
	Then $\mm G$ can be decomposed as $\oplus_1^d \: \gcd(\mm, m_j) \mbz_{m_j}$.
	Hence
	\(
		\abs{\mm G}
	=
		\prodt[d]{1} \rbr{ m_j/\gcd(\mm, m_j) }
	\ge
		\prodt[d]{1} \rbr{ m_j/\mm }
	=
		\abs G / \mm^d.
	\)
\end{Proof}

These lemmas combine to produce a simple, but key, corollary.

\begin{cor}
\label{res-p5:profile:pr-gcd}
	We have
	\[
		n \, \pr{ V \bcdot Z = 0, \: V \ne 0 \mid \typ }
	\le
		\ex{ \mfgcd^{d(G)} \, \one{V \ne 0} \mid \typ }.
	\]
\end{cor}

\begin{Proof}
The conditioning does not affect $Z$.
The corollary follows from \cref{res-p5:profile:unif-gcd,res-p5:profile:G/gammaG}.
\end{Proof}

Up to here, the proof has been very similar to that given in \cite[\S\ref{sec-p2:cutoff2:upper}]{HOt:rcg:abe:cutoff}.
Here it diverges somewhat, in the analysis of this gcd in \cref{res-p5:profile:pr-gcd}.
The conclusion is the same, though.

\medskip

Whether a coordinate is non-zero plays a major role in controlling this gcd.
We introduce the following notation, designed to handle the number of non-zero coordinates.
For $v \in \mbz^k$, write
\[
	\mci(v)
\cq
	\brb{ i \in [k] \midb v_i \not\equiv 0 \mod m_j \text{ for all } j = 1,...,d }.
\]
We always consider $V$ conditioned on typicality.
We have $\abs{V_i} \le 2r_* = \log k < m_j$ for all $i$ and $j$ local typicality and \cref{hyp-p5:profile}
Thus, conditioned on local typicality, we have
\[
	\mci(V) = \brb{ i \in [k] \midb V_i \ne 0 };
\Quad{abbreviate}
	\mci \cq \mci(V).
\]

We end up needing to separate the concepts of \emph{local} and \emph{global} typicality:
define
\[
	\typloc \cq \brb{ W, W' \in \mcw_\ell }
\Qand
	\typglo \cq \brb{ W, W' \in \mcw_g };
\Quad{then}
	\typ = \typloc \cap \typglo.
\]

We control the gcd via we determine the probability an individual coordinate is a multiple of a given number in the following auxiliary lemma;
it is taken from \cite[\S\ref{sec-p2:cutoff1:upper}]{HOt:rcg:abe:cutoff}.
Write $\alpha \wr \beta$ if $\alpha$~divides~$\beta$.


\begin{lem}
\label{res-p5:profile:divis}
	For
		all non-empty $I \subseteq [k]$ with $\bra{ \mci = I } \cap \typ \ne \emptyset$
	and
		all $\gamma \in \mbn$,
	we have
	\[
		\pr{ \gamma \wr V_i \: \forall \, i \in I \mid \mci = I, \: \typloc }
	\le
		\gamma^{-\abs I}.
	\]
\end{lem}

\begin{Proof}
The coordinates are independent and \emph{local} typicality merely conditions each coordinate to lie in a certain interval centred at 0.
The claim now follows immediately from \cite[\cref{res-p2:cutoff1:divis}]{HOt:rcg:abe:cutoff}.
\end{Proof}

From this, using the conditions of \cref{hyp-p5:profile}, we can deduce that
\(
	\ex{ \mfgcd^{d(G)} \one{V \ne 0} \mid \typ }
=
	1 + \oh1.
\)
We refer to it as a ``corollary'', since its proof is purely technical, not relying on any properties of the RW or the generators, just algebraic manipulations.
Its proof is briefly deferred.

\begin{cor}
\label{res-p5:profile:gcd-ex}
	Given \cref{hyp-p5:profile},
	we have
	\(
		\ex{ \mfgcd^{d(G)} \one{V \ne 0} \mid \typ }
	=
		1 + \oh1.
	\)
\end{cor}

\cref{res-p5:profile:l2} now follows immediately from \cref{res-p5:profile:V=0,res-p5:profile:pr-gcd,res-p5:profile:gcd-ex}.

\begin{Proof}[Proof of \cref{res-p5:profile:l2}]
By \cref{res-p5:profile:V=0,res-p5:profile:pr-gcd,res-p5:profile:gcd-ex},
we have
\[
	n \, \pr{ V \bcdot Z = 0 \mid \typ }
&
\le
	n \, \pr{ V = 0 \mid \typ }
+	n \, \pr{ V \bcdot Z = 0, \, V \ne 0 \mid \typ }
\\&
\le
	n \, \pr{ V = 0 \mid \typ }
+	\ex{ \mfgcd^{d(G)} \, \one{V \ne 0} \mid \typ }
=
	1 + \oh1.
\qedhere
\]
\end{Proof}

We now prove \cref{res-p5:profile:gcd-ex}.
Abbreviate $d \cq d(G)$.
We use the decomposition
\[
	\ex{ \mfgcd^d \one{ V \ne 0 } \midb \typ }
=
	\sumt{\emptyset \ne I \subseteq [k]}
	\ex{ \mfgcd^d \one{ V \ne 0 } \one{ \mci = I } \midb \typ }.
\]

First we analyse `large' $\mci$.
We defer the proof of the next lemma to the end of the subsection.

\begin{lem}
\label{res-p5:profile:gcd}
	For all $I \subseteq [k]$ with $\bra{ \mci = I } \cap \typ \ne \emptyset$ and $\abs I \ge d + 2$,
	we have
	\[
		\ex{ \mfgcd^d \midb \mci = I, \: \typ }
	\le
		1 + 3 \cdot 2^{d - \abs I} / \pr{\typglo \mid \mci = I, \: \typloc }.
	\]
\end{lem}

An easy corollary of this controls the expectation when $\mci$ is `large'.

\begin{cor}
\label{res-p5:profile:gcd-ex:large}
	Let $L \ge d+2$.
	We have
	\[
		n \sumt{\abs I \ge L}
		\ex{ \mfgcd^{d(G)} \one{ \mci = I } \midb \typ }
	\le
		1 + 3 \cdot 2^{d - L} / \prt{\typ}.
	\]
\end{cor}

\begin{Proof}[Proof of \cref{res-p5:profile:gcd-ex:large}]
Recall Bayes's rule,
specifically the fact that $\pr{B \mid C} / \pr{C \mid B} = \pr{B} / \pr{C}$ for non-null events $B$ and $C$.
By \cref{res-p5:profile:gcd},
for $L \ge d + 2$,
we deduce that
\[
&	n \sumt{\abs I \ge L}
	\ex{ \mfgcd^d \one{ \mci = I } \midb \typ }
\\&\qquad
\le
	\sumt{\abs I \ge L}
	\rbb{
		\pr{ \mci = I \midb \typ } + 3 \cdot 2^{d - \abs I} \, \pr{ \mci = I } / \prt{\typ}
	}
\\&\qquad
\le
	\pr{ \abs \mci \ge L \midb \typ }
+	3 \cdot 2^{d - L} \pr{ \abs \mci \ge L } / \prt{\typ}
\le
	1 + 3 \cdot 2^{d - L} / \pr{\typ}.
\qedhere
\]
\end{Proof}

Next we analyse `small' $\mci$.
This is somewhat more delicate: we use the trivial bound $\mfgcd \le r_* = \log k$ (conditional on typicality) and carefully $\pr{ \mci = I \mid \typ }$ using global typicality.
Again, we defer the proof of the next lemma to the end of the subsection.

\begin{lem}
\label{res-p5:profile:I}
	For all $I \subseteq [k]$,
	we have
	\[
		\pr{ \mci = I \midb \typ }
	\le
		n^{-1} e^{-\omega} / p_*^{\abs I} / \prt{\typ}
	=
		e^{-\omega} n^{-1} k^{2\abs I} / \prt{\typ}
	\label{eq-p5:profile:I-typ}
	\nt
	\]
\end{lem}

An easy corollary of this controls the expectation when $\mci$ is `small'.

\begin{cor}
\label{res-p5:profile:gcd-ex:small}
	Set $L \cq \tfrac15 \log n / \log k$.
	We have
	\[
		\sumt{1 \le \abs I \le L}
		\ex{ \mfgcd^{d(G)} \one{ \mci = I } \midb \typ }
	=
		\ohb{1/\prt{\typ}}.
	\]
\end{cor}

\begin{Proof}
Apply \cref{res-p5:profile:gcd,res-p5:profile:I} for $I$ with $1 \le \abs I \le L$:
\[
&	\prt{\typ}
	\sumt{1 \le \abs I \le L}
	\ex{ \mfgcd^d \one{ \mci = I } \midb \typ }
\le
	\sumt{1 \le \abs I \le L}
	(\log k)^d \cdot e^{-\omega} n^{-1} k^{2 \abs I}
\\&\qquad
\le
	L k^L \cdot (\log k)^d \cdot e^{-\omega} n^{-1} k^{2 L}
\le
	e^{-\omega} \cdot \log n \cdot n^{3/5} \cdot n^{1/5} \cdot n^{-1}
=
	\oh1,
\]
since the number of $I \subseteq [k]$ with $\abs I = \ell$ is $\binom k\ell \le k^\ell$ and $d \le \tfrac15 \log n / \log \log k$ by \cref{hyp-p5:profile}.
\end{Proof}

We now have all the ingredients required to prove \cref{res-p5:profile:gcd-ex} .

\begin{Proof}[Proof of \cref{res-p5:profile:gcd-ex}]
Combining
\cref{res-p5:profile:gcd-ex:large,res-p5:profile:gcd-ex:small}
with the decomposition
\[
	\ex{ \mfgcd^{d(G)} \one{ V \ne 0 } \midb \typ }
=
	\sumt{\emptyset \ne I \subseteq [k]}
	\ex{ \mfgcd^{d(G)} \one{ V \ne 0 } \one{ \mci = I } \midb \typ },
\]
gives
\(
	\ex{ \mfgcd^{d(G)} \one{ V \ne 0 } \midb \typ }
=
	1 + \oh{1/\prt{\typ}}.
\)
\cref{res-p5:profile:typ} states that $\prt{\typ} \asymp 1$.
\end{Proof}

It remains to give the deferred proofs of \cref{res-p5:profile:V=0,res-p5:profile:gcd,res-p5:profile:I}.

\begin{Proof}[Proof of \cref{res-p5:profile:V=0}]
By direct calculation, we have
\[
	\pr{ V = 0, \: \typ }
&
=
	\pr{ W = W', \: W \in \mcw }
\\&
=
	\sumt{w \in \mcw} \pr{W = w} \pr{W' = w}
=
	\sumt{w \in \mcw} \pr{W = w}^2,
\]
since $W$ and $W'$ are iid copies.
Recall global typicality: $\pr{W = w} \le n^{-1} e^{-\omega}$ for all $w \in \mcw$.
Thus
\[
	n \, \pr{ \mci = \emptyset \mid \typ }
\le
	n \sumt{w \in \mcw} \pr{W = w}^2 / \pr{\typ}
\le
	e^{-\omega} / \pr{\typ}.
\qedhere
\]
\end{Proof}

\begin{Proof}[Proof of \cref{res-p5:profile:gcd}]
Local typicality applies coordinatewise and so conditioning on it does not break independence.
This is not the case for global typicality.
We thus move from global to local:
\[
	\ex{ \mfgcd^d \midb \mci = I, \: \typ }
=
	1 + \ex{ \mfgcd^d - 1 \midb \mci = I, \: \typ }
=
	1 + \ex{ \mfgcd^d - 1 \midb \mci = I, \: \typloc } / \pr{ \typglo \midb \mci = I, \: \typloc }.
\]
Write $\pr[I]{\cdot}$ and $\ex[I]{\cdot}$ to denote probability and expectation, respectively, conditioned on $\mci = I$ and $\typloc$ (ie \emph{local} typicality).
If the gcd equals $\mm$, then certainly $\mm$ divides each coordinate.
Thus
\[
	\ex[I]{ \mfgcd^d }
=
	\sumt{\gamma \ge 2}
	\gamma^d \, \pr[I]{ \mfgcd = \gamma }
+	1
\le
	\sumt{\gamma \ge 2}
	\gamma^d \, \pr[I]{ \gamma \wr V_i \: \forall \, i \in I }
+	1.
\]
Applying \cref{res-p5:profile:divis},
for $I$ with $\bra{ \mci = I } \cap \typ \ne \emptyset$ and $\abs I \ge d + 2$,
we obtain
\[
	\ex{ \mfgcd^d \midb \mci = I, \: \typ }
\le
	1
+	\sumt{\gamma \ge 2} \gamma^{d - \abs I} / \pr[I]{\typglo}
\le
	1 + 3 \cdot 2^{d - \abs I} / \pr[I]{\typglo}.
\qedhere
\]
\end{Proof}

\begin{Proof}[Proof of \cref{res-p5:profile:I}]
Requiring $\mci = I$ places restrictions on the coordinates in $I^c$, but not on the coordinates of $I$ other than that they are non-zero; we ignore the latter to get an upper bound (see below).
For a vector $w \in \mbz^k$, write
\[
	\mcw_I(w) \cq \brb{ w' \in \mbz^k \mid \mci(w - w') = I }.
\]
Using the independence of $W$ and $W'$, we have
\[
	\pr{ \mci = I, \: W \in \mcw }
=
	\sumt{w \in \mcw}
	\pr{ W = w } \pr{ W' \in \mcw_I(w) }.
\]
Hence, using the independence of the coordinates of $W'$, given $w \in \mcw$ we have
\[
	\pr{ W' \in \mcw_I(w) }
=
	\pr{ W' = w } \cdot \prod_{i \in I} \frac{\pr{W_i' \ne w_i}}{\pr{W'_i = w_i}}
\le
	\pr{ W' = w } \cdot \prod_{i \in I} \frac1{\pr{W'_i = w_i}}.
\]
An immediate consequence of the definitions of $r$ and $p$, in \cref{def-p5:profile:rp:def}, is the following:
\[
	\text{for all $\alpha \in \mbr$},
\Quad{if}
	\absb{ w_1 - \ex{W_1(t_\alpha)} } \le r_\alpha(k,n)
\Quad{then}
	\pr{ W_1(t_\alpha) = w_1 } \ge p_\alpha(k,n).
\]
By \cref{def-p5:profile:rp:def,res-p5:profile:rp:res}, we have $p_\alpha \ge p_* = k^{-2}$.
Hence, for $w \in \mcw$, we obtain
\[
	\pr{ W' \in \mcw_I(w) }
\le
	\pr{ W' = w } / p_*^{\abs I}
\le
	n^{-1} e^{-\omega} / p_*^{\abs I}
=
	n^{-1} e^{-\omega} / k^{2 \abs I}
\]
From this and the sum above, \cref{res-p5:profile:I} follows by summing over all $w \in \mcw$:
\[
	\pr{ \mci = I, \: \typ }
\le
	\pr{ \mci = I, \: W \in \mcw }
\le
	n^{-1} e^{-\omega} k^{2 \abs I} \sumt{w \in \mcw} \pr{W = w}
\le
	n^{-1} e^{-\omega} k^{2 \abs I}.
\qedhere
\]
\end{Proof}

\section{Cutoff: A Detailed Investigation of $\mathbb Z_\pp^d$}
\label{sec-p5:p}

In this section we perform a detailed analysis of the behaviour of the mixing time for the random walk on the uniform random Cayley graph of degree $k$ of $\mbz_\pp^d$, proving \cref{res-p5:intro:p}.
In \cite{HOt:rcg:abe:cutoff} we established cutoff, under the assumption that $k - d \gg 1$.
If $G = \mbz_\pp^d$ for $\pp$ prime,~then
\[
	\brb{ \mm G \mid \mm \wr n }
=
	\brb{ \mm G \mid \mm \wr \pp }
=
	\bra{ G } \cup \bra{ \pp G }
=
	\bra{ G, \: \bra{ \id } };
\]
that is, the only options are the group itself and the trivial group, corresponding to $\mm = 1$ and $\mm = \pp$, respectively.
Thus, applying \cite[\cref{res-p2:cutoff2:res}]{HOt:rcg:abe:cutoff}, we deduce that there is cutoff at the entropic time $t_0(\pp, \abs G) = t_0(\pp, \pp^d)$, ie the time at which the entropy of the RW on $\mbz_\pp^k$ becomes $\log \abs G = d \log \pp$.

In this exposition, we consider some cases not covered in \cite{HOt:rcg:abe:cutoff}. In particular, we allow $k - d$ to be a fixed constant, not diverging.
When this is the case and $\pp$ diverges, it can be shown that choosing $k$ elements $Z_1, ..., Z_k \sim^\iid \Unif(\mbz_\pp^d)$ generates the group whp.
On the other hand, if $\pp$ is also a fixed constant, then this is not the case; we establish cutoff conditional on generating the group.
See \cref{res-p5:gen:dichotomy} for a dichotomy giving necessary and sufficient conditions to generate the~group.

To ease notation, we drop completely any $\pp$-s, and often drop the $\pp$; to be explicit, we state in the next subsection precisely what notation we are going to use.

\subsection{Entropic Times: Methodology, Definition and Concentration}
\label{sec-p5:p:ent}

We use an `entropic method'; for further details, see \cite{HOt:rcg:abe:cutoff}.
To make this as self-contained as possible, we now explain the specific application in a little more depth.

\smallskip

We define an auxiliary random process $(W(t))_{t\ge0}$, recording how many times, mod $\pp$, each generator has been used:
	for $t \ge 0$, for each generator $i = 1,...,k$, write $W_i(t)$ for the number of times that it has been picked by time $t$.
By independence, $W(\cdot)$ forms a rate-1 DRW on $\mbz_\pp^k$.
For the undirected case, recall that we either apply a generator or its inverse; when we apply the inverse of generator $i$, increment $W_i \to W_i - 1$ (rather than $W_i \to W_i + 1$).
In this case, $W(\cdot)$ is a SRW (rather than DRW) on $\mbz_\pp^k$.
Note that every element of $G = \mbz_\pp^d$ has order $\pp$, since $\pp$ is prime. Hence it suffices to look at the walk $W$ mod $\pp$, ie on $\mbz_\pp^k$, rather than on $\mbz^k$.

Since the underlying group is Abelian, the order in which the generators are applied is irrelevant and generator-inverse pairs cancel; hence we can write
\[
	S(t)
=
	\sumt[k]{i=1} W_i(t) Z_i
=
	W(t) \bcdot Z.
\]

Recall that the invariant distribution is uniform on $G$, giving mass $1/n$ to each vertex.
The proposed mixing time is then the time at which the auxiliary process $W$ obtains entropy $\log n$.
This time will be calculated fairly precisely in many situations; see Proposition~\ref{res-p5:p:t0a}.

Write $\mu_t$, respectively $\nu_s$, for the law of $W(t)$, respectively $W_1(sk)$;
so
\(
	\mu_t = \nu_{t/k}^{\otimes k}.
\)
Define
\[
	Q(t) \cq - \log \mu_t\rbb{W_i(t)}
\Quad{and}
	Q_i(t) \cq - \log \nu_{t/k}\rbb{W_i(t)};
\]
then, $Q_i$ forms an iid sequence over $i \in [k]$, and
\[
	Q(t) = \sumt[k]{i=1} Q_i(t),
\quad
	h(t) \cq \ex{ Q(t) }
\Qand
	H(s) \cq \ex{ Q_1(sk) }.
\]
So $h(t)$ and $H(s)$ are the entropies of $W(t)$ and $W_1(sk)$, respectively.
Note that $h(t) = k H(t/k)$ and that $h : [0,\infty) \to [0, \log(\pp^k))$ is a strictly increasing bijection.

While all our results can be phrased in terms of (Shannon) entropy,
from a technical point of view it will be convenient to define the \textit{relative entropy}:
\[
	R(s) \cq \log \pp - H(s).
\]
The maximal entropy of a random variable on $\mbz_\pp$ is $\log \pp$, obtained uniquely by the uniform distribution.
Since the RW converges to the uniform distribution, $R(s) \to 0$ \asinf s.
Of great importance will be the parameter
\[
	\zeta \cq \tfrac1k (k - d) \log \pp = \log \pp - \tfrac1k \log n
\Qwhere
	n
\cq
	\abs{\mbz_\pp^d} = \pp^d.
\]

\begin{defn}
\label{def-p5:p:entropic-time}
Define
\(
	\zeta
\cq
	\tfrac1k (k - d) \log \pp
=
	\log \pp - \tfrac1k \log n,
\)
and,
for $\alpha \in \mbr$,
define
\[
	\zeta_\alpha
\cq
	\zeta \rbb{ 1 - 2 \alpha / \sqrt{\zeta k (\zeta \vee 1)} },
\quad
	s_\alpha
\cq
	H^{-1}(\zeta_\alpha)
\Qand
	t_\alpha
\cq
	s_\alpha k.
\]
	%
\end{defn}

\subsection{Entropic Times: Evaluation and Concentration}

In this subsection, we estimate the entropic times in different regimes, and give a concentration result.
The proofs are given in \cite[\S\ref{sec-p0:re:p}]{HOt:rcg:supp}; precise references are given at the appropriate times.
Importantly, the arguments for Propositions~\ref{res-p5:p:t0a} and~\ref{res-p5:p:conc} \emph{do not} require $p$ to be prime.

\smallskip

The first proposition estimates the entropic times $t_0$ and the difference $t_\alpha - t_0$; the second gives concentration of the $Q$ random variable around these times.
Recall that $\zeta = \tfrac1k (k - d) \log \pp$.

\begin{subtheorem}{thm}
	\label{res-p5:p:t0a}

\begin{prop}[{\cite[\cref{res-p0:re:p:t0}]{HOt:rcg:supp}}]
\label{res-p5:p:t0}
	Suppose that $1 \ll k \lesssim d \log \pp$.
	The following hold:
	\begin{alignat*}{2}
		\text{if}
	\quad
		\zeta &\eqmathsbox{p:0}{\ll} 1,
	&\Quad{then}
		t_0/k = s_0 &\eqsim \tfrac12 \log(1/\zeta) / \rbb{ 1 - \cos(2\pi/\pp) };
	\\
		\text{if}
	\quad
		\zeta &\eqmathsbox{p:0}{\gtrsim} 1,
	&\Quad{then}
		t_0/k = s_0 &\asymp \pp^2 e^{-2\zeta} = (\pp^d)^{2/k};
	\intertext{further, if in fact $1 \ll k \ll d \log \pp$, then}
		\text{if}
	\quad
		\zeta &\eqmathsbox{p:0}{\gg} 1,
	&\Quad{then}
		t_0/k = s_0 &\eqsim \pp^2 e^{-2\zeta} / (2 \pi e) = (\pp^d)^{2/k} / (2 \pi e).
	\end{alignat*}
	Note that
	\(
		1 - \cos(2\pi/\pp)
	\eqsim_{\toinf \pp}
		2 \pi^2 / \pp^2
	=
		2 \pi^2 \pp^{-2d/k} e^{2\zeta}.
	\)
\end{prop}

\begin{prop}[{\cite[\cref{res-p0:re:p:ta}]{HOt:rcg:supp}}]
\label{res-p5:p:ta}
	Suppose that $1 \ll k \lesssim d \log \pp$ and $(k - d) \pp \gg 1$, ie $\zeta \gg 1/k$.
	Then,
	for all $\alpha \in \mbr$,
	we have
	\(
		t_\alpha \eqsim t_0
	\)
	and furthermore the following hold:
	\begin{alignat*}{3}
		\text{if}
	\quad
		\zeta &\eqmathsbox{p:a}{\lesssim} 1,
	\Quad{then}&
		(t_\alpha - t_0) / t_0 &\lesssim 1 / \rbb{ \sqrt{\zeta k} \log( (1/\zeta) \vee e ) } = \oh1;
	\\
		\text{if}
	\quad
		\zeta &\eqmathsbox{p:a}{\gg} 1,
	\Quad{then}&
		(t_\alpha - t_0) / t_0 &\lesssim 1 / \sqrt k = \oh1
	\quad
		\text{for the SRW}.
	\end{alignat*}
\end{prop}

\end{subtheorem}

\begin{prop}[Concentration, {\cite[\cref{res-p0:re:p:conc}]{HOt:rcg:supp}}]
\label{res-p5:p:conc}
	For $\alpha \in \mbr$, define
	\[
		Q_\alpha^+ \cq \bra{ Q(t_{ \alpha}) \ge d \log p + \alpha \sqrt{ k (\zeta \wedge 1) } }
	\Qand
		Q_\alpha^- \cq \bra{ Q(t_{-\alpha}) \le d \log p - \alpha \sqrt{ k (\zeta \wedge 1) } };
	\]
	For all $\alpha \in (0,\infty)$ with $\abs{ \zeta_\alpha - \zeta_0 } \le \tfrac12 \zeta_0$,
	we have
	\(
		\pr{ (Q_\alpha^\pm)^c } \lesssim \alpha^{-2}.
	\)
\end{prop}

\subsection{Precise Statement and Remarks}

Recall that $d_{G_k}(t)$ is the TV distance from uniform after time $t$ for the RW on $G_k$.

\begin{thm}[Cutoff]
\label{res-p5:p:res}
	Let $G$ be a finite, Abelian group admitting a decomposition $G \cq \mbz_\pp^d$ with $\pp$ prime.
	Assume that $1 \ll k \lesssim d \log \pp$.
	Define the entropic times $\bra{t_\alpha}_{\alpha \in \mbr}$
	as in \cref{def-p5:p:entropic-time}.
	The entropic times are asymptotically evaluated in Proposition~\ref{res-p5:p:t0a}.
	
	Suppose that $(k - d) \pp \gg 1$, ie $\zeta k \gg 1$.
	Then the RW on $G_k$ exhibits cutoff whp at $t_0$.
	More precisely,
		choose a sequence $(\beta_N)_\Ninn \subseteq \mbr_+^\mbn$ with $\beta_N \to \infty$ (arbitrarily slowly)
	and
		let $c \in \bra{\pm1}$.~%
	Then
	\[
		d^\pm_{G_k,N}\rbr{ t_{c \beta_N} }
	\to^\mbp
		\one{c = -1}
	\Quad{(in probability)}
		\asinf N.
	\]
	Moreover, the implicit lower bound holds deterministically, ie for all choices of generators.
	
	Also, if $0 \le (k - d)p \lesssim 1$, then, conditional that the uniformly chosen multisubset $[Z_1, ..., Z_k]$ generates the group, there is cutoff whp at time
	\(
		\tfrac12 d \log d / \rbr{1 - \cos(2\pi/\pp)}.
	\)
\end{thm}


\begin{rmkt}
\label{rmk-p5:p:extensions}
The outline of the proof is the same for all $\zeta$ with $\zeta k \gg 1$; we assume this initially.

We also consider the case where $0 \le k - d = \Oh1$ conditional on generating the group. This uses a standard argument for the case $k = d$, and then compares the walk which uses $Z = [Z_1, ..., Z_k]$ with another walk using a subset $Z'$ of size $d$ which generates the group.

We explain how to do this at the end of the section in \S\ref{sec-p5:p:ext}.
\end{rmkt}

\begin{rmkt*}
	Prior to our work, cutoff had already been established for any Abelian group when $k \gg \log n$, with an explicit mixing time; see \cite[\S\ref{sec-p2:intro:previous-work:ad-conj}]{HOt:rcg:abe:cutoff}.
	Although our technique can be adapted to allow $k \gtrsim \log n$ (when $n = \abs{\mbz_p^d}$, this is equivalent to $k \gtrsim d \log p$), we do not give details here.
\end{rmkt*}

\subsection{Lower Bound on Mixing for $\mbz_\pp^d$}
\label{sec-p5:p:lower}

In this subsection,
we prove the lower bound on mixing, which holds for every choice of $Z$.
(This argument is almost identical to the one which we give in \cite[\S\ref{sec-p2:cutoff1:lower}]{HOt:rcg:abe:cutoff}.)

\begin{Proof}[Proof of Lower Bound]
For this proof, we assume that $Z$ is given, and suppress it.

The concentration result \cref{res-p5:p:conc} gives
\(
	\pr{Q_\alpha^-} \to 1
\Quad{as}
	\alpha \to \infty.
\)
Consider the set
\[
	A_\alpha
\cq
	\brb{ x \in G \midb \exists \, w \in \mbz^d \st \mu_{t_\alpha}(w) \ge n^{-1} e^\omega \text{ and } x = w \bcdot Z }.
\]
Since we use $W$ to generate $S$, we have $\pr{ S(t_\alpha) \in A_\alpha \mid \mce_\alpha } = 1$.
Every element $x \in A_\alpha$ can be realised as $x = w_x \bcdot Z$ for some $w_x \in \mbz^k$ with $\mu_{t_\alpha}(w_x) \ge n^{-1} e^\omega$.
Hence, for all $x \in A_\alpha$, we have
\[
	\pr{ S(t_\alpha) = x }
\ge
	\pr{ W(t_\alpha) = w_x }
=
	\mu_{t_\alpha}(w_x)
\ge
	n^{-1} e^\omega.
\]
From this we deduce that
\[
	1
\ge
	\sumt{x \in A_\alpha}
	\pr{ S(t_\alpha) = x }
\ge
	\abs{A_\alpha} \cdot n^{-1} e^\omega,
\Quad{and hence}
	\abs{A_\alpha}/n \le e^{-\omega} = \oh1.
\]
Finally we deduce the lower bound from the definition of TV distance:
\[
	\tvb{ \pr{ S(t_\alpha) \in \cdot \mid Z } - \pi_G }
\ge
	\pr{ S(t_\alpha) \in A_\alpha } - \pi_G(A_\alpha)
\ge
	\pr{ Q_\alpha^- } - \tfrac1n \abs{A_\alpha}
=
	1 - \oh1.
\qedhere
\]
\end{Proof}

\subsection{Upper Bound on Mixing for $\mbz_\pp^d$ for $(k - d)p \gg 1$}
\label{sec-p5:p:upper}

This subsection is devoted to the upper bound.

\begin{Proof}[Outline of Proof]
\qedtriangle
	Consider $\alpha$ with $\alpha \to \infty$, but arbitrarily slowly.
	We show that the TV distance from uniform is $\oh1$ on the event $Q_\alpha^+$ (whp over $Z$) and that $\pr{Q_\alpha^+} = 1 - \oh1$, using similar techniques to those in \cite[\S\ref{sec-p2:cutoff1:upper}]{HOt:rcg:abe:cutoff}.
	\cref{res-p5:p:res} follows from this and Propositions~\ref{res-p5:p:t0a} and~\ref{res-p5:p:conc}.
\end{Proof}

We now make this outline precise and rigorous.
Herein, we frequently suppress the time and $\alpha$-subscripts, eg writing $W$ for $W(t_\alpha)$ or $W(t)$, depending on context.

Key is a `modified $L_2$ calculation'; cf \cite[\cref{res-p2:cutoff1:mod-l2}]{HOt:rcg:abe:cutoff}.
In short, one condition that $W$ is `typical' (in some precise sense), and then applies the standard TV--$L_2$ calculation on the conditioned law.

Let $W'$ be an independent copy of $W$; then $S' \cq W' \bcdot Z$ is an independent copy of $S$.

\begin{lem}
\label{res-p5:p:mod-l2}
	For all $t \ge 0$ and all $\mcw \subseteq \mbz_\pp^k$,
	the following inequalities hold:
	\begin{subequations}
	\begin{gather*}
		\tvb{ \pr{ S(t) \in \cdot \mid Z } - \pi_G }
	\le
		\tvb{ \pr[Z]{ S(t) \in \cdot \mid W(t) \in \mcw } - \pi_G }
	+	\pr{ W(t) \notin \mcw };
	\label{eq-p5:p:mod-l2:triangle}
	\nt
	\\
		4 \, \ex{ \tvb{ \pr[Z]{ S(t) \in \cdot \mid W(t) \in \mcw } - \pi_G }^2 }
	\le
		n \, \pr{ S(t) = S'(t) \mid W(t), W'(t) \in \mcw } - 1.
	\label{eq-p5:p:mod-l2:tv-l2}
	\nt
	\end{gather*}
	\end{subequations}
\end{lem}

\begin{Proof}
The first claim follows immediately from the triangle inequality.
For the second, using Cauchy--Schwarz, we upper bound the TV distance of the conditioned law by its $L_2$ distance:
\[
&
	4 \, \tvb{ \pr[Z]{ S \in \cdot \mid W \in \mcw } - \pi_G }^2
\le
	n \sumt{x} \rbb{ \pr[Z]{ S = x \mid W \in \mcw } - \tfrac1n }^2
\\&
\qquad
=
	n \sumt{x} \pr[Z]{ S = x \mid W \in \mcw }^2 - 1
=
	n \sumt{x} \pr[Z]{ S = S' = x \mid W,W' \in \mcw } - 1,
\]
as $S = W \bcdot Z$, $S' = W' \bcdot Z$ and $V = W - W'$.
The claim follows from Jensen's inequality.
\end{Proof}

We now make the specific choice of the `typical' set $\mcw$; we make a different choice for each $\alpha \in \mbr$. Cf \cite[\cref{def-p2:cutoff1:typ}]{HOt:rcg:abe:cutoff}.
The collection $\bra{\mcw_\alpha}_{\alpha \in \mbr}$ will satisfy
\[
	\pr{ W(t_\alpha) \notin \mcw_\alpha }
\approx
	0
\quad
	\text{for large $\alpha$},
\]
using the concentration result \cref{res-p5:p:conc}.
We show that the expression \cref{eq-p5:p:mod-l2:tv-l2} is $\oh1$.
Then applying \cref{eq-p5:p:mod-l2:triangle} gives
\(
	d_{G_k}(t_\alpha) \approx 0
\)
\whp,
for large $\alpha$.

\begin{defn}
\label{def-p5:p:typ}
	For all $\alpha \in \mbr$,
	define
	\(
		\omega_\alpha
	\cq
		\alpha \sqrt{ k (\zeta \wedge 1) } \gg 1,
	\)
	\[
		\mcw_\alpha
	\cq
		\brb{ w \in \mbz_\pp^k \mid \pr{ W(t_\alpha) = w } \le n^{-1} e^{-\omega_\alpha} }
	\Qand
		\typ_\alpha
	\cq
		\brb{ W(t_\alpha), W'(t_\alpha) \in \mcw_\alpha }.
	\]
\end{defn}

The following proposition determines the probability that $W(t_\alpha)$ lies in $\mcw_\alpha$, ie of typicality.

\begin{prop}
\label{res-p5:p:typ}
	For all $\alpha \in (0,\infty)$ with $\abs{ \zeta_\alpha - \zeta_0 } \le \tfrac12 \zeta_0$,
	we have
	\[
		\pr{ W(t_\alpha) \notin \mcw_\alpha }
	\lesssim
		\alpha^{-2}.
	\]
\end{prop}

\begin{Proof}
The lemma follows immediately from \cref{res-p5:p:conc}, since
\(
	\bra{ W(t_\alpha) \in \mcw_\alpha }
=
	Q^+_\alpha.
\)
\end{Proof}

Herein, inside proofs we often drop the time dependence and $\alpha$-subscripts from the notation, eg writing $W$ for $W(t_\alpha)$ and $\mcw$ for $\mcw_\alpha$ or $\typ$ for $\typ_\alpha$.
The `typical set' $\mcw$ is designed precisely so that the following lemma holds.

\begin{lem}
\label{res-p5:p:W=W'}
	For all $\alpha \in (0,\infty)$,
	we have
	\[
		\pr{ W(t_\alpha) = W'(t_\alpha) \mid \typ_\alpha }
	\le
		n^{-1} e^{-\omega_\alpha} / \pr{\typ_\alpha}
	\ll
		n^{-1}.
	\]
\end{lem}

\begin{Proof}
By direct calculation,
using independence of $W$ and $W'$,
we have
\[
	\pr{ W = W', \: \typ }
=
	\pr{ W = W', \: W \in \mcw }
=
	\sumt{w \in \mcw} \pr{W = w}^2
\le
	n^{-1} e^{-\omega},
\]
with the final inequality using global typicality.
The result follows by Bayes's rule.
\end{Proof}


When $W = W'$, we necessarily have $S = S'$ (since the group is Abelian).
Now consider when $W \ne W'$.
The following lemma is a special case of \cite[\cref{res-p2:cutoff1:unif-gcd}]{HOt:rcg:abe:cutoff}.

\begin{lem}[{\cite[\cref{res-p2:cutoff1:unif-gcd}]{HOt:rcg:abe:cutoff}}]
\label{res-p5:p:unif-gcd}
	For any $v \in \mbz_\pp^k \setminus \bra{0}$,
	we have
	\(
		v \bcdot Z \sim \Unif(G).
	\)
\end{lem}

\begin{cor}
\label{res-p5:p:W=/=W'}
	For all $\alpha \in \mbr$,
	we have
	\[
		\pr{ S(t_\alpha) = S'(t_\alpha), \: W(t_\alpha) \ne W'(t_\alpha) \mid \typ_\alpha }
	\le
		\tfrac1n.
	\]
\end{cor}

\begin{Proof}
Condition on $W = w$ and $W' = w'$ with $w \ne w'$ and $w,w' \in \mcw$.
This conditioning is independent of $Z$.
Hence
\(
	S - S' = (w - w') \bcdot Z \sim \Unif(G)
\)
by \cref{res-p5:p:unif-gcd}.
The claim follows.
\end{Proof}


\begin{Proof}[Proof of Upper Bound in \cref{res-p5:p:res} Given Propositions~\ref{res-p5:p:t0a} and~\ref{res-p5:p:conc}]
We are assuming that $(k - d) p \gg 1$, ie $\zeta k \gg 1$.
This means that $\abs{ \zeta_\alpha - \zeta_0 } \le \tfrac12 \zeta_0$ for all $\alpha \in \mbr$.
Choose $\alpha$ with $\alpha \to \infty$, arbitrarily slowly.
We need to show that $d_{G_k}(t_\alpha) = \oh1$ \whp.

Apply \cref{res-p5:p:W=W',res-p5:p:W=/=W'} to deduce that $\pr{S = S' \mid \typ} = \oh1$.
Apply the modified $L_2$ calculation of \cref{res-p5:p:mod-l2} using \cref{def-p5:p:typ} for the definition of typicality.
Bound the `error term' using \cref{res-p5:p:typ}.
This gives $d_{G_k}(t_\alpha) = \oh1$ \whp.
\end{Proof}

\subsection{Removing the Condition $(k - d) \pp \gg 1$}
\label{sec-p5:p:ext}

In this subsection, we explain how to remove the condition $(k - d) \pp \gg 1$ in \cref{res-p5:p:res} under the conditioning that the group is generated, as referenced in \cref{rmk-p5:p:extensions}.
There are two cases to consider:
	$k = d$ with $\pp$ arbitrary (allowed to diverge)
and
	$0 < k - d = \Oh1$ with $\pp$ a fixed prime.

\begin{Proof}[Case $k = d$]
Here we do not need to assume that $\pp$ is prime.
Note also that $d = k \gg 1$.
The occurrence of cutoff in this set-up is not difficult.
As we could not find a proof in the literature---note that $\pp$ is not assumed to be fixed---we give the details. 

A key observation is that if $Z'$ is a set of size $d$ that generates $\mbz_\pp^d$, then the Cayley graph with respect to it is isomorphic to the Cayley graph with respect to the standard basis $\bra{e_1, ..., e_d}$.
Namely, it is the $d$-fold Cartesian product chain of the $\pp$-cycle with itself.

\smallskip

For the lower bound, we combine the method of distinguishing statistics and Wilson's method \cite{W:wilson-method}.
Let $f_2(y) \cq \cos(2 \pi y/\pp)$ and $\lambda_\pp \cq \cos(2 \pi / \pp)$.
Then $\lambda_\pp$ is the second largest eigenvalue of the transition matrix of SRW on $\mbz_\pp$, and the corresponding eigenvector is $f_2$.
Then
\(
	f(x_1,...,x_d)
\cq
	\tfrac1d
	\sumt[d]{i=1} f_2(x_i)
\)
is an eigenvector of the transition matrix of SRW on $\mbz_\pp^d$ with eigenvalue
\[
	\Lambda_{\pp,d}
\cq
	(d-1)/d + \lambda_\pp/d
=
	1 - (1 - \lambda_\pp)/d.
\]

We use initial state $(0, ..., 0)$.
To apply the method of distinguishing statistics, we need to bound both the expectation and the variance of $f$, both under the uniform and the RW distributions (at time $t$); write these as $\pi$ and $\mu_t$, respectively.
Under the uniform distribution, since the coordinates are independent and $\abs{f_2(z)} \le 1$ for all $z \in \mbz_\pp$, we have
\(
	\VAR[\pi]{f} \le \tfrac1d;
\)
similarly, we have
\(
	\VAR[\mu_t]{f} \le \tfrac1d.
\)
Also, since $f$ is an eigenvector, we have
\(
	E_\pi(f) = 0
\)
and
\(
	E_{\mu_t}(f) = e^{-(1-\Lambda_{\pp,d})t}.
\)
Applying the method of distinguishing statistics, eg as stated in \cite[Proposition~7.12]{LPW:markov-mixing},
for all $\eps \in (0,1)$,
whenever
\(
	e^{-(1-\Lambda_{\pp,d})t} \ge C_\eps/\sqrt d,
\)
for a sufficient large constant $C_\eps$,
we have
\(
	t \le \tmix(1 - \eps).
\)

Rearranging
\(
	e^{-(1-\Lambda_{\pp,d})t} \ge C_\eps/\sqrt d
\)
and recalling the definitions of $\Lambda_{\pp,d}$ and $\lambda_\pp$,
we obtain
\[
	t
\le
	\frac{\frac12 \log d - \log C_\eps}{1 - \Lambda_{p,d}}
=
	\frac{\frac12 d \log d - d \log C_\eps}{1 - \lambda_\pp}
=
	\frac{\frac12 d \log d}{1 - \lambda_\pp} \cdot \rbb{ 1 - \oh1 }.
\]

\smallskip

We now prove a matching upper bound on the mixing time.
Let $P_t$ be the time-$t$ transition probabilities of the walk on $\mbz_\pp^d$.
We identify this walk with the aforementioned $d$-fold Cartesian product of the $\pp$-cycle with itself.
Using independence of the coordinates, we have
\[
	\pp^d \, P_{2t}(x,x) - 1
=
	\rbb{ \pp \, Q_{2t/d}(0,0) }^d - 1,
\]
where $Q_s$ is the time-$s$ transition kernel for a rate-1 SRW on $\mbz_\pp$.
This is the $L_\infty$ distance at time $2t$, and hence the square of the $L_2$ distance at time $t$.
Let $\eps > 0$ be a constant.
Using the fact that
\[
	\rbr{ 1 + \tfrac12 \eps^2/d }^d
\le
	\eps^2
\]
when $\eps$ is sufficiently small, we have
\[
	\pp^d \, P_{2t}(x,x) - 1
\le
	\eps
\Qwhere
	\pp \, Q_{2t/d}(0,0) - 1
\le
	\tfrac12 \eps^2/d.
\]
Using the eigenvalue representation,
\[
	\text{if}
\quad
	2t/d \ge (1 + \eps)\log(2d/\eps^2) / (1 - \lambda_\pp),
\Quad{then}
	\pp \, Q_{2t/d}(0,0) - 1
\le
	\tfrac12 \eps^2/d,
\]
and hence
\(
	\pp^d P_{2t}(x,x) - 1
\le
	\eps^2.
\)
Thus
\[
	\text{if}
\quad
	t
\ge
	\frac{\frac12 (1 + \eps) d \log d + \tfrac12 (1 + \eps) d \log(2/\eps^2)}{1 - \lambda_\pp}
=
	\frac{\frac12 d \log d}{1 - \lambda_\pp} \cdot \rbb{ 1 + \oh1 }
\Quad{then}
	d_2(t) \le \eps,
\]
provided $\eps$ is sufficiently small.
This upper bound matches our lower bound.
\end{Proof}

\begin{Proof}[Case $0 < k - d = \Oh1$]
When $\pp \gg 1$, we already established cutoff, and so the group is generated, whp. Thus we may assume that both $\pp$ and $k-d$ are fixed, but $k \ge d \gg 1$.

\smallskip

The following statement is key:
	if $Z_1,...,Z_k$ generate $\mbz_\pp^d$ for $\pp$ prime, then there exists a set $S \subseteq [k]$ such that $\abs S = d$ and $\bra{Z_\ell}_{\ell \in S}$ generate $\mbz_\pp^d$.
This is immediate by viewing $\mbz_\pp^d$ as a vector space over the field $\mbf_\pp$ and noting that a set generates $\mbz_\pp^d$ if and only if it spans it.

\smallskip

We begin by obtaining an upper bound.
Choose a subset of generators of size $d$ which generate the group.
We consider the walk on the Cayley graph corresponding to this subset of generators.
In the natural realisation of this walk, each coordinate is updated at rate $1/d$; we want it to be updated at rate $1/k$.
If this walk is mixed, then since the walk on $G_k$ is obtained by a random independent shift of that walk, the walk on $G_k$ is also mixed.
Hence the previous entropic upper bound on the mixing time from the case $k = d$ is still valid, after multiplication by
\(
	k/d = 1 + \Oh{1/d} = 1 + \oh1,
\)
due to replacing the rate $1/d$ by $1/k$.

\smallskip

We now turn to the lower bound.
Set $\zeta \cq \tfrac1k(k - d) \log \pp \asymp 1/k \ll 1$.
Since $k > d$, we can apply our argument from \S\ref{sec-p5:p:lower}.
To bound the variance, we can no longer assume that $\abs{ \zeta_\alpha - \zeta_0 } \le \tfrac12 \zeta_0 = \tfrac12 \zeta$, since $\zeta k \asymp 1$.
(Recall the definition of $\zeta_\alpha$ from \cref{def-p5:p:entropic-time}.)
This means that there is an extra factor of $\zeta_\alpha / \zeta_0 = 1 - 2 \alpha / \sqrt{\zeta k} \asymp 1$ multiplying the variance.
For the lower bound, we need only consider $\alpha < 0$ with $\abs \alpha$ large.
Multiplying the variance by a constant only changes the window by a constant; in particular it does not affect the occurrence of cutoff.
Hence the entropic time lower bound is still valid.
Also, as $k/d = 1 + \oh1$, multiplying it by $k/d$ does not affect the leading order.

\smallskip

Finally we must asymptotically evaluate this entropic time when $0 < (k - d)\pp = \Oh1$.
Up to subleading order terms, by Proposition~\ref{res-p5:p:t0a} it equals the desired time:
\[
	\tfrac12 k \log\rbb{ k (k - d)^{-1} / \log \pp } / \rbb{ 1 - \cos(2\pi/\pp) }
\eqsim
	\tfrac12 d \log d / \rbb{ 1 - \cos(2\pi/\pp) },
\qedhere
\]
\end{Proof}

\section{Cutoff: No Cutoff When $k$ Is Constant}
\label{sec-p5:const-k}

Throughout the paper we have always been assuming that $\kinf$ as $\ninf$.
It is natural to ask what happens when $k$ does not diverge.
This case has actually already been covered by \textcite{DSc:growth-rw}, using their concept of \textit{moderate growth}.
Here we give a short exposition of their results leading to the conclusion that, for nilpotent groups of bounded step, there is no cutoff---for any choice of generating set, not only when one draws the Cayley graph uniformly.

Recall that a group $G$ is called \textit{nilpotent of step at most $L$} if its lower central series terminates in the trivial group after at most $L$ steps:
	$G_0 \cq G$ and $G_\ell \cq [G_{\ell-1}, G]$ for $\ell \in \mbn$ with $G_L = \bra{\id}$.

\begin{defn}[\cite{DSc:growth-rw}]
	Let $G$ be a finite group.
	Let $Z$ be a symmetric generating subset; that is,
	\(
		\bra{ z_1 \cdots z_r \mid r \in \mbn_0, \: z_1, ..., z_r \in Z } = G
	\)
	and if $z \in Z$ then $z^{-1} \in Z$ also.
	For $R \in \mbn_0$, let $\mcb(R)$ denote the $R$-ball around the identity in $\mcg$.
	Write $\Delta \cq \inf\bra{ R \in \mbn_0 \mid \abs{\mcb(R)} = \abs G }$ for the \textit{diameter} of $G(Z)$.
	We say that $G(Z)$ is of \textit{$(A,d)$-moderate growth} if
	\(
		\abs{\mcb(R)} \ge A^{-1} \abs G (R / \Delta)^d
	\)
	for all $R \in \mbn_0$.
\end{defn}

The main abstract result of \textcite{DSc:growth-rw} considers simple random walks on general Cayley graphs of moderate growth; see \cite[Theorem~3.1]{DSc:growth-rw} for a slight extension, considering more general random walks, which fundamentally gives the same conclusion.

\begin{thm}[{\cite[Theorem~1.2]{DSc:growth-rw}}]
	Let $G(Z)$ be a Cayley graph of $(A,d)$-moderate growth; write $\Delta \cq \diam G(Z)$.
	For $t \in \mbn_0$, let $d_\TV(t)$ denote the TV distance between the law of the lazy SRW run for $t$ steps and the uniform distribution.
	Let $c > 0$.
	Then the following hold:
	\[
		d_\TV\rbb{ 2 (1 + c) \abs Z \Delta^2 }
	&\le
		B e^{-c}
	\Qwhere
		B \cq A^{1/2} 2^{d(d+3)/4};
	\\
		d_\TV\rbb{ c \Delta^2 / (2^{4d+1} A^2) }
	&\ge
		\tfrac12 e^{-c}.
	\]
	Further, the corresponding \textit{relaxation time} $\trel^*$ satisfies $\trel^* \ge \Delta^2 / (4^{2d+1} A^2)$.
\end{thm}

The claim on the spectral is not included in the statements of \textcite{DSc:growth-rw}. However, the lower bound is proved precisely via the standard eigenvalue analysis; \cite[(3.2)]{DSc:growth-rw} gives the required inequality (in the notation there $\beta_1$ is the largest non-trivial eigenvalue).

\citeauthor{DSc:growth-rw} then make the following observation, formalised below.

\begin{cor}
	Let $A, d > 0$.
	Let $(G_N(Z_{(N)}))_\Ninn$ be a sequence of finite, undirected Cayley graphs of $(A,d)$-moderate growth and with $\sup_N \abs{Z_{(N)}} < \infty$.
	Then the corresponding sequences of lazy simple random walks does not exhibit the cutoff phenomenon; in fact,
	\[
		\tmix\rbb{ G_N(Z_{(N)}) } / k_N
	\lesssim
		\rbb{ \diam G_N(Z_{(N)}) }^2
	\lesssim
		\trel^*\rbb{ G_N(Z_{(N)}) }
	\lesssim
		\tmix\rbb{ G_N(Z_{(N)}) }
	\quad
		\asinf N.
	\]
\end{cor}

\citeauthor{DSc:growth-rw} apply this to nilpotent groups of bounded step.

\begin{thm}[{\cite[Lemma~5.1 and Theorem~5.2]{DSc:growth-rw}}]
	Let $G$ be a nilpotent group of step $L$.
	Let $Z$ be a symmetric set of generators for $G$.
	Then $G(Z)$ is of $(A, \log_2 A)$-moderate growth for some $A \cq A(\abs Z, L)$, depending only on the number of generators $\abs Z$ and the step $L$.
\end{thm}

As a corollary of this, if the number of generators is bounded and the underlying group is nilpotent of bounded step, then the corresponding simple random walks do not exhibit cutoff.

For a Cayley graph $G(Z)$, use the following notation.
	Write $\Delta \cq \diam G(Z)$ for its diameter.
	For the lazy simple random walk on $G(Z)$,
	write
		$\trel \cq \trel\rbr{ G(Z) }$ for the relaxation time (ie inverse of the spectral gap)
	and
		$\tmix \cq \tmix\rbr{ \eps; G(Z) }$ for the (TV) $\eps$-mixing time, for $\eps \in (0,1)$.
When considering sequences $(G_N(Z_{(N)}))_\Ninn$, add an $N$-sub/superscript.

\begin{cor}[{cf \cite[Corollary~5.3]{DSc:growth-rw}}]
\label{res-p5:const-k:main}
	Let $(G_N)_\Ninn$ be a sequence of finite, nilpotent groups.
	For each $\Ninn$, let $Z_{(N)}$ be a symmetric generating set for $G_N$ and write $L_N$ for the step of $G_N$.
	Suppose that $\sup_N \abs{Z_{(N)}} < \infty$ and $\sup_N L_N < \infty$.
	Then
	\(
		\tmix^N / k_N
	\lesssim
		\Delta_N^2
	\lesssim
		\trel^N
	\lesssim
		\tmix^N
	\)
	\asinf N;
	in particular, $(\tmix^N)_\Ninn$ does not exhibit the cutoff phenomenon
\end{cor}

\section{Cutoff: Difficulties Arising for High-Dimensional Matrices}
\label{sec-p5:upd:p-small-k-close-to-d}

This section is an extension of \cite[\S\ref{sec-p1:conc-rmks:small}]{HOt:rcg:matrix}.
There we consider cutoff for the RW on $G_k$ with $G \in \bra{\ugr, \hgr}$ with either $p$ small (eg $p \asymp 1$) or $k$ close to $d(\gab)$.
It is not intended to be read independently of the referenced section.
In particular, we use terminology and notation with which a reader of \cite{HOt:rcg:matrix} should be familiar by that point of the article.

\renewcommand{\ugr}{\UU_{p,d}}
\renewcommand{\uab}{\UU_{p,d}^\ab}
\renewcommand{\uk}{(\UU_{p,d})_k}
\renewcommand{\hgr}{\HH_{p,d}}
\renewcommand{\hab}{\HH_{p,d}^\ab}
\renewcommand{\hk}{(\HH_{p,d})_k}

\medskip

Denote by $S^\ab(t)$ the location $S(t)$ of the RW projected to the Abelianisation $\gab$.
Since $\gab$ is Abelian, $S^\ab(t)$ depends only on the final count $W(t)$, not the order in which the generators are chosen.
Now, $\gab = \mbz_p^r$, where $r \in \bra{d-1, 2d-4}$ (depending on whether $G = \UU$ or $G = \HH$). In particular, every element is of order $p$. Thus we need only know $W(t)$ with each coordinate taken mod $p$.
Defining $\widebar W$ to be a rate-1 RW on $\mbz_p^k$ (while $W$ is a RW on $\mbz^k$), an identical lower bound to that given in \S\ref{sec-p1:cutoff:lower} shows that one must wait for $\widebar W$ to have entropy $\log \abs \gab$ (while before we waited for $W$ to have this entropy).
If $p$ is small, eg $p = 2$, then the entropy of $\widebar W$ may grow significantly more slowly than that of $W$.
Even when $k - d(\gab) \asymp k$, the entropic time for $\widebar W$ can be a constant factor larger than that for $W$.
When $k - d(\gab) \ll k$, it can be of larger order.
Thus it is not sufficient to simply work with $W$; one needs $\widebar W$.
(In \cite[\cref{res-p2:intro:tv}]{HOt:rcg:abe:cutoff}, we extend this concept further and establish cutoff for all Abelian groups $A$ when $k - d(A) \gg 1$.)

If $p \gg 1$ (diverging arbitrarily slowly) and $k - d(\gab) \asymp k$, then one can check that the entropic times agree asymptotically.
\emph{Potentially}, one can then with $W$ rather than $\widebar W$.
This is could very beneficial for the Heisenberg group, since we can already handle any $\HH_{m,d}$ with $m$ fixed (not necessarily prime) and $k - d(\hab) \gg 1$.
For the upper triangular group, however, we do not have such a result.
However, it is distinctly possible that there are cases where the entropic times agree asymptotically and yet one should still work with $\widebar W$, rather than $W$.

\smallskip

Another obstacle, potentially more substantial, is the following.
As part of our typicality conditions, we ensured that $\abs{C_{i,j}} < p$, at least for many pairs.
This was crucial in controlling the error probability $q(t)$.
If $p \asymp 1$, then this is not actually an issue.
Indeed, when $k \gtrsim \log \abs \gab$, we found that in fact at least a constant proportion of generators are picked precisely once.
In this case, sufficiently many $(i,j)$ satisfy $\abs{C_{i,j}} = 1$.
We used this when extending from $p$ prime to general $m$; see \S\ref{sec-p1:ext:comp}.
If $p \asymp 1$, then $k \ge d(\hab) = 2d-4 \asymp \log \abs \hab$, so we always fall into this set-up.

For the Heisenberg group $\hgr$, we already handled $p \asymp 1$.
Unfortunately, another issue arises for the upper triangular group $\ugr$.
We partitioned the generators according to the columns, so as to get some desired independence:
	see the analysis leading up to \cref{eq-p1:cutoff:d:decomp};
	key in this is the factor $2$ in \cref{eq-p1:cutoff:d:Sab=S'ab}.
If we could remove this factor $2$ (or replace it with $1 + \oh{1/d}$), then we believe that the analysis from \S\ref{sec-p1:cutoff:d} would pass through to handle $p \asymp 1$. The only change would be to restrict to those generators picked at most once, exactly as we did in \S\ref{sec-p1:ext:comp}.
However, we do not know how to replace $2$ with $1 + \oh1$, never mind $1 + \oh{1/d}$.
(It seems plausible that $2/p$ could be replaced with $1/p + 1/p^2$; however, for $p \asymp 1$, this is not sufficient.)
It thus appears that if one can do something more creative than simply partitioning the generators, then the case of $p \asymp 1$---even $\UU_{m,d}$ with $m \asymp 1$ not necessarily prime---could be handled.
Such an argument eludes us---for now at least.

Recall from \cref{res-p1:cutoff:d:matrix-product} the expression for a product of upper triangular matrices. We only analysed the first (corresponding to the Abelianisation) and second (corresponding to $C_{i,j}$) order terms, discarding the remaining `higher order terms'.
For $d = 3$, these non-existent.
There are more and more of them, however, as $d$ grows.
Discarding them may be too crude when $p \asymp 1$.

For any $p \gg 1$, the regime $1 \ll k \ll \log \abs \gab$ is non-empty.
If $p$ diverges sufficiently slowly and, eg, $r > \tfrac12 k$, then the assumption does not hold.
Indeed, $\ex{\abs{W_1}} \asymp t_0(k, p^r)/k \asymp p^{2r/k} \ge p$ if $k - d(\gab) \asymp k$.
Thus is cannot be the case that $\ex{\abs{C_{1,2}}} \ll p$.
This is a serious problem for our method of controlling $q(t)$ in this regime---and it exists whichever of $W$ or $\widebar W$ is used.

\smallskip

Multiple of the above suggested solutions require $k - d(\gab) \asymp k$.
When this does not hold, things can be significantly more difficult.
Using the RW $\widebar W$ on $\mbz_p^k$, rather than the RW $W$ on $\mbz^k$, seems unavoidable. Indeed, the entropic time for the former will, if $k - d(\gab) \ll k / \log p$, be of larger order for $\widebar W$ than for $W$.
This is perhaps not such a substantial problem though: it is clear what the correction is (ie replace $W$ with $\widebar W$).
More substantial is that in many cases one cannot assume that $\abs{C_{i,j}} < p$.
One needs to study the `wrap around' effect of taking $C_{i,j}$ mod $p$.
If one tries to replace prime $p$ by general $m$, the gcd analysis may get even more complicated.

\section{Cutoff: From Upper Triangular Matrices to Nilpotent Groups}
\label{sec-p5:heis-nil}

Most of the following discussion is based on observations made by P\'eter Varj\'u during discussions of our work with him.
%
A group is \textit{nilpotent of step at most $\ell$} if all iterated commutators of order at least $\ell+1$ vanish necessarily.
For example, step-1 is Abelian; step-2 has $[[g_1, g_2], g_3] = \id$ for all $g_1$, $g_2$ and $g_3$, ie the commutator subgroup is central.
Our analysis has focussed on unit upper triangular matrix groups; these are a canonical class of nilpotent groups.
However, some of our analysis does extend somewhat to more general nilpotent.

\smallskip

Recall that we wrote $S$ for the location of the walk and $W$ for its auxiliary variable; let $W'$ be an independent copy of $W$, and define $S'$ correspondingly.
As previously, we work in the directed regime; so in the word $S$ there are no inverses.
Recall the definition of $C_{i,j}$ from \cite[\cref{eq-p1:cutoff:3:Cij-def}]{HOt:rcg:matrix}:
\[
	C_{i,j}
\cq
	\sumt[N]{\ell=1}
	\oneb{ G_\ell = j }
	\sumt[\ell-1]{m=1}
	\oneb{ G_m = i }
\Qand
	C_{i,i}
\cq
	0
\Quad{for all}
	i,j \in [k],
\]
where there are $N$ steps and $G_m$ is the index of the generator chosen in step $m$.

\begin{lem}
	Up to multiplication by an element of $[G, G^\com]$, we can express $S$ as
	\[
		S
	=
		\rbb{ \prodt[k]{1} Z_i^{W_i} }
	\cdot
		\rbb{ \prodt{i < j} [Z_i^{-1}, Z_j^{-1}]^{-C_{i,j}} }
	\]
	If $G$ is step-2 nilpotent then $[G, G^\com] = \bra{\id}$ is the trivial group.
	
	(The second product is unordered, since we are working up to an element of $[G, G^\com]$, and so we may assume that commutators commute with any element of $G$; the first is ordered $i = 1, ..., k$.)
\end{lem}

\begin{Proof}[Sketch of Proof]
Writing a rigorous proof of this lemma is technical, and can obscure what is going on; we use an example to demonstrate how to prove the lemma.
In essence, we wish to move all the $Z_1$-s to the left, then all the $Z_2$-s to the left-but-one and so on.
To reverse the order terms, we use the fact that
\(
	h g = g h h^{-1} g^{-1} g h = g h [h^{-1}, g^{-1}]
\)
and
\(
	[h^{-1}, g^{-1}]
=
	[g^{-1}, h^{-1}]^{-1}.
\)
For example,
\[
	g h h g
=
	g h \cdot g h [g^{-1}, h^{-1}]^{-1}
=
	g \cdot g h [g^{-1}, h^{-1}]^{-1} \cdot h [g^{-1}, h^{-1}]^{-1}
=
	g^2 h^2 [g^{-1}, h^{-1}]^{-2}.
\]
To move $Z_i$ past $Z_j$, with $i < j$, for each occurrence of $Z_i$ we need to count the number of times that $Z_j$ appears before it in the word; this is precisely (the definition of) $C_{i,j}$.
\end{Proof}

Expressing $S^{-1} S'$ as a similar product, it is straightforward to see what we get when $W = W'$.
(We actually only need $W_i \equiv W'_i \mod \ord Z_i$ for each $i$, but $W = W'$ is generally easier to analyse.)

\begin{cor}
	If $W = W'$, then $S^{-1} S' \in [G, G] / [G, [G, G]]$.
	The converse holds in the free group, ie when considering $Z_1, ..., Z_k$ as formal variables (ie with no relations between them).
	
	If $W = W'$, then, up to multiplication by an element of $[G, G^\com]$, we can express $S^{-1} S'$ as
	\[
		S^{-1} S'
	=
		\prodt{i < j} [Z_i^{-1}, Z_j^{-1}]^{D_{i,j}}
	\Qwhere
		D_{i,j} \cq C_{i,j} - C'_{i,j};
	\Quad{write}
		D \cq (D_{i,j})_{i,j}.
	\]
	In particular, if $C_{i,j} = C'_{i,j}$ for all $i$ and $j$ (which implies that $W_i = W'_i$ for all $i$ by taking $i = j$), then $S^{-1} S' \in [G, G^\com]$; if the group is step-2 nilpotent, then $[G, G^\com] = \bra{\id}$, and hence $S = S'$.
\end{cor}

Consider now step-2 nilpotent groups, of which the $3 \times 3$ matrices is an example.
We are interested in analysing $\pr{S = S' \mid \typ}$; typicality will primarily involve entropic considerations.
For ease of presentation, here we drop $\typ$ from the notation.
As in \cite[\S\ref{sec-p1:cutoff}]{HOt:rcg:matrix}, we note that
\[
	\pr{ S = S' }
\le
	\pr{ S = S' \mid W = W' } \pr{ W = W' }
+	\pr{ S = S' \mid W \ne W' }.
\]
Typicality (entropy) bounds $\pr{W = W'} \ll 1/\abs{G^\ab} = \abs{G^\com}/\abs G$, as for upper triangular groups.

Assume that $t \ll k$, and that every generator is picked at most once---eg, this is the case if $k \gg \log n$.
The assumption means that some generator is picked once in $S$ and never in $S'$ (or vice versa); this will allow us to deduce that $S^{-1} S' \sim \Unif(G)$, and hence
\(
	\pr{S = S' \mid W \ne W'}
=
	1/n.
\)

Since $S = S'$ when $D_{i,j} = 0$ for all $i$ and $j$, we have
\[
	\pr{ S = S' \mid W = W' }
\le
	\pr{ S = S' \mid W = W', \, D \ne 0}
+	\pr{ D = 0 }.
\]
When the nilpotent group is of higher step, the bound
\(
	\pr{S = S' \mid D = 0} \le 1
\)
may be too crude.
We analysed $\pr{D = 0}$ in \cite[\S\ref{sec-p1:cutoff}]{HOt:rcg:matrix}, obtaining $\pr{D = 0} \approx 1/t!$.
We desire this to be close to $1/\abs{G^\com}$.

We wish to get
\(
	\pr{ S = S' \mid W = W', \, D \ne 0 }
\Quad{close to}
	1/\abs{ G^\com }.
\)
To do this, write
\[
	S^{-1} S'
=
	\prodt{i < j : D_{i,j} \ne 0} \sbb{ Z_i^{-1}, Z_{j'}^{-1} }{}^{D_{i,j}}.
\]
While these commutators are neither uniformly random nor independent, we aim to have suitably many $D_{i,j} \ne 0$ so that the commutator product is sufficiently close to uniform (on $G^\com$).

If ``close'' can mean ``up to a sufficiently small factor'', then combining all these bounds gives
\(
	\pr{S = S', \, W = W'} \ll 1/n.
\)
The modified $L_2$ distance is then given by
\(
	n \, \pr{S = S'} - 1 = \oh1.
\)

\medskip

We can apply the method for nilpotent groups of greater step, by quotienting out $[G, G^\com]$.
However, as the step increases the bounds become more crude:
	we could have
	\(
		\pr{S = S'}
	\ll
		\pr{S^{-1} S' \in [G, G^\com]},
	\)
	which would be bad for this method.
	This is, in essence, what we did for the $d \times d$ matrices; it is one of the main reasons why our analysis does not work well for large $d$.
The analysis also applies to non-nilpotent groups, for which such issues can be even worse.

\section{Typical Distance: Generalised Graph Distance}
\label{sec-p5:typdist}

This section focusses on distances from a fixed point in the uniform random Cayley graph of degree $k$ of an Abelian group $G$.
The analysis is very similar to that of \cite[\S\ref{sec-p3:typdist2}]{HOt:rcg:abe:geom} where the same statistic was studied; here we are more general.
In particular, there we only considered $k \asymp \log \abs G$.
Here we adapt that analysis to consider $1 \ll k \ll \log \abs G$; we also extend the concept of graph distance from an $L_1$-type concept to an $L_\qq$-type, for general $\qq \in [1,\infty]$.

\subsection{Definition of $L_\qq$ Typical Distance}
\label{sec-p5:typdist:def}

Graphs distances in Cayley graphs have some special properties.
Consider a collection $z = [z_1, ..., z_k]$ of generators and distances in the Cayley graph $G(z)$.
For a path $\rho$ in $G(z)$,
for each $i \in [k]$,
write
	$\rho_{i,+}$ for the number of times $z_i$ is used,
	$\rho_{i,-}$ for the number of times $z_i^{-1}$ is used (if in the undirected case otherwise $\rho_{i,-} \cq 0$)
and
	$\rho_i \cq \rho_{i,+} - \rho_{i,-}$.
The path connects the identity with $\rho \bcdot z$.
Then the length, in the usual graph distance, of $\rho$ is $\norm{\rho}_1 \cq \sumt[k]{1} \rbr{ \rho_{i,+} + \rho_{i,-} }$.

\smallskip

For any $\qq \in [1,\infty)$,
define the \textit{$L_\qq$ graph distance} of $\rho$ by
\(
	\norm{\rho}_\qq^\qq
\cq
	\sumt[k]{1} \rbr{ \rho_{i,+}^\qq + \rho_{i,-}^\qq }.
\)
For the \textit{$L_\infty$-graph distance},
define
\(
	\norm{\rho}_\infty
\cq
	\maxt{i} \rbr{ \rho_{i,+} + \rho_{i,-}}.
\)
(The usual graph distance is given by $\qq = 1$.)

For Abelian groups, clearly for any $\qq \in [1,\infty)$ an \textit{$L_\qq$ geodesic}, ie a path of minimal length, will only use either $z_i$ or $z_i^{-1}$, not both (since the terms in the product can be reordered), ie $\rho_{i,+} \rho_{i,-} = 0$ for all $i$.
Thus $\norm{\rho}_\qq^\qq = \sumt[k]{1} \abs{\rho_i}^\qq$.
Similarly, any $L_\infty$-geodesic $\rho$ can be adjusted into a new path $\rho'$ with $\norm{\rho}_\infty = \norm{\rho'}_\infty$ and $\rho'_{i,+} \rho'_{i,-} = 0$ for all $i$.

We define the \textit{$L_\qq$ typical distance} $\mcd_{G(z),\qq}(\cdot)$ analogously to $\mcd_{G(z)}(\cdot)$, ie the $\qq = 1$ case.
When the $k$ generators are chosen uniformly at random, we write $\mcd_{G_k,\qq}^\pm(\cdot)$, with the $\pm$-superscript indicating whether or not the Cayley graph is directed.

\subsection{Precise Statement}
\label{sec-p5:typdist:statement}

For an Abelian group, we define the \textit{dimension} and \textit{minimal side-length}, respectively, as follows:
\[
	d(G) &\cq \min\brb{ d \in \mbn \midb \oplus_1^d \: \mbz_{m_j} \text{ is a decomposition of } G };
\\
	m_* &\cq \max\brb{ \mint{j=1,...,d} m_j \midb \oplus_1^d \: \mbz_{m_j} \text{ is a decomposition of } G }.
\]
It can be shown that there exists an optimal decomposition $\bra{m_j}_1^d$ for $m_*$ with $d = d(G)$.
Our main constraints will be $\limsup d/k < 1$ and $k^{1/\qq} n^{1/k} / m^* \ll 1$.

\begin{hyp}
\label{hyp-p5:typdist}
	The sequence
		$(k_N, G_N)_\Ninn$
	and
		$\qq \in [1,\infty]$
	jointly satisfy \textit{\cref{hyp-p5:typdist}} if the following conditions hold (defining $k^{1/\infty} \cq 1$ for $k \in \mbn$):%
	\begin{gather*}
		\LIM{\Ninf} k_N = \infty,
	\quad
		\LIM{\Ninf} k_N / \log \abs{G_N} = 0
	\Qand
		\LIM{\Ninf} k_N^{1/\qq} \abs{G_N}^{1/k_N} / m_*(G_N) = 0;
	\\
		\text{if\quad $\qq \in (1,\infty)$\quad then additionally\quad $k_N \le \log \abs{G_N} / \log\log \abs{G_N}$ for all $\Ninn$};
	\\
		\LIMSUP{\Ninf} \frac{d_N}{k_N}
	<
		\begin{cases}
			1			&\text{for undirected graphs}, \\
			\tfrac12	&\text{for directed graphs}.
		\end{cases}
	\end{gather*}
\end{hyp}

Finally we set up a little more notation.
Make the following definitions:
\[
	C^-_\qq \cq 2 \, \Gamma(1/\qq+1) (\qq e)^{1/\qq},
\quad
	C_\qq^+ \cq \tfrac12 C^-_\qq
\Qand
	\mfd^\pm_\qq(k,n) \cq k^{1/\qq} n^{1/k} / C^\pm_\qq,
\]
where the case $\qq = \infty$ is to be interpreted as the limit $\qq \to \infty$;
eg,
\(
	C^-_\infty
=
	2
\)
and
\(
	\mfd^+_\infty(k,n)
=
	n^{1/k}.
\)
When these are sequences $(k_N, G_N)_\Ninn$,
for $\Ninn$ and $q \in [1,\infty]$,
write $\mfd^\pm_{N,\qq} \cq \mfd^\pm_\qq(k_N,\abs{G_N})$.

Similarly,
for a sequence $(G_N)_\Ninn$ of finite groups with corresponding multisubsets $(Z_{(N)})_\Ninn$ of sizes $(k_N)_\Ninn$,
for $\Ninn$, $\beta \in [0,1]$ and $\qq \in [1,\infty]$,
define
\(
	\mcd_{N,\qq}^\pm \cq \mcd_{G_N^\pm(Z_{(N)})}(\beta).
\)

\begin{thm}
\label{res-p5:typdist:res}
	Let $(k_N)_\Ninn$ be a sequence of positive integers and $(G_N)_\Ninn$ a sequence of finite, Abelian groups;
	for each $\Ninn$, define $Z_{(N)} \cq [Z_1, ..., Z_{k_N}]$ by drawing $Z_1, ..., Z_{k_N} \sim^\iid \Unif(G_N)$.
	
	Suppose that the sequence $(k_N, G_N)_\Ninn$ satisfies \cref{hyp-p5:typdist}.
	Then,
	for all $\beta \in (0,1)$,
	we~have
	\[
		\mcd^\pm_{N,\qq}(\beta) /\mfd^\pm_{N,\qq}
	\to^\mbp
		1
	\Quad{(in probability)}
		\asinf N.
	\]
	Moreover, the implicit lower bound holds for all choices of generators and for all Abelian groups, only requiring the conditions in \cref{hyp-p5:typdist} which depend only on $(k_N, \abs{G_N})_\Ninn$ and $q$.
\end{thm}

\begin{rmkt*}
	We initially prove this theorem for undirected Cayley graphs.
	In \S\ref{sec-p5:typdist:directed}, we explain how to adapt the proof from the undirected case to the directed case.
	Doing this, rather than making every statement apply for both the un- and directed cases, significantly increases the readability.
	In particular, when we speak of $\mbz$ we are referring to the set of all integers, positive and negative.
\end{rmkt*}

\begin{rmkt*}
	We use the same methodology as \cite[\S\ref{sec-p3:typdist2}]{HOt:rcg:abe:geom}.
	An outline of the proof is given in \cite[\S\ref{sec-p3:typdist2:outline}]{HOt:rcg:abe:geom}.
\end{rmkt*}

\subsection{Size of Ball Estimates and Lower Bound}
\label{sec-p5:typdist:size-lower-bound}

In the lemmas below, used to prove this theorem, instead of writing one lemma with multiple parts, we split into separate lemmas according to $\qq$ and $k$, eg $\qq \in (1,\infty)$ or $k \asymp \log n$; these parts are indexed with letters, eg Lemmas~\ref{res-p5:typdist:size:1}, \ref{res-p5:typdist:size:q} and \ref{res-p5:typdist:size:inf}.

\medskip

We wish to determine the size of the $L_\qq$ balls in $\mbr^k$.
This is done by Lemmas \ref{res-p5:typdist:size} and \ref{res-p5:typdist:M};
the statements are given below, with proofs are deferred to
the supplementary material, \cite[\S\ref{sec-p0:balls}]{HOt:rcg:supp}.

For $\qq \in [1,\infty)$, write $V_{k,\qq}(R)$ for the (Lebesgue) volume of the $L_\qq$ ball of radius $R$ in $\mbr^k$, ie
\[
	V_{k,\qq}(R) \cq \vol\brb{ x \in \mbr^k \midb \norm{x}_\qq \le R };
\]
also write $V_{k,\qq} \cq V_{k,\qq}(1)$ and note that $V_{k,\qq}(R) = R^k V_{k,\qq}$.
It is known (see \cite{W:Lp-balls}) that
\[
	V_{\ell,\qq} = 2^\ell \Gamma(1/\qq+1)^\ell / \Gamma(\ell/\qq+1).
\label{eq-p5:typdist:vol}
\nt
\]
We can use this, along with \cref{res-p5:typdist:size:q} below, to well-approximate $\abs{B_{k,\qq}(R)}$ when $\qq \notin \bra{1,\infty}$;
for $\qq = 1$ we directly bound $\abs{B_{k,1}(\cdot)}$, while for $\qq = \infty$ we have an exact expression.

\begin{subequations}
	\label{eq-p5:typdist:size}
\begin{subtheorem}{thm}
	\label{res-p5:typdist:size}
	
\begin{lem}
\label{res-p5:typdist:size:1}
	For $\qq = 1$ and all $R \ge 0$,
	we have
	\[
		2^k \binomt{\floor{R}}{k} \one{R \ge k}
	\le
		\absb{ B_{k,1}(R) }
	\le
		2^k \binomt{\floor{R}+k}{k}.
	\label{eq-p5:typdist:size:1}
	\nt
	\]
\end{lem}

\begin{lem}
\label{res-p5:typdist:size:q}
	For $\qq \in (1,\infty)$ and all $R \ge k^{1+1/\qq}$,
	we have
	\[
		\absb{ B_{k,\qq}(R) }
	=
		V_{k,\qq}(R) \, \rbb{ 1 + \Ohb{ k^{1+1/\qq}/R } }.
	\label{eq-p5:typdist:size:q}
	\nt
	\]
\end{lem}

\begin{lem}
\label{res-p5:typdist:size:inf}
	For $\qq = \infty$ and all $R \ge 0$,
	we have
	\[
		\absb{ B_{k,\infty}(R) }
	=
		\rbb{2\floor{R} + 1}^k.
	\label{eq-p5:typdist:size:inf}
	\nt
	\]
\end{lem}

\end{subtheorem}
\end{subequations}

We use this lemma to find an $M$ so that $\abs{B_{k,\qq}(M)} \approx n$.

\begin{defn}
\label{def-p5:typdist:M}
	Set $\omega \cq \max\bra{ \logk[2], k/n^{1/(2k)} }$, and choose $M_{k,\qq}$ to be the minimal integer satisfying $\abs{ B_{k,\qq}(M_{k,\qq}) } \ge n e^\omega$.
	Note that $\omega$ satisfies $1 \ll \omega \ll k$ if $k \ll \log n$.
\end{defn}

Recall that $\mfm_{k,\qq} = k^{1/\qq} n^{1/k} / C_\qq$, and that $C_\qq = 2 \, \Gamma(1/\qq + 1) (\qq e)^{1/\qq}$.
The next lemma shows that the difference between $M$ and $\mfm$ is only by subleading order terms.
Also, let $K$ be a constant, assumed to be as large as required, and let $\xi \cq 1 - e^{-K\omega/k}$.

\begin{subequations}
	\label{eq-p5:typdist:M}
\begin{subtheorem}{thm}
	\label{res-p5:typdist:M}

\begin{lem}
\label{res-p5:typdist:M:1}
	For $k \ll \log n$ and $\qq = 1$,
	we have
	\[
		M_{k,1} \le \ceilb{ \mfm_{k,1} (1 + \xi) }
	\Qand
		\absb{ B_{k,1}\rbb{ \mfm_{k,1} (1 - \xi) } } \ll n.
	\label{eq-p5:typdist:M:1}
	\nt
	\]
\end{lem}

\begin{lem}
\label{res-p5:typdist:M:q}
	For $k \le \log n/\log\log n$ and all $\qq \in [1,\infty)$,
	we have
	\[
		M_{k,\qq} \le \floorb{ \mfm_{k,\qq} (1 + \xi) }
	\Qand
		\absb{ B_{k,\qq}\rbb{ \mfm_{k,\qq} (1 - \xi) } } \ll n.
	\label{eq-p5:typdist:M:q}
	\nt
	\]
\end{lem}

\begin{lem}
\label{res-p5:typdist:M:inf}
	For $\qq = \infty$,
	we have
	\[
		M_{k,\infty}
	=
		\ceilb{ \tfrac12 n^{1/k} e^{\omega/k} - \tfrac12 }
	\Qand
		\absb{ B_{k,\infty}\rbb{ \mfm_{k,\infty} (1 - \xi) } } \ll n.
	\label{eq-p5:typdist:M:inf}
	\nt
	\]
	Moreover, if $k \ll \log n$ then $M_{k,\infty} \eqsim \mfm_{k,\infty}$.
\end{lem}

\end{subtheorem}
\end{subequations}

\subsection{Lower Bound on Typical Distance}

From this lemma, it is straightforward to deduce the lower bound in \cref{res-p5:typdist:res}.

\begin{Proof}[Proof of Lower Bound in \cref{res-p5:typdist:res}]
Observe that
\(
	\abs{ \mcb_{k,\qq}(M) }
\le
	\abs{ B_{k,\qq}(M) }.
\)
By Lemma \ref{res-p5:typdist:M}, the right-hand side is $\oh n$ when $M \cq M_{k,\qq} (1-\xi)$ when $k \ll \log n$.
Thus $\mcd_{G_k,\qq}(\beta) \ge M$ for all~$Z$.
\end{Proof}

\subsection{Upper Bound on Typical Distance}

The outline of this subsection follows closely that of \cite[\S\ref{sec-p3:typdist2:upper}]{HOt:rcg:abe:geom}.

%


\begin{prop}
\label{res-p5:typdist:l2-bound}
	Let $\qq \in [1,\infty]$.
	Suppose that $k \ll \log n$.
	If $\qq \in (1,\infty)$, then further restrict to $k \le \log n/\log\log n$.
	Suppose also that $\limsup d/k < 1$.
	Then
	\(
		\ex{ \norm{ \pr{ \WW \bcdot Z = \cdot \mid Z } - \pi_G }_2^2 }
	=
		\oh1.
	\)
\end{prop}

Once we prove these propositions, we have all we need to prove \cref{res-p5:typdist:res}.

\begin{Proof}[Proof of \cref{res-p5:typdist:res} Given Lemma~\ref{res-p5:typdist:M} and \cref{res-p5:typdist:l2-bound}]
If $\norm{\pr{A \bcdot Z = \cdot \mid Z} - \pi_G}_2 \le \eps$, then the support $\mcs$ of $A \bcdot Z$ satisfies $\pi_G(\mcs^c) \le \eps$.
Combining this with Lemma~\ref{res-p5:typdist:M} and \cref{res-p5:typdist:l2-bound}, we deduce the upper bound in \cref{res-p5:typdist:res}.
\end{Proof}

\begin{rmkt*}
\cref{res-p5:typdist:l2-bound} actually holds even if $\eta \cq 1 - d/k \downarrow 0$, provided it does so sufficiently slowly and $k/\log n$ is sufficiently small.
It turns out that $k/\log n \ll \eta$ and $\eta \gg 1/\sqrt k$ is sufficient; this allows $k$ very close to both $d$ and $\log n$.
\end{rmkt*}

Let $\WW,\WW' \sim^\iid \Unif\rbr{ B_{k,\qq}(M) }$, and let $\VV \cq \WW - \WW'$.
Then we have
\[
	\ex{ \normb{ \pr[G_k]{ \WW \bcdot Z \in \cdot } - \pi_G }_2^2 }
=
	\ex{ n \, \pr{ \VV \bcdot Z = 0 \mid Z } - 1 }
=
	n \, \pr{ \VV \bcdot Z = 0 } - 1.
\]
%

First, it is immediate to see that
\(
	\pr{ \WW = \WW' } = \abs{ B_{k,\qq}(M) }^{-1} \le n^{-1} e^{-\omega}.
\)
Analogously in \S\ref{sec-p5:profile}, the side-lengths $\bra{m_j}_1^d$ satisfy $\min_j m_j > 2M$.
Then we have
\[
	\mci
\cq
	\brb{ i \in [k] \midb \VV_i \not\equiv 0 \MOD m_j \: \forall \, j = 1,...,d }
=
	\brb{ i \in [k] \midb \WW_i \ne \WW'_i }.
\]

\begin{subequations}
	\label{eq-p5:typdist:mi=I}
\begin{subtheorem}{thm}
	\label{res-p5:typdist:mi=I}

\begin{lem}
\label{res-p5:typdist:mi=I:empty}
	For all $k$ and all $\qq$,
	we have
	\[
		\pr{ \mci = \emptyset }
	\le
		n^{-1} e^{-\omega}.
	\label{eq-p5:typdist:mi=I:empty}
	\nt
	\]
\end{lem}

\begin{lem}
\label{res-p5:typdist:mi=I:fin}
	Suppose that $k \ll \log n$ and $\qq \in [1,\infty)$.
	If $\qq \in (1,\infty)$, then restrict further to $k \le \log n/\log\log n$.
	Then, for all $I \subseteq [k]$, we have
	\[
		\pr{ \mci = I }
	\le
		e^{k(1/(e \qq) + \xi_\qq)} n^{-1 + \abs I/k}
	\label{eq-p5:typdist:mi=I:fin}
	\nt
	\]
	where $\xi_\qq \cq K_\qq \omega / k \ll 1$, for some constant $K_\qq$.
\end{lem}

\begin{lem}
\label{res-p5:typdist:mi=I:inf}
	For $\qq = \infty$,
	for all $I \subseteq [k]$,
	we have
	\[
		\pr{ \mci = I }
	\le
		e^{-\omega(1 - \abs I/k)} n^{-1 + \abs I/k}
	\le
		n^{-1 + \abs I/k}
	\label{eq-p5:typdist:mi=I:inf}
	\nt
	\]
\end{lem}

\end{subtheorem}
\end{subequations}

While we have been stating results for undirected graphs, Lemma~\ref{res-p5:typdist:mi=I} holds in the directed case too.
Contrastingly, the following lemma distinguishes between the directed and undirected graphs at one point.
A proof of the lemma can be found in \cite[\cref{res-p3:typdist2:gcd-ex}]{HOt:rcg:abe:geom} in the main article.
(There, while we studied both undirected and directed graphs, it was sufficient to use the worst-case bound for both; there we need the slightly more refined statement. The identical proof works.)
Define
\[
	\mfgcd
\cq
	\gcd\rbb{ V_1, ..., V_k, n }.
\]

\begin{lem}
\label{res-p5:typdist:gcd}
	For all $I \subseteq [k]$,
	we have
	\[
		n \, \pr{ V \bcdot Z = 0 \mid \mci = I }
	\le
		\exb{ \mfgcd^d \mid \mci = I }.
	\label{eq-p5:typdist:gcd:pr-ex}
	\nt
	\]
	Further, there exists a constant $C$ so that,
	for all $I \subseteq [k]$,
	we have
	\begin{subequations}
		\label{eq-p5:typdist:gcd:ex}
	\begin{empheq}[left = {%
		\exb{ \mfgcd^d \mid \mci = I }
	\le
		\empheqlbrace}]%
	{alignat=2}
		&C 2^d (2M)^{d - \abs I + 2}
		&&\Qwhere \abs I \le d+1;
	\label{eq-p5:typdist:gcd:I<d}
	\\
		&1 + 3 \cdot 2^{d - \abs I}
		&&\Qwhere \abs I \ge d+2 \quad \text{for undirected grahs},
	\label{eq-p5:typdist:gcd:I>d:-}
	\\
		&1 + 5 \cdot (\tfrac32)^{2d - \abs I}
		&&\Qwhere \abs I \ge d+2 \quad \text{for directed graphs}.
	\label{eq-p5:typdist:gcd:I>d:+}
	\end{empheq}
	\end{subequations}
\end{lem}

The idea behind \cref{eq-p5:typdist:gcd:pr-ex} is that linear combinations of independent uniform random variables are uniform on their support.
Writing $G = \oplus_1^d \: \mbz_{m_j}$, we obtain $V \bcdot Z \sim \Unif(\oplus_1^d \: \mfgcd_j \mbz_{m_j})$ where $\mfgcd_j \cq \gcd(V_1, ..., V_k, m_j) \le \mfgcd$.
For a rigorous argument, see \cite[\cref{res-p3:typdist2:unif-gcd}]{HOt:rcg:abe:geom} in the main article.

When $\abs I$ is large, if $\mfgcd > 1$ then we are asking that a large number of coordinates have a common divisor; naturally this decays exponentially in $\abs I$.
Using this decay, we can sum over all ``large $I$''.

\begin{rmkt}
\label{rmk-p5:typdist:gcd:+vs-}
	We firmly believe that the stronger \cref{eq-p5:typdist:gcd:I>d:+} should hold for both the undirected and directed graphs (ie \cref{eq-p5:typdist:gcd:I>d:-} is unnecessary).
	It is merely a technical hurdle which is holding us back from proving this.
	When $\qq = \infty$, the coordinates of $V$ are independent; in this case, we can prove that \cref{eq-p5:typdist:gcd:I>d:-} holds for both the undirected and directed graphs.
	As a result, we can relax $\limsup d/k < \tfrac12$ to $\limsup d/k < 1$ for $\qq = \infty$.
\end{rmkt}

\begin{cor}
\label{res-p5:typdist:sum>L}
	For any $L$ with $L \ge d+2$, we have
	\begin{subequations}
		\label{eq-p5:typdist:sum>L}
	\begin{empheq}[left = {%
		n \sumt{\abs I \ge L}
		\pr{ V \bcdot Z = 0, \: \mci = I }
	\le
		\empheqlbrace}]%
	{alignat=2}
		&1 + 3 \cdot 2^{d - L}
		&&\quad \text{for undirected graphs}
	\\
		&1 + 5 \cdot (\tfrac32)^{2d - L}
		&&\quad \text{for directed graphs}
	\label{eq-p5:profile:sum>L:+}
	\end{empheq}
\end{subequations}
\end{cor}

\begin{Proof}
	This proof is a direct calculation.
	By \cref{eq-p5:typdist:gcd:pr-ex,eq-p5:typdist:gcd:I>d:-}, using Bayes's rule,
	specifically the fact that $\pr{B \mid C} / \pr{C \mid B} = \pr{B} / \pr{C}$ for non-null events $B$ and $C$,
	for $L \ge d + 2$ we deduce that
	\[
	&	n \sumt{\abs I \ge L}
		\pr{ V \bcdot Z = 0, \: \mci = I \mid \typ }
	=
		n \sumt{\abs I \ge L}
		\pr{ V \bcdot Z = 0 \mid \mci = I, \: \typ }
		\pr{ \mci = I \mid \typ }
	\\&\qquad
	\le
		\sumt{\abs I \ge L}
		\rbb{
			\pr{ \mci = I \mid \typ } + 3 \cdot 2^{d - \abs I} \, \pr{ \mci = I } / \pr{\typ}
		}
	\\&\qquad
	\le
		\pr{ \abs \mci \ge L \mid \typ }
	+	3 \cdot 2^{d - L} \pr{ \abs \mci \ge L } / \pr{\typ}
	\le
		1 + 3 \cdot 2^{d - L} / \pr{\typ}
	\]
	for undirected graphs.
	The case of directed graphs follows analogously.
\end{Proof}

We first prove the results on $\pr{\mci = I}$.
For a set $I \subseteq [k]$ and $\WW \in \mbz^k$, write $\WW_I = (\WW_i)_{i \in I}$ and $\WW_{\setminus I} = \WW_{I^c}$.
Recall that if $\mcc \subseteq \mcc'$ and $U \sim \Unif(\mcc')$, then $(U \mid U \in \mcc) \sim \Unif(\mcc)$.
Hence we have
\[
	\pr{ \WW_{\setminus I} = \WW'_{\setminus I} }
=
	\frac{ \pr{ \WW = \WW' } }{ \pr{ \WW_I = \WW'_I \mid \WW_{\setminus I} = \WW'_{\setminus I} } }
=
	\frac{ \abs{ B_{k,\qq}(M) }^{-1} }{ \ex{ \abs{ B_{\abs I,\qq}\rbr{ M - \norm{A_{\setminus I}}_1 } }^{-1} } }
\le
	\frac{ \abs{ B_{\abs I,\qq}(M) } }{ \abs{ B_{k,\qq}(M) } }.
\label{eq-p5:typdist:WnotI=W'notI}
\nt
\]
Write $\ell \cq \abs I$.
Recall that, by choice of $M$, we have $\abs{ B_{k,\qq}(M) } \ge n e^\omega$, and so
\[
	\pr{ \WW_{\setminus I} = \WW'_{\setminus I} }
\le
	n^{-1} e^{-\omega} \absb{ B_{\ell,\qq}(M) }.
\]

\begin{Proof}[Proof of \cref{res-p5:typdist:mi=I:empty}]
Recall the choice of $M_{k,\qq}$, from \cref{def-p5:typdist:M}. Then \cref{eq-p5:typdist:mi=I:empty} follows:
\[
	\pr{ \mci = \emptyset }
=
	\pr{ \WW = \WW' }
=
	\absb{ B_{k,\qq}(M_{k,\qq}) }^{-1}
\le
	n^{-1} e^{-\omega}.
\qedhere
\]
\end{Proof}

\begin{Proof}[Proof of \cref{res-p5:typdist:mi=I:fin}]
Consider first $\qq = 1$.	
From \cref{res-p5:typdist:M:1}, recall that $M_1 \le (2e)^{-1} k n^{1/k} e^\xi$ with $\xi \asymp \omega/k$.
Using \cref{res-p5:typdist:size:1}, for $\ell \le k$, we have
\[
	\absb{ B_{\ell,1}(M_1) }
\le
	2^\ell \binomt{M_1+\ell}{\ell}
\le
	\rbb{ 2e(M_1/\ell + 1) }^\ell
\le
	e^{\xi \ell} (k/\ell)^\ell n^{\ell/k}
\le
	e^{k(1/e+\xi)} n^{\ell/k},
\]
using the fact
	that $\binom N\ell \le (eN/\ell)^\ell$,
	that $\ell \mapsto (k/\ell)^\ell$ is maximised by $\ell = k/e$
and
	that $1 + x \le e^x$.
The proof is completed by noting that
\(
	\bra{ \mci = I } \subseteq \bra{ A_{\setminus I} = A'_{\setminus I} },
\)
and applying \cref{eq-p5:typdist:WnotI=W'notI}.

\smallskip

Now consider $\qq \in (1,\infty)$.
Justified by \cref{res-p5:typdist:size:q} and \cref{res-p5:typdist:M:q}, which shows that $M_{k,\qq} \gg k^{1+1/\qq}$ for all $\qq$, we replace this discrete ball by the continuous ball, and lose only a factor $1 + \oh1$; for readability, we do not carry this factor in subsequent formulae.

Using Stirling's formula and the upper for $M_{k,\qq}$ from \cref{res-p5:typdist:M:q} gives
\[
	V_{\ell,\qq}(M_{k,\qq})
\le
	V_{\ell,\qq} \cdot \rbb{ (1+\xi) k^{1/\qq} n^{1/k} / C_\qq }^\ell
\le
	\qq^{1/2} e^{K_\qq\omega} (k/\ell)^{\ell/\qq} n^{\ell/k}.
\]
From this, similarly to in \cref{res-p5:typdist:mi=I:empty}, using \cref{eq-p5:typdist:WnotI=W'notI}, we deduce that
\[
	\pr{ \WW_{\setminus I} = \WW'_{\setminus I} }
\le
	\qq^{1/2} e^{K_\qq \omega} (k/\ell)^{\ell/\qq} n^{-1+\ell/k}
\le
	e^{k(1/(e\qq)+\xi)} n^{-1+\ell/\qq},
\]
where $\xi \cq K_\qq \omega / k \ll 1$,
using again the fact
	that $\binom N\ell \le (eN/\ell)^\ell$
and
	that $\ell \mapsto (k/\ell)^\ell$ is maximised by $\ell = k/e$
The proof is completed by noting that
\(
	\bra{ \mci = I } \subseteq \bra{ \WW_{\setminus I} = \WW'_{\setminus I} }.
\)
\end{Proof}

\begin{Proof}[Proof of \cref{res-p5:typdist:mi=I:inf}]
The coordinates of $\WW$ satisfy $\WW_i \sim^\iid \Unif( \bra{0,\pm1, ..., \pm M_\infty } )$, for $i = 1, ..., k$.
Write $\ell \cq \abs I$.
Hence, by \cref{eq-p5:typdist:WnotI=W'notI} and \cref{eq-p5:typdist:mi=I:empty}, we have
\[
	\pr{ \WW_{\setminus I} = \WW'_{\setminus I} }
\le
	\absb{ B_{\ell,\infty}(M_\infty) } / \absb{ B_{k,\infty}(M_\infty) }
=
	(2M_\infty + 1)^{\ell-k}.
\]
By \cref{eq-p5:typdist:M:inf}, we have $2M_\infty+1 \ge n^{1/k}e^{\omega/k}$. Hence
\[
	\pr{ \mci = I }
\le
	\pr{ A_{\setminus I} = A'_{\setminus I} }
\le
	e^{\omega(-1 + \ell/k)} n^{-1+\ell/k}.
\qedhere
\]
\end{Proof}


We have now done all the hard work in proving \cref{res-p5:typdist:l2-bound}, from which we deduced \cref{res-p5:typdist:res}.
It remains to go through the details of how to combine the previous results; there are no more interesting ideas to prove the propositions, but the details are quite technical.


Similarly to the mixing proof, we use an $L_2$ calculation:
\[
	\ex{ \normb{ \pr{ \WW \bcdot Z = \cdot \mid Z } - \pi_G }_2^2 }
=
	n \sumt{I \subseteq [k]} \pr{ \VV \bcdot Z = 0, \: \mci = I } - 1.
\label{eq-p5:typdist:l2}
\nt
\]

\begin{Proof}[Proof of \cref{res-p5:typdist:l2-bound} \mdseries\unboldmath (when $\qq < \infty$)]
Recall that here $k \ll \log n$.
For undirected graphs, we have $\limsup d/k < 1$; for directed, $\limsup d/k < \tfrac12$.
Set $\eta^- \cq 1 - \limsup d/k > 0$ and $\eta^+ \cq \tfrac12 - \limsup d/k$;
set $L^- \cq d + \tfrac14 \eta^- (k - d)$ and $L^+ \cq 2d + \tfrac14 \eta^+ (k - 2d)$.
Use $L^-$ for undirected graphs and $L^+$ for directed.
Then
\[
	\limsup L^\pm/k
\le
	\tfrac14 \eta^\pm + (1 - \tfrac14 \eta^\pm) (1 - \eta^\pm)
\le
	1 - \tfrac23 \eta^\pm
<
	1;
\]
also, $L^- - d \gg 1$ and $L^+ - 2d \gg 1$.
Suppress the $\pm$-superscript: write $L \cq L^\pm$.
By Lemma~\ref{res-p5:typdist:M},
recalling that
\(
	C_\qq = 2 \, \Gamma(1/\qq+1) \, (\qq e)^{1/\qq},
\)
for some
\(
	\eps_\qq = \Oh{ \omega/k } = \oh1,
\)
we can write
\[
	M = (1 + \eps_\qq) k^{1/\qq} n^{1/k} / C_\qq.
\]
It can be shown that $C_\qq \ge 2$ for all $\qq \in [1,\infty]$, and hence
\[
	2M \le e^{\eps_\qq} k^{1/\qq} n^{1/k}.
\label{eq-p5:typdist:2M}
\nt
\]
Recall that when we consider $\qq = 1$, we only require $k \ll \log n$; when we consider $\qq \in (1,\infty)$, we require further that $k \le \log n/\log\log n$.
Note that if $\mci = \emptyset$ then $B = 0$, and so $B \bcdot Z = 0$. Hence
\[
	n \, \pr{ B \bcdot Z = 0 \mid \mci = \emptyset }
=
	n \, \pr{ \mci = \emptyset }
\le
	e^{-\omega},
\]
by the choice of the radius $M_{k,\qq}$.

Consider $I \subseteq [k]$ with $1 \le \ell = \abs I \le d+1$.
There are at most $2^k$ such sets $I$.
Recall $\xi_\qq$ given in \cref{res-p5:typdist:mi=I:fin}, and that $\xi_\qq = \Oh{\omega/k} = \oh1$.
Applying
(\ref{eq-p5:typdist:mi=I:fin}, \ref{eq-p5:typdist:gcd:pr-ex}, \ref{eq-p5:typdist:gcd:I<d}, \ref{eq-p5:typdist:2M}),
we obtain
\[
	n \, \pr{ \VV \bcdot Z = 0, \: \mci = I }
\le
	C 2^d e^{k \eps_\qq} e^{k \xi_\qq} k^{(d+2-\ell)/\qq} n^{(d+2-\ell)/k} \cdot e^{k/(e\qq)} n^{-1+\ell/k};
\]
algebraic manipulations using the fact that $\limsup d/k < 1$ and $2^d = n^{\oh1}$ then give
\[
	n \, \pr{ \VV \bcdot Z = 0, \: \mci = I }
=
	2^{-k} \oh1.
\label{eq-p5:typdist:l2mixing:proof:<:fin:small1}
\nt
\]
Consider $I \subseteq [k]$ with $d+2 \le \ell = \abs I \le L = d$.
Applying
(\ref{eq-p5:typdist:mi=I:fin}, \ref{eq-p5:typdist:gcd:pr-ex}, \ref{eq-p5:typdist:gcd:I>d:-}, \ref{eq-p5:typdist:gcd:I>d:+}),
we obtain
\[
	n \, \pr{ \VV \bcdot Z = 0, \: \mci = I }
\le
	15 \cdot e^{k(1/(e\qq) + \xi_\qq)} n^{-1+\ell/k};
\]
similar algebraic manipulations to those used when $1 \le \abs I \le d+1$ give
\[
	n \, \pr{ \VV \bcdot Z = 0, \: \mci = I }
=
	2^{-k} \oh1.
\label{eq-p5:typdist:l2mixing:proof:<:fin:small2.2}
\nt
\]
We now sum over all $I$ with $1 \le \abs I \le L$, using \cref{eq-p5:typdist:l2mixing:proof:<:fin:small1,eq-p5:typdist:l2mixing:proof:<:fin:small2.2}:
\[
	n \sumt{1 \le \abs I \le L}
	\pr{ \VV \bcdot Z = 0, \: \mci = I }
=
	\oh1.
\label{eq-p5:typdist:l2mixing:proof:<:fin:sum-small}
\nt
\]

Finally we consider $I \subseteq [k]$ with $L \le \abs I \le k$.
By \cref{res-p5:typdist:sum>L},
we have
\begin{subequations}
	\label{eq-p5:typdist:l2mixing:proof:sum-large:q<inf}
\begin{empheq}[left = {%
	n \sumd{L \le \abs I \le k}
	\pr{ V \bcdot Z = 0, \: \mci = I }
\le
	\empheqlbrace}]%
{alignat=2}
	&1 + 3 \cdot 2^{d - L} = 1 + \oh1
	&&\quad \text{for undirected graphs},
\label{eq-p5:typdist:l2mixing:proof:sum-large:q<inf:-}
\\
	&1 + 5 \cdot (\tfrac32)^{2d - L} = 1 + \oh1
	&&\quad \text{for directed graphs}.
\label{eq-p5:typdist:l2mixing:proof:sum-large:q<inf:+}
\end{empheq}
\end{subequations}
This last result actually holds for all $\qq \in [1,\infty]$ and all $1 \ll k \ll \log n$.

The proof is completed by combining \cref{eq-p5:typdist:l2mixing:proof:<:fin:sum-small,eq-p5:typdist:l2mixing:proof:sum-large:q<inf} with \cref{eq-p5:typdist:l2}.
\end{Proof}

\begin{Proof}[Proof of \cref{res-p5:typdist:l2-bound} \mdseries\unboldmath (when $\qq = \infty$)]
Recall that here $k \ll \log n$.
As discussed in \cref{rmk-p5:typdist:gcd:+vs-}, here \cref{eq-p5:typdist:gcd:I>d:+} holds for both the undirected and directed balls (ie \cref{eq-p5:typdist:gcd:I>d:-} is unnecessary); we only assume that $\limsup d/k < 1$ in either case.
Set $\eta \cq 1 - \limsup d/k > 0$.

By \cref{eq-p5:typdist:M:inf}, we have $2 M_{k,\infty} \le n^{1/k} e^{\omega/k} + 1$.
Consider $I \subseteq [k]$ with $1 \le \ell = \abs I \le d+1$.
There are at most $2^k$ such sets $I$.
Applying
(\ref{eq-p5:typdist:M:inf}, \ref{eq-p5:typdist:mi=I:inf}, \ref{eq-p5:typdist:gcd:pr-ex}, \ref{eq-p5:typdist:gcd:I<d}),
we obtain
\[
	n \, \pr{ \VV \bcdot Z = 0, \: \mci = I }
\le
	C 2^d n^{(d-\ell+2)/k} e^{\omega(d + 2 - \ell)/k} (1 + e^{-\omega/k}/n^{1/k})^{d+2-\ell}
\cdot
	e^{-\omega(1 - \ell/k)} n^{-1+\ell/k};
\]
algebraic manipulations using the fact that $\limsup d/k < 1$ and $2^d = n^{\oh1}$ then give
\[
	n \, \pr{ \VV \bcdot Z = 0, \: \mci = I }
=
	2^{-k} \oh1.
\label{eq-p5:typdist:l2mixing:proof:<:inf:small1}
\nt
\]
For $I \subseteq [k]$ with $d+2 \le \ell = \abs I \le (1 - \eta) k \eqqcolon L$, applying
(\ref{eq-p5:typdist:mi=I:inf}, \ref{eq-p5:typdist:gcd:pr-ex}, \ref{eq-p5:typdist:gcd:I>d:-})
gives
\[
	n \, \pr{ \VV \bcdot Z = 0, \: \mci = I }
\le
	15 \, n^{-1 + \ell/k}
\le
	2^{-k} n^{ -1 + L/k + \oh1 },
\label{eq-p5:typdist:l2mixing:proof:<:inf:small2}
\nt
\]
since $k \ll \log n$.
We now sum over the $I$ with $1 \le \abs I \le L = (1 - \eta)k$, using \cref{eq-p5:typdist:l2mixing:proof:<:inf:small1,eq-p5:typdist:l2mixing:proof:<:inf:small2}:
\[
	n \sumt{1 \le \abs I \le L}
	\pr{ \VV \bcdot Z = 0, \: \mci = I }
\le
	2^k \cdot 2^{-k} n^{-1 + L/k + \oh1}
\le
	n^{ -\eta + \oh1 }
=
	\oh1.
\label{eq-p5:typdist:l2mixing:proof:<:inf:sum-small}
\nt
\]
Finally we consider $I \subseteq [k]$ with $L \le \abs I \le k$.
As above,
by \cref{res-p5:typdist:sum>L},
we have
\begin{subequations}
	\label{eq-p5:typdist:l2mixing:proof:sum-large:q=inf}
\begin{empheq}[left = {%
	n \sumd{L \le \abs I \le k}
	\pr{ V \bcdot Z = 0, \: \mci = I }
\le
	\empheqlbrace}]%
{alignat=2}
	&1 + 3 \cdot 2^{d - L} = 1 + \oh1
	&&\quad \text{for undirected graphs},
\label{eq-p5:typdist:l2mixing:proof:sum-large:q=inf:-}
\\
	&1 + 5 \cdot (\tfrac32)^{2d - L} = 1 + \oh1
	&&\quad \text{for directed graphs}.
\label{eq-p5:typdist:l2mixing:proof:sum-large:q=inf:+}
\end{empheq}
\end{subequations}
The proof is completed by combining this with \cref{eq-p5:typdist:l2mixing:proof:<:inf:sum-small} and \cref{eq-p5:typdist:l2}.
\end{Proof}

\subsection{Adapting Proof to Directed Cayley Graphs}
\label{sec-p5:typdist:directed}

Where the random variable $A$ was uniform on a certain undirected lattice ball, it is now uniform on a directed ball (of a different radius).
Other than this, the only adaptation that needs be made is in determining the sizes of the discrete lattice balls: now instead of being a subset of $\mbz^k$, for some $k$, they are restricted to the first quadrant, ie to $\mbz_+^k$.
Assuming that their radius is large enough, this simply reduces their size by a factor (roughly) $2^k$.

Since all the sizes in question scale like $R^k$ when the ball-radius is $R$, when $k \ll \log n$ (and so $R \gg 1$), the desired radius for the directed ball is twice that of the undirected ball. When $k \asymp \log n$ (and we consider the $L_1$ ball), the directed ball has size $\binom{R+k}{k}$, so we are still interested in $R \asymp k \asymp \log n$, just the constant is different for directed compared with directed.

Finally, for directed graphs, we have a slightly weakened bound on the expected gcd, ie $\ex{\mfgcd^d \mid \mci = I}$; see \cref{res-p5:typdist:gcd}.
We addressed this in the proof of \cref{res-p5:typdist:l2-bound} at the time.

\section{Generating the Group with Uniform, Independent Generators}
\label{sec-p5:gen}

In this short section, we discuss (almost exclusively for Abelian groups) conditions on the number $k$ of generators chosen so that these elements \emph{actually generate} the group---that is, the graph $G_k$ is connected.
First we discuss the independent, uniform choice of generators, which has been our focus in the entire random Cayley graphs project.
Then we briefly discuss conditions used by other authors who have studied such Cayley graphs when $k \asymp 1$ (ie does not grow \asinf{\abs G}).

\subsection{Independent, Uniform Generators}
\label{sec-p5:gen:abelian}

Write $\varphi_k(G)$ for the probability that $k$ independent, uniform elements of $G$ generates the group.
To start, we state a worst-case result proved by \textcite[Lecture~1, Theorem~6]{P:gen-group-notes}:
if $\abs G \le 2^d$ then
\(
	\varphi_k(G) \ge \varphi_k(\mbz_2^d)
\)
for all $k$;
in words, $\mbz_2^d$ is the hardest group to generate.
As $\mbz_2^d$ forms a vector space,
one can calculate the probability $\varphi_k(\mbz_2^d)$ explicitly;
see \cref{eq-p5:gen:exact:full-group,eq-p5:gen:exact:p-group} below.

\medskip

We now move onto Abelian groups.
Recall that we write $d(G)$ for the minimal size of a generating set of $G$; abbreviate $d \cq d(G)$.
Clearly (by definition) we need $k \ge d$.
Draw $Z_1, Z_2, ... \sim \Unif(G)$.
\textcite{P:gen-abelian} shows that the expected number of independent, uniform generators required to generate the group is at most $d(G) + 3$.
As such, by Markov's inequality, $k - d(G) \gg 1$ is always sufficient to generate the group whp.
It is trivial that this is necessary for some groups:
	using the explicit expression for $\varphi_k(\mbz_2^d)$, it is easy to see that there exists a continuous, decreasing function $f : [0,\infty) \to (0,1)$ so that $\varphi_{d + C}(\mbz_2^d) \ge f(\ceil C)$ for all $C \in \mbn_0$.
Interestingly, for general Abelian groups, when $k - d(G) \gg 1$ is not required, $k = d(G) + 1$ is sufficient---there is no middle ground.

Recall that any group can be written as a direct product of $p$-groups; when the group is Abelian, this turns into a direct sum of Abelian $p$-groups.
Further, up to reordering, this decomposition is unique.
Use the following notation for such a decomposition:
\(
	G
=
	\oplus_{p \in \mcp} \: G_p;
\)
abbreviate $d_p \cq d(G_p)$.

\begin{lem}
\label{res-p5:gen:dichotomy}
	Use the notation established above.
	The following dichotomy holds.
	\begin{itemize}[itemsep = 0pt, topsep = \smallskipamount, label = \bcdot]
		\item 
		If there exists $p \in \mcp$ with $p \asymp 1 \asymp d - d_p$, then
		$k - d \gg 1$ is necessary for $\varphi_k(G) = 1 - \oh1$.
		
		\item 
		If there exists $\omega \gg 1$ with $\min_{p \in \mcp \cap [1,\omega]} (d - d_p) \gg 1$, then $\varphi_{d+1}(G) = 1 - \oh1$.
	\end{itemize}
\end{lem}

\begin{rmkt*}
	Our proof is strongly based on \textcite[Equations~2 and~4]{P:gen-abelian}; see \cref{eq-p5:gen:exact:full-group,eq-p5:gen:exact:p-group} below.
	The argument given in \cite{P:gen-abelian} is simple and well-known; \cite[Equation~2]{P:gen-abelian} is in essence \textcite[Lemma~1]{A:gen-group}.
	To be self-contained, we repeat it here.
	Our proof should not be thought of as novel, but rather combining already-established ideas to give a statement of the form we desire.
\end{rmkt*}

\begin{Proof}[Proof of \cref{res-p5:gen:dichotomy}]
Draw $Z_1, ..., Z_k \sim^\iid \Unif(G)$.
For each $i \in [k]$ and $p \in \mcp$, let $Z_i(p)$ be the projection of $Z_i$ to the $p$-group $G_p$.
Then $Z_i = (Z_i(p))_{p \in \mcp}$ where we identify $G$ with $\oplus_{p \in \mcp} \: G_p$.
Then $(Z_i(p) \mid i \in [k], p \in \mcp)$ are jointly independent.
For different $p \in \mcp$, the events
\[
	A_{p,k}
\cq
	\bra{ G_p \text{ is generated by } [Z_1(p), ..., Z_k(p)] }
\]
are independent.
Since
\[
	\bra{ G \text{ is generated by } [Z_1, ..., Z_k] }
=
	\cap_{p \in \mcp} \: A_{p,k},
\]
we deduce that
\[
	\varphi_k(G)
=
	\prodt{p \in \mcp}
	\varphi_k(G_p).
\label{eq-p5:gen:exact:full-group}
\nt
\]

We now find an expression for $\varphi_k(G_p)$ for a (prime) Abelian $p$-group.
We claim that
\[
	\varphi_k(G_p)
=
	\prodt[d_p]{\ell=1}
	\rbb{ 1 - p^{- (k-d_p) - \ell} }.
\label{eq-p5:gen:exact:p-group}
\nt
\]
To see this, note first that $G_p$ is isomorphic to the vector space of dimension $d_p$ over (the field) $\mbz_p$, under vector addition.
The process of generating this may be thought of as passing a series of tests:
	first we must choose a non-zero vector;
	next we must choose a vector not in the subspace generated by the previous one;
	next we must choose a vector not in the subspace generated by the previous two;
	etc.
If the dimension of the subspace already generated is $\ell$ and $\ell < d_p$, then the probability of choosing a vector not in this subspace is $1 - p^{\ell-d_p}$.
From this, one obtains \cref{eq-p5:gen:exact:p-group}.

\smallskip

We now consider the first claim of the lemma.
From \cref{eq-p5:gen:exact:full-group,eq-p5:gen:exact:p-group}, it is immediate that
\[
	\varphi_k(G)
\le
	\varphi_k(G_p)
\le
	1 - p^{-(k-d_p)-1}.
\]
Thus, using the hypotheses of this part, we have
\(
	1 - \varphi_k(G)
\ge
	p^{-(k-d_p)-1}
\asymp
	1,
\)
as required.

For the second part, we desire a lower bound on $\varphi_k(G_p)$.
From \cref{eq-p5:gen:exact:p-group}, we obtain
\[
	1 - \varphi_k(G_p)
=
	1
-	\prodt[d_p]{\ell=1}
	\rbb{ 1 - p^{- (k-d_p) - \ell} }
\le
	p^{-(k-d_p)-1}
	\sumt[d_p-1]{\ell=0}
	p^{-\ell}
\le
	2 p^{-(k-d_p)-1}.
\]
(In essence this, and the following, are union bounds.)
Plugging this into \cref{eq-p5:gen:exact:full-group},
we obtain
\[
	1 - \varphi_{d+1}(G)
=
	1 - \prodt{p \in \mcp} \varphi_{d+1}(G_p)
\le
	\sumt{p \in \mcp}
	\rbb{ 1 - \varphi_{d+1}(G_p) }
\le
	2 \sumt{p \in \mcp}
	p^{-(d-d_p)-2}.
\]
Separating this sum into $p \in \mcp$ with $p \le \omega$ and with $p \ge \omega$ and applying the hypotheses of this part, we see that it is $\oh1$.
That is,
\(
	\varphi_{d+1}(G)
=
	1 - \oh1,
\)
as required.
\end{Proof}

\subsection{Generating a Group with Order 1 Generators}
\label{sec-p5:gen:k=1}

As mentioned in
	\cite[\S\ref{sec-p1:intro:previous-work:typdist}]{HOt:rcg:matrix}
and
	\cite[\S\ref{sec-p3:intro:previous-work:typdist}]{HOt:rcg:abe:geom},
previous work on distance-metrics in random Abelian Cayley graphs had focussed on the case where $k$ is some fixed number, not diverging \asinf{\abs G}.
The previous section shows that in general the Abelian group will not be generated by this many independently and uniformly chosen generators.
Take, eg, $\mbz_{2n}$: to generate this group, at least one of the generators must be odd; however, the probability that none are odd is $2^{-k} \asymp 1$.

For each of the works which we referenced, we describe briefly the adaptations which they make to ensure that the group is generated (whp).

\begin{itemize}[label = \bcdot]
	\item 
	\textcite{AGg:diam-cayley-Zq} studied cyclic groups of prime order.
	In this case (due to the primality of the order), every element generates the group.
	Thus no conditions are required.
	
	\item 
	\textcite{MS:diameter-cayley-cycle} studied cycling groups of order $n$ (for a random $n$) without any primality assumption.
	The elements which they draw are chosen uniformly at random conditional on being jointly coprime with each other and with $n$.
	(Note that $d(\mbz_n) = 1$.)
	
	\item 
	\textcite{SZ:diam-cayley} studied Abelian groups of arbitrary (fixed, ie not diverging) rank without any primality assumption.
	The elements which they draw are chosen uniformly at random conditional on jointly generating the group.
	(Note that they may have $d(G) > 1$, but it must be bounded---$d(G) \le \rank(G)$, by definition.)
	
	\item 
	\textcite{EbP:cayley-diam-nil} determine, for fixed $k$, the diameter of various finite nilpotent groups of fixed step, including the upper triangular group which we study.
	They require that the rank be fixed and that $k$ be strictly larger than the rank.
	They sample uniformly conditional on generating the group.
	(It is standard that for a collection $Z$ to generate a nilpotent group $G$, it suffices for the projection of $Z$ to the Abelianisation $G^\ab$ to generate~$G^\ab$.)
\end{itemize}

\renewcommand{\bibfont}{\sffamily}
\renewcommand{\bibfont}{\sffamily\small}
\printbibliography[heading=bibintoc]

\end{document}

%% file: RCG_preamble_biblio.tex
\usepackage[doi=false,isbn=false,url=false,
	hyperref=auto,
	sorting=nyt,
	maxnames=10,
	maxcitenames=3,
	backend=biber,
	texencoding=auto,
	giveninits=true,
	block=space,
	style=numeric,
]{biblatex}

\DefineBibliographyStrings{english}{%
	andothers = {et al}
} 
 
\DeclareFieldFormat{eprint:arxiv}{%
	\href{https://arxiv.org/abs/#1}{\emph{arXiv:\allowbreak #1}}}

\DeclareFieldFormat{eprint:mrnumber}{%
	\href{http://www.ams.org/mathscinet-getitem?mr=MR#1}{MR\allowbreak #1}}

\DeclareFieldFormat{eprint:jstor}{%
	\href{http://www.jstor.org/stable/#1}{JSTOR\allowbreak #1}}

\DeclareFieldFormat{eprint:customeprint}{%
	\em Available at \upshape \ttfamily \href{https://#1}{#1}}

\DeclareFieldFormat{eprint:online}{%
	Available \href{https://#1}{online}}

\DeclareFieldFormat{eprint:inprep}{%
	\emph{In preparation}}

\DeclareFieldFormat{eprint:manual}{%
	\emph{#1}}

\DeclareFieldFormat{eprint:onarxiv}{%
	\emph{#1}}

\DeclareFieldFormat{eprint:toappear}{%
	\mbox{\emph{To appear}}}

\DeclareFieldFormat{eprint:accepted}{%
	\emph{Accepted at {#1}; to appear}}

\DeclareSourcemap{
	\maps[datatype=bibtex]{
		\map{
			\step[fieldsource=mrnumber,		fieldtarget=eprint, final]
			\step[fieldset=eprinttype,		fieldvalue=mrnumber]
		}
		\map{
			\step[fieldsource=arxiv,		fieldtarget=eprint, final]
			\step[fieldset=eprinttype,		fieldvalue=arxiv]
		}
		\map{
			\step[fieldsource=jstor,		fieldtarget=eprint, final]
			\step[fieldset=eprinttype,		fieldvalue=jstor]
		}
		\map{
			\step[fieldsource=customeprint,	fieldtarget=eprint, final]
			\step[fieldset=eprinttype,		fieldvalue=customeprint]
		}
		\map{
			\step[fieldsource=online,		fieldtarget=eprint, final]
			\step[fieldset=eprinttype,		fieldvalue=online]
		}
		\map{
			\step[fieldsource=inprep,		fieldtarget=eprint, final]
			\step[fieldset=eprinttype,		fieldvalue=inprep]
		}
		\map{
			\step[fieldsource=manual,		fieldtarget=eprint, final]
			\step[fieldset=eprinttype,		fieldvalue=manual]
		}
		\map{
			\step[fieldsource=onarxiv,		fieldtarget=eprint, final]
			\step[fieldset=eprinttype,		fieldvalue=onarxiv]
		}
		\map{
			\step[fieldsource=toappear,		fieldtarget=eprint, final]
			\step[fieldset=eprinttype,		fieldvalue=toappear]
		}
		\map{
			\step[fieldsource=accepted,		fieldtarget=eprint, final]
			\step[fieldset=eprinttype,		fieldvalue=accepted]
		}
	}
}

\DeclareFieldFormat{doi}{%
	\href{https://doi.org/#1}{DOI}%
}

\DeclareFieldFormat{url}{%
	\href{https://#1}{\texttt{#1}}%
}

\DeclareFieldFormat{comments}{%
	\emph{#1}%
} 
 
\DeclareFieldFormat
	[article,inbook,incollection,inproceedings,patent,thesis,unpublished]
	{title}{{#1\isdot}}

\DeclareFieldFormat
	[article,book,inbook,incollection,inproceedings,patent,thesis,unpublished, online]
	{date}{\mkbibparens{#1}}
\DeclareFieldFormat
	[article,book,inbook,incollection,inproceedings,patent,thesis,unpublished]
	{volume}{\mkbibbold{#1}}
\DeclareFieldFormat
	{pages}{\mkbibparens{#1}}

\DeclareFieldFormat{edition}{%
	\ifinteger{#1}
	{\mkbibordedition{#1}~\bibstring{edition}\addcomma}
	{#1~\bibstring{edition}\addcomma}}

\renewbibmacro*{volume+number+eid}{%
	\printfield{volume}%
\setunit*{\adddot}%
	\printfield{number}%
\setunit{\space}%
	\printfield{eid}%
}

\DeclareBibliographyDriver{article}{%
	\usebibmacro{bibindex}%
	\usebibmacro{begentry}%
	\usebibmacro{author}%
\setunit{\addspace}%
	\usebibmacro{date}%
\newunit\newblock
	\usebibmacro{title}%
\newunit\newblock
	\printfield{journaltitle}\addperiod%
\setunit*{\addspace}%
	\usebibmacro{volume+number+eid}%
\setunit*{\addspace}\newblock
	\printfield{pages}
\setunit*{\addspace}
	\printfield{comments}%
	\usebibmacro{eprint}%
\setunit*{\addspace}%
	\printfield{doi}%
	\usebibmacro{finentry}%
}

\DeclareBibliographyDriver{online}{%
	\usebibmacro{bibindex}%
	\usebibmacro{begentry}%
	\usebibmacro{author}%
\setunit{\addspace}%
	\usebibmacro{date}%
\newunit\newblock
	\usebibmacro{title}%
\newunit\newblock
	\usebibmacro{eprint}%
	\usebibmacro{finentry}%
}

\DeclareBibliographyDriver{book}{%
	\usebibmacro{bibindex}%
	\usebibmacro{begentry}%
	\usebibmacro{author}%
\setunit{\addspace}\newblock
	\usebibmacro{date}%
\newunit\newblock
	\usebibmacro{title}%
\setunit*{\addspace}\newblock
	\printfield{edition}%
\newunit\newblock
	\printlist{publisher}
\setunit*{\addspace}%
	\printfield{volume}%
\setunit*{\addspace}\newblock%
	\printfield{comments}%
	\usebibmacro{eprint}%
\setunit*{\addspace}
	\printfield{doi}
	\usebibmacro{finentry}%
}

\DeclareBibliographyDriver{thesis}{%
	\usebibmacro{bibindex}%
	\usebibmacro{begentry}%
	\usebibmacro{author}%
\setunit{\addspace}\newblock
	\usebibmacro{date}%
\newunit\newblock
	\usebibmacro{title}.%
\setunit*{\addspace}\newblock
	\space Thesis, \printlist{institution}
	\usebibmacro{eprint}%
	\printfield{comments}%
\setunit*{\addspace}
	\printfield{doi}%
	\usebibmacro{finentry}%
}

\DeclareBibliographyDriver{inproceedings}{%
	\usebibmacro{bibindex}%
	\usebibmacro{begentry}%
	\usebibmacro{author}%
\setunit{\addspace}%
	\usebibmacro{date}%
\newunit\newblock
	\usebibmacro{title}%
\newunit\newblock
	\printfield{booktitle}\addcomma%
\setunit*{\addspace}%
	\printfield{series}\addcomma%
\setunit*{\addspace}\newblock
	\usebibmacro{editor}\addcomma
\setunit*{\addspace}\newblock
	\printlist{publisher}
\setunit*{\addspace}%
	\printfield{volume}%
\setunit*{\addspace}%
	\printfield{pages}
\setunit*{\addspace}
	\usebibmacro{eprint}%
	\printfield{comments}%
\setunit*{\addnnspace}
	\printfield{doi}
	\usebibmacro{finentry}%
}

\DeclareBibliographyDriver{inbook}{%
	\usebibmacro{bibindex}%
	\usebibmacro{begentry}%
	\usebibmacro{author}%
\setunit{\addspace}%
	\usebibmacro{date}%
\newunit\newblock
	\usebibmacro{title}%
\newunit\newblock
	\printfield{booktitle}\addcomma%
\setunit*{\addspace}%
	\printfield{series}\addcomma%
\setunit*{\addspace}\newblock
	\usebibmacro{editor}\addcomma
\setunit*{\addspace}\newblock
	\printlist{publisher}
\setunit*{\addspace}%
	\printfield{volume}%
\setunit*{\addspace}%
	\printfield{pages}
\setunit*{\addspace}
	\usebibmacro{eprint}%
	\printfield{comments}%
\setunit*{\addspace}
	\printfield{doi}
	\usebibmacro{finentry}%
}

\DeclareBibliographyDriver{incollection}{%
	\usebibmacro{bibindex}%
	\usebibmacro{begentry}%
	\usebibmacro{author}%
\setunit{\addspace}%
	\usebibmacro{date}%
\newunit\newblock
	\usebibmacro{title}%
\newunit\newblock
	\printfield{booktitle}\addcomma%
\setunit*{\addspace}%
	\printfield{series}\addcomma%
\setunit*{\addspace}\newblock
	\printlist{publisher}
\setunit*{\addspace}%
	\printfield{volume}%
\setunit*{\addspace}%
	\printfield{pages}
\setunit*{\addspace}
	\usebibmacro{eprint}%
	\printfield{comments}%
\setunit*{\addspace}
	\printfield{doi}
	\usebibmacro{finentry}%
}

\AtEveryBibitem{%
	\clearfield{day}%
	\clearfield{month}%
} 
 
\DeclareNameAlias{sortname}{given-family}
\DeclareNameAlias{default}{given-family}

%% file: RCG_preamble_main.tex
\newcommand{\acknofootnote}{%
	\blfootnote{%
		\hspace*{1em}%
		Jonathan Hermon%
	\hfill%
		Sam Olesker-Taylor%
		\hspace*{1em}%
	\\%
		\href{mailto:jhermon@math.ubc.ca}{jhermon@math.ubc.ca},
		\href{http://www.math.ubc.ca/~jhermon/}{math.ubc.ca/$\sim$jhermon/}%
	\hfill%
		\href{mailto:sam.ot@posteo.co.uk}{sam.ot@posteo.co.uk},
		\href{https://mathematicalsam.wordpress.com}{mathematicalsam.wordpress.com}%
	\\%
		University of British Columbia, Vancouver, Canada%
	\hfill%
		Department of Mathematical Sciences, University of Bath, UK%
	\\%
		Supported by EPSRC EP/L018896/1 and an NSERC Grant%
	\hfill%
		Supported by EPSRC Grants 1885554 and EP/N004566/1%
	\par\smallskip\par
	\centering%
		The vast majority of this work was undertaken while both authors were at the University of Cambridge%
	}
}

\newcommand*{\mm}{\ensuremath{L}}
\newcommand*{\pp}{\ensuremath{p}}
\newcommand*{\qq}{\ensuremath{q}}

\newcommand*{\WW}{\ensuremath{W}}

\newcommand*{\VV}{\ensuremath{V}}

\newcommand{\printtoc}[1]{%
	\ifthenelse%
		{\equal{#1}{1}}%
		{\sffamily\boldmath\tableofcontents\unboldmath\normalfont}%
		{\newpage\small\sffamily\boldmath\tableofcontents\unboldmath\normalfont\normalsize}%
	}

\newcommand{\nextresult}{%
	\setcounter{introthm}{\value{introthm}}
	\setcounter{introcor}{\value{introthm}}
	\setcounter{introconj}{\value{introthm}}
	\setcounter{introdefn}{\value{introthm}}
	\setcounter{intrormkT}{\value{introthm}}
}

\usepackage{indentfirst, environ}
\usepackage[utf8]{inputenc}

\usepackage{empheq}

\usepackage{manyfoot}
\SetFootnoteHook{\hspace*{-1.8em}}
\DeclareNewFootnote{bl}[gobble]
\newcommand{\blfootnote}[1]{\footnotebl{\sffamily#1}}

\usepackage{graphicx}
\graphicspath{{../Figures/}}
\usepackage[labelsep=period]{caption}

\usepackage{chngcntr}
\counterwithin{figure}{section}

\usepackage[UKenglish]{babel}
\usepackage{amsmath, amsthm, amssymb, mathrsfs, bm, mathtools}
\usepackage{wasysym}
\usepackage{mathtools}

\allowdisplaybreaks

\usepackage{cases}

\usepackage{titlesec}
\usepackage{titletoc}

\titleformat{\title}
{\sffamily \Huge \bfseries \scshape \boldmath}{}{}{}
\titleformat{\section}
{\sffamily \Large \bfseries \scshape \boldmath}{\thesection}{1em}{}
\titleformat{\subsection}
{\sffamily \large \bfseries \scshape \boldmath}{\thesubsection}{1em}{}
\titleformat{\subsubsection}
{\sffamily \normalsize \bfseries \scshape \boldmath}{\thesubsubsection}{1em}{}
\titleformat{\paragraph}
{\sffamily \normalsize \bfseries \scshape \boldmath}{\theparagraph}{1em}{}
\titleformat{\subparagraph}[runin]
{\sffamily \normalsize \bfseries \scshape \boldmath}{\thesubparagraph}{1em}{}

\def\IfAmpersandUseAlign#1#2&#3\EndIfAmpersandUseAlign
{%
	\if\relax\detokenize{#3}\relax
	\begin{equation*}%
	#1%
	\end{equation*}%
	\else
	\begin{align*}%
	#1%
	\end{align*}%
	\fi
}
\def\[#1\]%
{%
	\IfAmpersandUseAlign{#1}#1&\EndIfAmpersandUseAlign
}

\usepackage{eqparbox}
\newcommand{\eqmathsbox}[3][\mathrel]{%
	#1{\eqmakebox[#2]{$\displaystyle#3$}}%
}

\usepackage{scrextend}

\usepackage{enumitem}

\newcommand{\numberingroman}{%
	\renewcommand{\labelenumi}{(\roman{enumi})}%
	\renewcommand{\theenumi}{(\roman{enumi})}%
}

\setlist[description]{%
	topsep		= 0pt,		
	noitemsep,				
	font		= {\mdseries\itshape},	
}

\newcommand{\gab}{G^\ab}

\newcommand{\UU}{U}
\newcommand{\HH}{H}

\newcommand{\ugr}{\UU_{p,d}}
\newcommand{\uab}{\UU_{p,d}^\ab}

\newcommand{\uk}{(\UU_{p,d})_k}

\newcommand{\hgr}{\HH_{p,d}}
\newcommand{\hab}{\HH_{p,d}^\ab}

\newcommand{\hk}{(\HH_{p,d})_k}

\usepackage[dvipsnames,hyperref]{xcolor}

\newcommand{\nt}{\addtocounter{equation}{1}\tag{\theequation}}

\newcommand{\bcdot}{\ensuremath{\bm{\cdot}}}
\let\mod\relax
\DeclareMathOperator{\mod}{\, mod}
\DeclareMathOperator{\MOD}{mod}

\DeclareMathOperator{\ord}{ord}

\newcommand{\Quad}[1]{
	\mathchoice
	{\quad\text{#1}\quad}
	{\text{ #1 }}
	{\text{ #1 }}
	{\text{ #1 }}
}

\newcommand{\whp}{\text{whp}\xspace}

\newcommand{\Qforall}{\Quad{for all}}
\newcommand{\Qfor}{\Quad{for}}
\newcommand{\Qand}{\Quad{and}}
\newcommand{\Qwhere}{\Quad{where}}

\newcommand{\id}{\mathsf{id}}

\newcommand{\typ}{\mathsf{typ}}
\newcommand{\typloc}{\mathsf{typ}_\ell}
\newcommand{\typglo}{\mathsf{typ}_g}

\newcommand{\cq}{\coloneqq}
\renewcommand{\epsilon}{\varepsilon}

\newcommand{\eps}{\epsilon}

\usepackage{xifthen}

\newcommand{\binomt}[2]{ \textstyle \binom{#1}{#2} \displaystyle }

\newcommand{\maxt}[1]{ \textstyle \max_{#1} \displaystyle }

\newcommand{\mint}[1]{ \textstyle \min_{#1} \displaystyle }

\usepackage{calc}
\newlength{\halfplusheight}
\newcommand{\MAX}[1]{\mathop{\raisebox{\halfplusheight}{\(\displaystyle\max_{#1}\)}}}

\newcommand{\LIM}[1]{\mathop{\raisebox{\halfplusheight}{\(\displaystyle\lim_{#1}\)}}}
\newcommand{\LIMSUP}[1]{\mathop{\raisebox{\halfplusheight}{\(\displaystyle\limsup_{#1}\)}}}
\newcommand{\LIMINF}[1]{\mathop{\raisebox{\halfplusheight}{\(\displaystyle\liminf_{#1}\)}}}

\DeclareMathOperator*{\sumTT}{\textstyle\sum}
\newcommand{\sumT}[2][]{
	\ifthenelse{\isempty{#1}}
	{\sumTT_{#2}}
	{\sumTT_{#2}^{#1}}
}
\newcommand{\sumt}[2][]{
	\ifthenelse{\isempty{#1}}
	{\textstyle \sum_{#2}      \displaystyle}
	{\textstyle \sum_{#2}^{#1} \displaystyle}
}
\newcommand{\sumd}[2][]{
	\ifthenelse{\isempty{#1}}
	{\displaystyle \sum_{#2}}
	{\displaystyle \sum_{#2}^{#1}}
}

\newcommand{\intt}[2][]{
	\ifthenelse{\isempty{#1}}
	{\textstyle \int_{#2}      \displaystyle}
	{\textstyle \int_{#2}^{#1} \displaystyle}
}

\newcommand{\prodt}[2][]{
	\ifthenelse{\isempty{#1}}
	{\textstyle \prod_{#2}      \displaystyle}
	{\textstyle \prod_{#2}^{#1} \displaystyle}
}
\newcommand{\prodd}[2][]{
	\ifthenelse{\isempty{#1}}
	{\prod_{#2}}
	{\prod_{#2}^{#1}}
}

\usepackage{xspace}

\let\originalexp\exp
\let\exp\relax

\DeclareRobustCommand{\exp} [1]{\originalexp(#1)}

\newcommand{\abs}  [1]{| #1 |}
\newcommand{\absb} [1]{\big| #1 \bigr|}

\newcommand{\norm}  [1]{\lVert #1 \rVert}
\newcommand{\normb} [1]{\big\lVert #1 \bigr\rVert}

\newcommand{\rbr} [1]{ ( #1 ) }
\newcommand{\rbb} [1]{\bigl( #1 \bigr)}

\newcommand{\sbb} [1]{\bigl[ #1 \bigr]}

\newcommand{\bra} [1]{ \{ #1 \} }
\newcommand{\brb} [1]{\bigl\{ #1 \bigr\}}

\newcommand{\dist}{\mathrm{dist}}
\newcommand{\TV}{\mathrm{TV}}

\DeclareMathOperator{\diam}{diam}

\newcommand{\mix}{\mathrm{mix}}

\newcommand{\rel}{\mathrm{rel}}

\newcommand{\st}{{ \ \mathrm{st} \ }}

\DeclareMathOperator{\vol}{\mathrm{vol}}

\newcommand{\ab} {\mathrm{ab}}
\newcommand{\com}{\mathrm{com}}
\DeclareMathOperator{\rank}{rank}
\DeclareMathOperator{\step}{step}

\newcommand{\Unif}{\mathrm{Unif}}
\newcommand{\iid}{\mathrm{iid}}

\newcommand{\tmix}{t_\mix}

\newcommand{\trel}{t_\rel}

\newcommand{\Ninn}{{N\in\mathbb{N}}}

\newcommand{\floor}[1]{\lfloor #1 \rfloor}
\newcommand{\floorb}[1]{\bigl \lfloor #1 \bigr \rfloor}

\newcommand{\ceil}[1]{\lceil #1 \rceil}
\newcommand{\ceilb}[1]{\bigl\lceil #1 \bigr\rceil}

\newcommand{\midb}{\bigm\vert}

\newcommand{\one}  [1]{\bm1( #1 )}
\newcommand{\oneb} [1]{\bm1\bigl( #1 \bigr)}

\newcommand{\tvb}[1]{\bigl\lVert #1 \bigr\rVert_{\mathrm{TV}}}

\newcommand{\logk}[1][]{
	\ifthenelse{\equal{}{#1}}
	{\log k}
	{(\log k)^{#1}}
}

\newcommand{\logn}[1][]{
	\ifthenelse{\equal{}{#1}}
	{\log n}
	{(\log n)^{#1}}
}

\newcommand{\logm}[1][]{
	\ifthenelse{\equal{}{#1}}
	{\log m}
	{(\log m)^{#1}}
}

\newcommand{\loglogn}[1][]{
	\ifthenelse{\equal{}{#1}}
	{\log\log n}
	{(\log\log n)^{#1}}
}

\newcommand{\prt}[2][]{
	\ifthenelse{\equal{}{#1}}
	{\mathbb{P}(#2)}
	{\mathbb{P}_{#1}(#2)}
}
\newcommand{\pr}[2][]{
	\mathchoice
	{\ifthenelse{\isempty{#1}}
		{\mathbb{P}\bigl(#2\bigr)}
		{\mathbb{P}_{#1}\bigl(#2\bigr)}}
	{\ifthenelse{\isempty{#1}}
		{\mathbb{P}(#2)}
		{\mathbb{P}_{#1}(#2)}}
	{\ifthenelse{\isempty{#1}}
		{\mathbb{P}(#2)}
		{\mathbb{P}_{#1}(#2)}}
	{\ifthenelse{\isempty{#1}}
		{\mathbb{P}(#2)}
		{\mathbb{P}_{#1}(#2)}}
}
\newcommand{\prb}[2][]{
	\ifthenelse{\equal{}{#1}}
	{\mathbb{P}\bigl( #2 \bigr)}
	{\mathbb{P}_{#1}\bigl( #2 \bigr)}
}
\newcommand{\prB}[2][]{
	\ifthenelse{\equal{}{#1}}
	{\mathbb{P}\Bigl( #2 \Bigr)}
	{\mathbb{P}_{#1}\Bigl( #2 \Bigr)}
}
\newcommand{\prbb}[2][]{
	\ifthenelse{\equal{}{#1}}
	{\mathbb{P}\biggl( #2 \biggr)}
	{\mathbb{P}_{#1}\biggl( #2 \biggr)}
}
\newcommand{\prBB}[2][]{
	\ifthenelse{\equal{}{#1}}
	{\mathbb{P}\Biggl( #2 \Biggr)}
	{\mathbb{P}_{#1}\Biggl( #2 \Biggr)}
}
\newcommand{\prs}[2][]{
	\ifthenelse{\equal{}{#1}}
	{\mathbb{P}\left( #2 \right)}
	{\mathbb{P}_{#1}\left( #2 \right)}
}

\newcommand{\qr}[2][]{
	\mathchoice
	{\ifthenelse{\isempty{#1}}
		{\mathbb{Q}\bigl(#2\bigr)}
		{\mathbb{Q}_{#1}\bigl(#2\bigr)}}
	{\ifthenelse{\isempty{#1}}
		{\mathbb{Q}(#2)}
		{\mathbb{Q}_{#1}(#2)}}
	{\ifthenelse{\isempty{#1}}
		{\mathbb{Q}(#2)}
		{\mathbb{Q}_{#1}(#2)}}
	{\ifthenelse{\isempty{#1}}
		{\mathbb{Q}(#2)}
		{\mathbb{Q}_{#1}(#2)}}
}
\newcommand{\qrb}[2][]{
	\ifthenelse{\equal{}{#1}}
	{\mathbb{Q}\bigl( #2 \bigr)}
	{\mathbb{Q}_{#1}\bigl( #2 \bigr)}
}
\newcommand{\qrB}[2][]{
	\ifthenelse{\equal{}{#1}}
	{\mathbb{Q}\Bigl( #2 \Bigr)}
	{\mathbb{Q}_{#1}\Bigl( #2 \Bigr)}
}
\newcommand{\qrbb}[2][]{
	\ifthenelse{\equal{}{#1}}
	{\mathbb{Q}\biggl( #2 \biggr)}
	{\mathbb{Q}_{#1}\biggl( #2 \biggr)}
}
\newcommand{\qrBB}[2][]{
	\ifthenelse{\equal{}{#1}}
	{\mathbb{Q}\Biggl( #2 \Biggr)}
	{\mathbb{Q}_{#1}\Biggl( #2 \Biggr)}
}
\newcommand{\qrs}[2][]{
	\ifthenelse{\equal{}{#1}}
	{\mathbb{Q}\left( #2 \right)}
	{\mathbb{Q}_{#1}\left( #2 \right)}
}

\newcommand{\ext}[2][]{
\ifthenelse{\equal{}{#1}}
{\mathbb{E}(#2)}
{\mathbb{E}_{#1}(#2)}
}
\newcommand{\ex}[2][]{
	\mathchoice
	{\ifthenelse{\isempty{#1}}
		{\mathbb{E}\bigl(#2\bigr)}
		{\mathbb{E}_{#1}\bigl(#2\bigr)}}
	{\ifthenelse{\isempty{#1}}
		{\mathbb{E}(#2)}
		{\mathbb{E}_{#1}(#2)}}
	{\ifthenelse{\isempty{#1}}
		{\mathbb{E}(#2)}
		{\mathbb{E}_{#1}(#2)}}
	{\ifthenelse{\isempty{#1}}
		{\mathbb{E}(#2)}
		{\mathbb{E}_{#1}(#2)}}
}
\newcommand{\exb}[2][]{
	\ifthenelse{\equal{}{#1}}
	{\mathbb{E}\bigl( #2 \bigr)}
	{\mathbb{E}_{#1}\bigr( #2 \bigr)}
}
\newcommand{\exB}[2][]{
	\ifthenelse{\equal{}{#1}}
	{\mathbb{E}\Bigl( #2 \Bigr)}
	{\mathbb{E}_{#1}\Bigl( #2 \Bigr)}
}
\newcommand{\exbb}[2][]{
	\ifthenelse{\equal{}{#1}}
	{\mathbb{E}\biggl( #2 \biggr)}
	{\mathbb{E}_{#1}\biggl( #2 \biggr)}
}
\newcommand{\exBB}[2][]{
	\ifthenelse{\equal{}{#1}}
	{\mathbb{E}\Biggl( #2 \Biggr)}
	{\mathbb{E}_{#1}\Biggl( #2 \Biggr)}
}

\newcommand{\fx}[2][]{
	\ifthenelse{\equal{}{#1}}
	{\mathbb{F}(#2)}
	{\mathbb{F}_{#1}(#2)}
}
\newcommand{\fxb}[2][]{
	\ifthenelse{\equal{}{#1}}
	{\mathbb{F}\bigl( #2 \bigr)}
	{\mathbb{F}_{#1}\bigr( #2 \bigr)}
}
\newcommand{\fxB}[2][]{
	\ifthenelse{\equal{}{#1}}
	{\mathbb{F}\Bigl( #2 \Bigr)}
	{\mathbb{F}_{#1}\Bigl( #2 \Bigr)}
}
\newcommand{\fxbb}[2][]{
	\ifthenelse{\equal{}{#1}}
	{\mathbb{F}\biggl( #2 \biggr)}
	{\mathbb{F}_{#1}\biggl( #2 \biggr)}
}
\newcommand{\fxBB}[2][]{
	\ifthenelse{\equal{}{#1}}
	{\mathbb{F}\Biggl( #2 \Biggr)}
	{\mathbb{F}_{#1}\Biggl( #2 \Biggr)}
}

\newcommand{\Var}[1]{\mathbb{V}\mathrm{ar}(#1)}
\newcommand{\Varb}[2][]{
	\ifthenelse{\equal{}{#1}}
	{\mathbb{V}\mathrm{ar} \bigl(#2\bigr)}
	{\mathbb{V}\mathrm{ar}_{#1} \bigl(#2\bigr)}
}
\newcommand{\VAR}[2][]{
	\ifthenelse{\equal{}{#1}}
	{\mathrm{Var}(#2)}
	{\mathrm{Var}_{#1}(#2)}
}

\newcommand{\Oh}  [1]{\mathcal{O}( #1 )}
\newcommand{\Ohb} [1]{\mathcal{O}\bigl( #1 \bigr)}

\newcommand{\oh}  [1]{o( #1 )}
\newcommand{\ohb} [1]{o\bigl( #1 \bigr)}

\newcommand{\mbf}{\mathbb{F}}

\newcommand{\mbn}{\mathbb{N}}

\newcommand{\mbp}{\mathbb{P}}

\newcommand{\mbr}{\mathbb{R}}

\newcommand{\mbz}{\mathbb{Z}}

\newcommand{\mcb}{\mathcal{B}}
\newcommand{\mcc}{\mathcal{C}}
\newcommand{\mcd}{\mathcal{D}}
\newcommand{\mce}{\mathcal{E}}

\newcommand{\mcg}{\mathcal{G}}

\newcommand{\mci}{\mathcal{I}}

\newcommand{\mcp}{\mathcal{P}}

\newcommand{\mcs}{\mathcal{S}}

\newcommand{\mcw}{\mathcal{W}}

\newcommand{\mfd}{\mathfrak{D}}

\newcommand{\mfgcd}{\mathfrak{g}}

\newcommand{\mfm}{\mathfrak{M}}

\newcommand{\toinf}[1]{\ensuremath{#1\to\infty}}
\newcommand{\asinf}[1]{\text{as \ensuremath{#1\to\infty}}}

\newcommand{\kinf}{{k\to\infty}}

\newcommand{\ninf}{{n\to\infty}}

\newcommand{\Ninf}{{N\to\infty}}

\makeatletter
\newenvironment{subtheorem}[1]{%
	\def\subtheoremcounter{#1}%
	\refstepcounter{#1}%
	\protected@edef\theparentnumber{\csname the#1\endcsname}%
	\setcounter{parentnumber}{\value{#1}}%
	\setcounter{#1}{0}%
	\expandafter\def\csname the#1\endcsname{\theparentnumber\alph{#1}}%
	\expandafter\def\csname theH#1\endcsname{thm.\theparentnumber\alph{#1}}%
	\unskip\ignorespaces
}{%
	\setcounter{\subtheoremcounter}{\value{parentnumber}}%
	\ignorespacesafterend
}
\makeatother
\newcounter{parentnumber}

\makeatletter
\newenvironment{subtheorem-num}[1]{%
	\def\subtheoremcounter{#1}%
	\refstepcounter{#1}%
	\protected@edef\theparentnumber{\csname the#1\endcsname}%
	\setcounter{parentnumber}{\value{#1}}%
	\setcounter{#1}{0}%
	\expandafter\def\csname the#1\endcsname{\theparentnumber.\arabic{#1}}%
	\expandafter\def\csname theH#1\endcsname{thm.\theparentnumber.\arabic{#1}}%
	\unskip\ignorespaces
}{%
	\setcounter{\subtheoremcounter}{\value{parentnumber}}%
	\ignorespacesafterend
}
\makeatother

\usepackage[colorlinks, pdfauthor = {Jonathan\ Hermon\ and\ Sam\ Olesker-Taylor}, pdfusetitle]{hyperref}
\hypersetup{
	colorlinks,
	linkcolor	= {blue},
	filecolor	= {red},
	urlcolor	= {blue},
	citecolor	= {ForestGreen}
}

\newcommand{\qedtriangle}{\renewcommand{\qedsymbol}{\ensuremath{\triangle}}}

\usepackage[capitalise]{cleveref}

\crefformat{equation}{\textup{(#2#1#3)}}
\crefmultiformat{equation}{\textup{(#2#1#3}}{\textup{, #2#1#3)}}{\textup{, #2#1#3}}{\textup{, #2#1#3)}}
\crefrangeformat{equation}{\textup{(#3#1#4--#5#2#6)}}

\newenvironment{Proof}[1][\proofname]{%
	\proof[\upshape\bfseries\sffamily\boldmath#1]
}{\endproof}

\usepackage{etoolbox}
\usepackage{needspace}

\newtheoremstyle{sfsl}
{1\baselineskip}		
{1\baselineskip}		
{\slshape}				
{}						
{\bfseries\sffamily}	
{.}						
{0.5em}					
{\thmname{#1}\thmnumber{ #2}\thmnote{ {\mdseries(#3)}}}

\newtheoremstyle{sfup}
{1\baselineskip}		
{1\baselineskip}		
{\upshape}				
{}						
{\bfseries\sffamily}	
{.}						
{0.5em}					
{\thmname{#1}\thmnumber{ #2}\thmnote{ {\mdseries(#3)}}}

\theoremstyle{sfsl}

\newtheorem*{thm*}{Theorem}
\newtheorem{thm} {Theorem}[section]
\crefname{thm}{Theorem}{Theorems}

\newtheorem*{introthm*}{Theorem}
\newtheorem{introthm}{Theorem}

\crefname{introthm}{Theorem}{Theorems}

\newtheorem*{cor*}{Corollary}
\newtheorem{cor} [thm]{Corollary}
\crefname{cor}{Corollary}{Corollaries}

\newtheorem*{introcor*}{Corollary}

\crefname{introcor}{Corollary}{Corollaries}

\newtheorem*{introconj*}{Conjecture}

\crefname{introconj}{Conjecture}{Conjectures}

\newtheorem*{introques*}{Question}

\crefname{introques}{Question}{Questions}

\newtheorem*{lem*}    {Lemma}
\newtheorem{lem} [thm]{Lemma}
\crefname{lem}{Lemma}{Lemmas}

\newtheorem*{introlem*}{Lemma}

\crefname{introlem}{Lemma}{Lemmas}

\newtheorem*{prop*}    {Proposition}
\newtheorem{prop} [thm]{Proposition}
\crefname{prop}{Proposition}{Propositions}

\newtheorem*{clm*}    {Claim}

\crefname{clm}{Claim}{Claims}

\newtheorem*{defn*}    {Definition}
\newtheorem{defn} [thm]{Definition}
\crefname{defn}{Definition}{Definitions}

\newtheorem*{introdefn*}{Definition}
\newtheorem{introdefn}{Definition}

\crefname{introdefn}{Definition}{Definitions}

\newtheorem{innercustomgenericsl}{\customgenericnamesl}
\providecommand{\customgenericnamesl}{}
\newcommand{\newcustomtheoremsl}[2]{%
	\newenvironment{#1}[1]
	{%
		\renewcommand\customgenericnamesl{#2}%
		\renewcommand\theinnercustomgenericsl{##1}%
		\innercustomgenericsl
	}
	{\endinnercustomgenericsl}
}

\newcustomtheoremsl{customthm}{Theorem}
\newcustomtheoremsl{customprop}{Proposition}
\newcustomtheoremsl{customlemma}{Lemma}
\newcustomtheoremsl{customrmk}{Remark}
\newcustomtheoremsl{customconj}{Conjecture}
\newcustomtheoremsl{customdefn}{Definition}
\newcustomtheoremsl{customquestion}{Question}
\newcustomtheoremsl{customopenproblem}{Open Problem}
\newcustomtheoremsl{customhyp}{Hypothesis}

\newtheorem*{conj*}   {Conjecture}

\crefname{conj}{Conjecture}{Conjectures}

\newenvironment{conj-ind*}
	{\begin{quote}\textsf{\textbf{Conjecture.}}\slshape}
	{\end{quote}}
\newenvironment{conj-ind}
	{\begin{quote}\vspace{-\glueexpr\baselineskip+\topsep}\begin{customconj}}
	{\end{customconj}\end{quote}}

\newenvironment{question-ind*}
	{\begin{quote}\textsf{\textbf{Question.}}\slshape}
	{\end{quote}}
\newenvironment{question-ind}
	{\begin{quote}\vspace{-\glueexpr\baselineskip+\topsep}\begin{customquestion}}
	{\end{customquestion}\end{quote}}

\newenvironment{openproblem-ind*}
	{\begin{quote}\textsf{\textbf{Open Problem.}}\slshape}
	{\end{quote}}
\newenvironment{openproblem-ind}
	{\begin{quote}\vspace{-\glueexpr\baselineskip+\topsep}\begin{customopenproblem}}
	{\end{customopenproblem}\end{quote}}

\newtheorem*{hypothesis*}{Hypothesis}

\newtheorem*{hyp*}{Hypothesis}
\newtheorem{hyp}{Hypothesis}

\crefname{hyp}{Hypothesis}{Hypotheses}

\newtheorem*{rmk*}{Remark}

\theoremstyle{sfup}

\newtheorem{innercustomgenericup}{\customgenericnameup}
\providecommand{\customgenericnameup}{}
\newcommand{\newcustomtheoremup}[2]{%
	\newenvironment{#1}[1]
	{%
		\renewcommand\customgenericnameup{#2}%
		\renewcommand\theinnercustomgenericup{##1}%
		\innercustomgenericup
	}
	{\endinnercustomgenericup}
}

\crefname{exm} {Example}{Examples}
\crefname{exmT}{Example}{Examples}

\newtheorem{rmkT}[thm]{Remark}
	\newenvironment{rmkt}
	{\pushQED{\qed}\renewcommand{\qedsymbol}{\ensuremath{\triangle}}\rmkT}
	{\popQED\endrmkT}
\crefname{rmk} {Remark}{Remarks}
\crefname{rmkT}{Remark}{Remarks}

\newtheorem*{rmkTT}{Remark}
\newenvironment{rmkt*}
	{\pushQED{\qed}\renewcommand{\qedsymbol}{\ensuremath{\triangle}}\rmkTT}
	{\popQED\endrmkTT}

\newcustomtheoremup{customrmkT}{Remark}

\crefname{rmks} {Remarks}{Remarks}
\crefname{rmksT}{Remarks}{Remarks}

\newtheorem*{rmks*} {Remarks}
\newtheorem*{rmksTT}{Remarks}
\newenvironment{rmkst*}
	{\pushQED{\qed}\renewcommand{\qedsymbol}{\ensuremath{\triangle}}\rmksTT}
	{\popQED\endrmksTT}

\newtheorem{intrormkT}{Remark}
	\newenvironment{intrormkt}
	{\pushQED{\qed}\renewcommand{\qedsymbol}{\ensuremath{\triangle}}\intrormkT}
	{\popQED\endintrormkT}

\crefname{intrormk} {Remark}{Remarks}
\crefname{intrormkT}{Remark}{Remarks}

\newtheorem*{intrormk*} {Remark}
\newtheorem*{intrormkTT}{Remark}
\newenvironment{intrormkt*}
	{\pushQED{\qed}\renewcommand{\qedsymbol}{\ensuremath{\triangle}}\intrormkTT}
	{\popQED\endintrormkTT}

\newtheorem*{exm*} {Example}
\newtheorem*{exmTT}{Lemma}
	\newenvironment{exmt*}
	{\pushQED{\qed}\renewcommand{\qedsymbol}{\ensuremath{\triangle}}\exmTT}
	{\popQED\endexmTT}

\newtheorem*{note*} {Note}
\newtheorem*{noteTT}{Note}
	\newenvironment{notet*}
	{\pushQED{\qed}\renewcommand{\qedsymbol}{\ensuremath{\triangle}}\noteTT}
	{\popQED\endnoteTT}

\catcode`\^^A=14

\usepackage{xparse}

\newcounter{mixedsubequations}

\ExplSyntaxOn

\int_new:N \g_mixedsubeq_int

\NewDocumentEnvironment{mixedsubequations}{o}
{
	\IfNoValueTF { #1 }
	{
		\addtocounter{equation}{-\g_mixedsubeq_int}
		\stepcounter{mixedsubequations}
	}
	{
		\int_gset:Nn \g_mixedsubeq_int { \clist_count:n { #1 } }
		\clist_map_inline:nn { #1 }
		{
			\refstepcounter{equation}\label{##1}
		}
		\addtocounter{equation}{-\g_mixedsubeq_int}
		\setcounter{mixedsubequations}{1}
	}
	\domixedsubequations
}
{\ignorespacesafterend}

\NewDocumentCommand{\domixedsubequations}{}
{
	\cs_set:Npx \theequation
	{
		\exp_not:o { \theequation }
		\exp_not:n { \alph{mixedsubequations} }
	}
	\ignorespaces
}
\ExplSyntaxOff

\makeatletter
\let\save@mathaccent\mathaccent
\newcommand*\if@single[3]{%
  \setbox0\hbox{${\mathaccent"0362{#1}}^H$}%
  \setbox2\hbox{${\mathaccent"0362{\kern0pt#1}}^H$}%
  \ifdim\ht0=\ht2 #3\else #2\fi
  }
\newcommand*\rel@kern[1]{\kern#1\dimexpr\macc@kerna}
\newcommand*\widebar[1]{\@ifnextchar^{{\wide@bar{#1}{0}}}{\wide@bar{#1}{1}}}
\newcommand*\wide@bar[2]{\if@single{#1}{\wide@bar@{#1}{#2}{1}}{\wide@bar@{#1}{#2}{2}}}
\newcommand*\wide@bar@[3]{%
  \begingroup
  \def\mathaccent##1##2{%
    \let\mathaccent\save@mathaccent
    \if#32 \let\macc@nucleus\first@char \fi
    \setbox\z@\hbox{$\macc@style{\macc@nucleus}_{}$}%
    \setbox\tw@\hbox{$\macc@style{\macc@nucleus}{}_{}$}%
    \dimen@\wd\tw@
    \advance\dimen@-\wd\z@
    \divide\dimen@ 3
    \@tempdima\wd\tw@
    \advance\@tempdima-\scriptspace
    \divide\@tempdima 10
    \advance\dimen@-\@tempdima
    \ifdim\dimen@>\z@ \dimen@0pt\fi
    \rel@kern{0.6}\kern-\dimen@
    \if#31
      \overline{\rel@kern{-0.6}\kern\dimen@\macc@nucleus\rel@kern{0.4}\kern\dimen@}%
      \advance\dimen@0.4\dimexpr\macc@kerna
      \let\final@kern#2%
      \ifdim\dimen@<\z@ \let\final@kern1\fi
      \if\final@kern1 \kern-\dimen@\fi
    \else
      \overline{\rel@kern{-0.6}\kern\dimen@#1}%
    \fi
  }%
  \macc@depth\@ne
  \let\math@bgroup\@empty \let\math@egroup\macc@set@skewchar
  \mathsurround\z@ \frozen@everymath{\mathgroup\macc@group\relax}%
  \macc@set@skewchar\relax
  \let\mathaccentV\macc@nested@a
  \if#31
    \macc@nested@a\relax111{#1}%
  \else
    \def\gobble@till@marker##1\endmarker{}%
    \futurelet\first@char\gobble@till@marker#1\endmarker
    \ifcat\noexpand\first@char A\else
      \def\first@char{}%
    \fi
    \macc@nested@a\relax111{\first@char}%
  \fi
  \endgroup
}
\makeatother